\documentclass[11pt]{article}

\usepackage{amsmath, amsfonts, amssymb,latexsym,wasysym,graphicx,stmaryrd,tikz}
\usepackage{pgfplots}
\usepackage{float}
\input epsf.sty

\usetikzlibrary{decorations.pathmorphing}
\usetikzlibrary{arrows.meta,arrows}
\usetikzlibrary{shapes,shapes.geometric}

\newtheorem{Theorem}{Theorem}[section]
\newtheorem{Lemma}{Lemma}[section]
\newtheorem{Proposition}[Lemma]{Proposition}

\newtheorem{Definition}[Lemma]{Definition}

\newcommand{\BEQ}{\begin{equation}}     % Gleichungen Anfang ..
\newcommand{\BEA}{\begin{eqnarray}}
\newcommand{\BD}{\begin{displaymath}}
\newcommand{\EEQ}{\end{equation}}       % .. und Ende
\newcommand{\EEA}{\end{eqnarray}}
\newcommand{\ED}{\end{displaymath}}
\newcommand{\del}{\delta}
\newcommand{\Del}{\Delta}
\newcommand{\eps}{\varepsilon}          % epsilon

\newcommand{\supp}{{\mathrm{supp}}}

\newcommand{\LF}{{\mathrm{LF}}}

\newcommand{\R}{\mathbb{R}}

\newcommand{\Z}{\mathbb{Z}}
\newcommand{\N}{\mathbb{N}}
\newcommand{\D}{\mathbb{D}}
\renewcommand{\P}{\mathbb{P}}

\newcommand{\Id}{{\mathrm{Id}}}

\newcommand{\KPZ}{{\mathrm{KPZ}}}

\def\proba{{\mathbb{P}}}

\def\esper{{\mathbb{E}}}
\def\T{{\mathbb{T}}}
\def\F{{\mathbb{F}}}
\def\G{{\mathbb{G}}}

\def\LF{{\mathrm{LF}}}

\newcommand{\eop}{\hfill $\Box$}        % quod erat demonstrandum ...

\newcommand{\II}{{\rm i}}               % gerades i fuer komplexe Einheit
\renewcommand{\Re}{{\rm Re\ }}          % Realteil
\renewcommand{\Im}{{\rm Im\ }}          % Imaginaerteil
\newcommand{\half}{{1\over 2}}          % 1/2 als Bruch
\renewcommand{\vec}[1]{\boldsymbol{#1}} % Vektoren fettgedruckt
\catcode`\@=11
\def\numberbysection{\@addtoreset{equation}{section}
        \def\theequation{\thesection.\arabic{equation}}}
\numberbysection

\begin{document}

\title
{\bf The scaling limit of the KPZ equation in space dimension 3 and 
higher }

\date {}
\vskip -2cm
\maketitle

\vskip -2cm
\begin{center}
{\bf  Jacques Magnen$^a$ and J\'er\'emie Unterberger$^b$}
\end{center}

%\maketitle

%\vskip 0.5 cm
\centerline {\small $^a$Centre de Physique Th\'eorique,\footnote{Laboratoire associ\'e au CNRS UMR 7644}
Ecole Polytechnique,} 
\centerline{91128 Palaiseau Cedex, France}
\centerline{jacques.magnen@cpht.polytechnique.fr}
\vskip 0.5 cm
\centerline {\small $^b$Institut Elie Cartan,\footnote{Laboratoire 
associ\'e au CNRS UMR 7502. J. Unterberger acknowledges the support of the ANR, via the ANR project ANR-16-CE40-0020-01.} Universit\'e de Lorraine,} 
\centerline{ B.P. 239, 
F -- 54506 Vand{\oe}uvre-l\`es-Nancy Cedex, France}
\centerline{jeremie.unterberger@univ-lorraine.fr}

\vspace{2mm}
\begin{quote}

\renewcommand{\baselinestretch}{1.0}
\footnotesize
{We study in the present article the Kardar-Parisi-Zhang  (KPZ) equation
$$ \partial_t h(t,x)=\nu\Del h(t,x)+\lambda  |\nabla h(t,x)|^2 +\sqrt{D}\, \eta(t,x), \qquad (t,x)\in\R_+\times\R^d $$
in $d\ge 3$ dimensions in the perturbative regime, i.e. for $\lambda>0$ small enough and a  smooth, bounded, integrable  initial condition
$h_0=h(t=0,\cdot)$. The forcing term $\eta$ in the right-hand side is
 a regularized space-time white noise. The exponential of $h$ -- its so-called Cole-Hopf
 transform --  is known to satisfy a
 linear PDE with multiplicative noise. 
We  prove a large-scale diffusive limit for the solution, in particular a 
time-integrated heat-kernel behavior for the covariance in a parabolic scaling.

The proof is  based on  a rigorous implementation of K. Wilson's renormalization group
scheme. A double cluster/momentum-decoupling expansion  allows for  perturbative
estimates of the bare resolvent of the Cole-Hopf linear PDE in the small-field region where the noise is not too large, following the broad lines of Iagolnitzer-Magnen \cite{IagMag}. Standard large deviation estimates for $\eta$ make it possible to extend the above estimates to the large-field region. Finally, 
we show, by resumming all the by-products of the expansion, that the solution $h$  may be written in the large-scale limit (after a suitable Galilei transformation) as a small perturbation of the  solution  of the underlying linear Edwards-Wilkinson model ($\lambda=0$) with renormalized coefficients $\nu_{eff}=\nu+O(\lambda^2),D_{eff}=D+O(\lambda^2)$.

}
\end{quote}

\vspace{4mm}
\noindent



\newpage
\tableofcontents

%\medskip

\section{Introduction}

The KPZ equation \cite{KPZ} is a stochastic partial differential equation describing the growth by normal deposition of an interface  in $(d+1)$ space dimensions, see e.g. \cite{Bar,Car}. 
By definition the time evolution of the height $h(t,x)$, $x\in\R^d$, is given by
\BEQ \partial_t h(t,x)=\nu\Del h(t,x)+\lambda  |\nabla h|^2 +\sqrt{D}\, \eta(t,x), \qquad x\in\R^d \label{eq:KPZ} \EEQ
 where $\eta(t,x)$ is a regularized white noise, and $\nu,\lambda,D>0$ are constant. Three terms contribute to eq. (\ref{eq:KPZ}): a viscous term proportional to the {\em viscosity} $\nu$, leading to a smoothening
 of the interface; a growth by normal deposition with rate $\lambda$,
 called {\em deposition rate}, and playing the r\^ole of a {\em coupling constant}; and a random rise or lowering of the
 interface modelling molecular diffusivity, with coefficient $D$
 called  {\em noise strength}. In a related context, $h$ also represents  the free energy of  directed polymers in a random environment \cite{ImbSpe,CarHu,Cor}.
   It makes sense to consider more
 general nonlinearities of the form $V(\nabla h)$ with $V$, say, positive and
 convex, instead of $|\nabla h|^2$, which is in any case an approximation of 
 $2(\sqrt{1+|\nabla h(t,x)|^2}-1)$, assuming that   
 the gradient $|\nabla h|$ (the slope of the interface)  remains throughout small enough so that the evolution makes physically sense, precluding
e.g overhangs.

\medskip\noindent The interest is here  in the large-scale limit of this equation, for $t$ and/or $x$ large. 
A well-known naive rescaling argument gives some ideas about the dependence on the dimension of this limit. 
 Namely, the linearized equation, a stochastic heat (or infinite-dimensional Ornstein-Uhlenbeck \cite{Oks}) equation called 
 {\em Edwards-Wilkinson model} \cite{Bar} in the physics literature,
\BEQ \partial_t \phi(t,x)=\nu \Del \phi(t,x)+ \sqrt{D}\,  \eta(t,x), \quad (t,x)\in\R_+\times\R^d  \label{eq:EW} 
\EEQ
-- where $\eta$ requires no regularization -- 
 is invariant under the  rescaling $\phi(t,x)\mapsto \phi^{\eps}(t,x):= \eps^{-\half(\frac{d}{2}-1)}\phi(\eps^{-1} t,\eps^{-\half} x)$; we used here the equality in distribution,
$\eta(\eps^{-1} t,\eps^{-\half}x)\stackrel{(d)}{=} \eps^{\half(1+\frac{d}{2})}\eta(t,x)$. 
Assuming that $\phi$ is a solution of the KPZ equation instead yields after rescaling 
\BEQ \partial_t \phi^{\eps}(t,x)=\nu \Del \phi^{\eps}(t,x)+ \eps^{\half(\frac{d}{2}-1)} \frac{\lambda}{2}|\nabla \phi^{\eps}(t,x)|^2 + \sqrt{D}\, \eta^{\eps}(t,x), \label{eq:rescaled-EW} \EEQ
where (up to change of regularization) $\eta^{\eps}\stackrel{(d)}{=}\eta$.
 For $d>2$, $\eps^{\half(\frac{d}{2}-1)}$ vanishes in the limit $\eps\to 0$; in other terms, 
the KPZ equation is {\em infra-red super-renormalizable}, hence {\em (power-like) asymptotically free at large scales}  in $\ge 3$ dimensions, i.e.  expected to behave, in a {\em small coupling} (also called {\em small disorder}) regime where $\lambda\ll 1$, like the corresponding linearized equation up to a redefinition
(called {\em renormalization}) of the 
{\em diffusion constant} $\nu$ and of the {\em noise strength} $D$. 

\medskip
 Let us emphasize the striking difference  with the one-dimensional $\KPZ_1$ equation.
For this equation, scaling behaviors, see (\ref{eq:rescaled-EW}), are reversed with
respect to $d\ge 3$, in other words,
KPZ$_1$ is (power-like) asymptotically free at small scales (i.e. in the ultra-violet), or
equivalently (in the PDE analysts' terminology) {\em sub-critical}. A large part of the interest for this
equation  comes from the fact that the large-scale strongly coupled theory 
\cite{ACQ,Cor} is
understood by comparison with  integrable discrete statistical physics models \cite{DEHP,PraSpo,Sasspo,TraWid} relating to weakly asymmetric
exclusion process \cite{BerGia} or the Tracy-Widom distribution of the largest eigenvalue of random matrices connected with Bethe Ansatz \cite{TraWid},  free fermions and
 determinantal processes \cite{Joh},... Note that $\KPZ_2$ is believed by
perturbative QFT arguments to be strongly coupled at large scales \cite{Bar,Car} and its large-scale limit is not at all understood. 

\bigskip\noindent 
{\em We prove the diffusive limit of  $d$-dimensional KPZ $(d\ge 3)$ with small coupling  in the present work}, 
thus establishing on firm mathematical ground old predictions of
physicists, see e.g. Cardy \cite{Car}. The 
space dimension $d$ does not really matter as long as $d\ge 3$. In
  the small-coupling regime, contrary to the $1d$-case, we fall into the Edwards-Wilkinson universality class.

\bigskip
**************************************************************************

\medskip\noindent In comparison with the achievements made in the study of strongly
coupled large-scale $\KPZ_1$, this problem looks at first sight of lesser importance and difficulty. We believe that the interest of our result lies in the {\em precision
of our asymptotics}, and in the potential wide scope of applicability of our methods.

\medskip\noindent Namely, the KPZ model is one particular instance of a large
variety of   dynamical problems in statistical physics, modelized as  interacting particle systems,  or
as  parabolic SPDEs heuristically derived by some mesoscopic limit, which have been  turned into a functional integral form analogous to the Gibbs measure of 
{\em equilibrium} statistical mechanics, $e^{-\int {\cal L}_0- g\int {\cal L}_{int}}$,  using the so-called {\em response field} (RF), or {\em Martin-Siggia-Rose} (MSR) 
{\em formalism} and
studied by using standard perturbative expansions originated from quantum field theory (QFT); for reviews see e.g. \cite{Car} or \cite{AltBen}. Despite the lack of mathematical rigor,
this formalism yields a correct description of  the qualitative behaviour of such dynamical problems in the large scale limit.

 \medskip\noindent The {\em Feynman  perturbative} approach, see   e.g.
 \cite{LeB}, consists
in expanding $\exp -g\int {\cal L}_{int}$ into a series in $g$ and making a clever resummation of some truncation of it into so-called {\em counterterms}, represented 
in terms of a sum of diagrams;  as such, it is non-rigorous, since it yields $N$-point functions in terms
of an asymptotic expansion in the coupling parameter $g$ which is divergent in all interesting cases (at least for bosonic theories).  A few years ago, however, Gubinelli, M. Hairer,
H. Weber,... \cite{Gub,Hai,Hai2,BaiBer,BruHaiZam, CanFriGas,ChaHai,ChoAll,HaiLab,MouWeb,CatCho}, drawing sometimes on a dynamical approach
to the construction of equilibrium measures advocated by Nelson \cite{Nel2}, 
Parisi-Wu \cite{ParWu}, and Jona-Lasinio, Mitter and S\'en\'eor \cite{Jon1,Jon2,Jon3},
  have started developing this philosophy in a systematic way to
solve sub-critical parabolic SPDEs rigorously, i.e. beyond perturbation theory. Such SPDEs have only a finite number of counterterms,
each counterterm being the sum of a finite number of terms (that can be
interpreted in terms of Feynman diagrams), which makes the task considerably easier,
but still far from trivial.

\bigskip
**************************************************************************
\medskip

\medskip\noindent {\em Constructive approaches}  developed in the context of
statistical physics  by mathematical physicists from the mid-60es, see e.g.  \cite{Erice,GalNic,GawKup,FMRS,FS,GJa,IagMag,MagUnt1,MagUnt2,Nel} and
surveys \cite{GJb,Mas,Riv,Sal,Unt-mode},
have developed sophisticated, systematic truncation methods making it possible to control the error terms. The partial resummations are interpreted in the manner of K. Wilson \cite{Wilb,Wilc} as a
{\em scale-by-scale, finite renormalization} of the parameters $\nu,\Del,\lambda$ of the Lagrangian ${\cal L}_0+g{\cal L}_{int}$.  In many instances it
has proved possible to subtract {\em scale counterterms} explicitly by hand and prove that
the remainder is finite, yielding some description of the effective, large-scale
theory, see e.g. works in diverse contexts -- random walks in random environment,
 KAM theory, etc. -- by Bricmont, Gawedzki, Kupiainen and 
coauthors \cite{BricKup,BricGawKup,BKL}, and recent extensions to the study
of sub-critical parabolic PDEs \cite{Kup,KupMar}, as an alternative to the "global
counterterm" strategy mentioned in the last paragraph. However, the implementation
of a full-fledged, multi-scale constructive scheme is  for the moment 
 limited to equilibrium statistical
physics models.

\medskip\noindent The present work is, to the best of our knowledge, the first attempt
to use such a scheme in the context of non-equilibrium statistical mechanics, here
for a parabolic SPDE. Instead of using the MSR formalism,  we develop (as all  previously mentioned mathematically rigorous approaches do) a more straightforward approach, starting directly from the equation and cutting the propagator $e^{t\nu\Del}$ into scales. We actually work on the following model.

\bigskip
**************************************************************************

\medskip\noindent
{\bf The model.} {\em Let $d\ge 3$.
We consider the following equation on $\R_+\times \R^d$,
\BEQ (\partial_t-\nu^{(0)}\Del)h(t,x)=\lambda |\nabla h(t,x)|^2+ \sqrt{D^{(0)}}\, 
(\eta(t,x)-v^{(0)}), \qquad h\big|_{t=0}=h_0 \label{eq:model} \EEQ
where $\eta$ is a  white noise regularized in time and in space; 
$h_0$ is a  smooth, bounded, integrable initial condition, i.e.
$||h_0||_{L^{\infty}}:=\sup_{x\in\R^d} |h_0(x)|, ||h_0||_{L^1}:=\int_{\R^d} dx\,
 |h_0(x)|$ are $<\infty$;  $\lambda>0$ is small enough; and $v^{(0)}$ is a 
constant, average interface velocity which we shall fix later on.}

\medskip\noindent  The precise choice of regularization for the white noise is unimportant; one should just keep in mind that local (in
time and space)
solvability of (\ref{eq:KPZ}) in a strong sense requires that, for every compact set $\bar{\Del}\subset
\R^d$ (equivalently, for any $\bar{\Del}\in\bar{\D}^0$ as in Definition 
\ref{def:phase-space} (iii)),  $t\mapsto \sup_{x\in\bar{\Del}} \left(|\eta(t,x)|
+|\nabla\eta(t,x)| \right)$ is locally  integrable.
For simplicity of exposition,  we define $\eta$ to be a smooth, stationary Gaussian noise with short-range covariance.  To be definite:

\medskip\noindent {\em we fix a smooth, isotropic (i.e. invariant under
space rotations) function $\omega:\R\times\R^d\to\R$ with support $\subset B(0,\half)$ and $L^1$-norm $\int dt\, dx\, \omega(t,x)=1$, and let 
 \BEQ \langle \eta(t,x)\eta(t',x')\rangle:=(\omega\ast\omega)(t-t',x-x')=
 \int dt''\, \int dx''\, \omega(t-t'',x-x'')\omega(t''-t',x''-x').
 \label{eq:cov-eta} \EEQ }

\medskip \noindent
Our main result is the following. Gaussian expectation with respect to $\eta$ is denoted either  by $\langle \ \cdot\ \rangle$, or $\langle \ \cdot \rangle_{\lambda}$ or also $\langle \ \cdot \rangle_{\lambda;\nu^{(0)},D^{(0)},v^{(0)}}$
 if one wants to emphasize the dependence
on the parameters $\nu^{(0)},D^{(0)},\lambda,v^{(0)}$; the result also depends obviously on the initial condition $h_0$. {\em By convention, $\langle \cdot \rangle_{0;\nu,D}$ refers to the expectation with respect to the  measure of the
Edwards-Wilkinson equation $(\partial_t-\nu\Del)\phi(t,x)=\sqrt{D}\, \eta(t,x)$ with
zero initial condition, where
$\eta$ is a standard (unregularized) space-time white noise; for this equation we implicitly set $v=0$.} By definition, $\phi(t,x)=\sqrt{D} \int_0^t ds\, \left(e^{(t-s)\nu\Del} \eta_s\right)(x)$ is a centered Gaussian process.

\begin{Theorem} \label{th:main} {\bf (Main Theorem).}\\
Let $d\ge 3$. Fix $D^{(0)},\nu^{(0)}>0$ and a smooth, bounded, integrable initial condition $h_0$. Let  $\lambda\ge 0$ be small enough, $\lambda\le \lambda_{max}=\lambda_{max}(||h_0||_{L^1}, ||h_0||_{L^{\infty}})$. Then  there exist three coefficients $D_{eff}=
D^{(0)}+O(\lambda^2)$, $\nu_{eff}=\nu^{(0)}+O(\lambda^2)$ and $v^{(0)}=v^{(0)}(\lambda)=O(\lambda)$, all independent of
the initial condition $h_0$,
 such that the solution $h$ of the KPZ equation (\ref{eq:model}) satisfies the following asymptotic properties:
 
\begin{enumerate}
\item for all $(t,x)$ with $t>0$,
\BEQ\langle h_{\eps^{-1}t}(\eps^{-\half}x)\rangle_{\lambda;\nu^{(0)},D^{(0)},v^{(0)}}= O_{\eps\to 0}(\eps^{d/2}); \label{eq:0.5} \EEQ

\item  for all $(t_1,x_1),\ldots,(t_{2N},x_{2N})$, $N\ge 1$ with $t_i>0$, $i=1,\ldots,2N$ and $(t_i,x_i)\not=(t_j,x_j),  i\not=j$, letting $h_i:=\langle h_{\eps^{-1}t_i}(\eps^{-\half}x_i)\rangle_{\lambda;v^{(0)},\nu^{(0)},D^{(0)}}$,
\BEQ \Big\langle \prod_{i=1}^{2N} \left( h_{\eps^{-1} t_i}(\eps^{-\half} x_i) -h_i\right)\Big\rangle_{\lambda;v^{(0)},\nu^{(0)},D^{(0)}} \sim_{\eps\to 0}
\eps^{N(\frac{d}{2}-1)}
\Big\langle \prod_{i=1}^{2N} h_{ t_i}(x_i)\Big\rangle_{0;\nu_{eff},D_{eff}}. \label{eq:0.6}\EEQ

\end{enumerate}

\end{Theorem}

\noindent Since $\langle\ \cdot\ \rangle_{0;\nu_{eff},D_{eff}}$ is a Gaussian measure, {\em 2.} may
be rephrased as follows. Let 
\BEQ 
K_{eff}(t_1,x_1;t_2,x_2):=\lim_{\eps\to 0}
\eps^{-(\frac{d}{2}-1)}
\big\langle\left( h_{\eps^{-1}t_1}(\eps^{-1/2}x_1)-h_1\right) \, 
\left( h_{\eps^{-1}t_2}(\eps^{-1/2}x_2)-h_2\right) \big\rangle_{\lambda;
 v^{(0)},\nu^{(0)},D^{(0)}} \label{eq:intro-Keff}
 \EEQ
  ($t,t'>0$, $(t,x)\not=(t',x')$). Then
\BEQ K_{eff}(t,x;t',x')=  \langle h(t,x)h(t',x')\rangle_{0;\nu_{eff},D_{eff}}   
\EEQ

and
  
 \BEQ \Big\langle \prod_{i=1}^{2N} \left( h_{\eps^{-1} t_i}(\eps^{-\half} x_i) -h_i\right)\Big\rangle_{\lambda;v^{(0)},\nu^{(0)},D^{(0)}} \sim_{\eps\to 0}
\eps^{N(\frac{d}{2}-1)} \sum_{{\mathrm{pairings}}} \prod_{j=1}^N K_{eff}(t_{i_{2j-1}},x_{i_{2j-1}};t_{i_{2j}},
x_{i_{2j}}) \label{eq:0.7}
\EEQ
where the sum ranges over all pairings $(i_1,i_2),\ldots,(i_{2N-1},i_{2N})$ of
the $2N$ indices $1,2,\ldots,2N$.

\medskip

In other words, {\em up to a Galilei transformation $h_t(x)\mapsto h_t(x)-
\sqrt{D^{(0)}}\, v^{(0)}t$, the  $N$-point functions of
the KPZ equation $(\partial_t-\nu^{(0)}\Del)h=\lambda|\nabla h|^2+\sqrt{D^{(0)}}
\, \eta$  behave asymptotically in the large-scale limit as the $N$-point functions of the 
solution of the  Edwards-Wilkinson equation with renormalized coefficients $D_{eff},\nu_{eff}$},
\BEQ (\partial_t -\nu_{eff} \Del) h_t(x)=\sqrt{D_{eff}}\, \eta(t,x) \qquad (t\ge 0), \qquad 
h_0\equiv 0 \EEQ
where $\eta$ requires no regularization. Generally speaking, main corrections to
the above asymptotic behaviour (\ref{eq:0.5},\ref{eq:0.7}) are smaller by $O(\eps^{(1/2)^-})$ as proved in \S 5.3 {\bf D.}
Effective coefficients $D_{eff},\nu_{eff}$ have a
(diverging) asymptotic expansion in terms of $\lambda$;  lowest-order corrections in $O(\lambda^2)$ are computed in (\ref{eq:leading-nu-eff}) and (\ref{eq:leading-D-eff}).
The $O(\eps^{d/2})$-term in (\ref{eq:0.5}) is a contribution due to the
initial condition; further contributions of the initial condition to  $N$-point functions
come with an extra multiplicative factor in $O(\lambda\eps^{\frac{d}{2}-1})$, which is the
scaling of the vertex. Corrections to Gaussianity of $N$-point
functions, of order $O(\lambda^2 \eps^{\frac{d}{2}-1})$, are examined in (2) a few pages below. Furthermore, our multi-scale scheme actually involves an effective propagator
differing slightly from the effective Edwards-Wilkinson propagator $e^{(t-s)\nu_{eff}\Del}$, see section 7; this implies 
a correction w.r. to the r.h.s. of (\ref{eq:0.7}) with a small extra prefactor,
which is proved to be a  $O(\eps)$ but could easily be improved to $O(\eps^n)$ with
$n$ arbitrary large. 

\medskip\noindent {\em Remark.} A more common choice of regularization for $\eta$ is to take a discretized "kick force", namely, we 
pave $\R_+$ by unit size intervals $[n,n+1)$, $n\ge 0$,   and let
$\xi_{n+\half}:=\eta\big|_{[n,n+1)}$, $n=0,1,\ldots$ be independent, centered Gaussian fields on $\R^d$
which are constant in time and have smooth, space-translation invariant covariance kernel
with finite range, for instance. This does not change the conclusion of Theorem 0.1,
except that, the law of $\eta$ being now only $\Z$-periodic in time, $h_{\infty}(t):=\lim_{n\to +\infty} \langle h_{n+t}(0)\rangle$ is
 now a $1$-periodic function instead of the constant $0$. This
  regularization has several advantages (see section \ref{section:Cole-Hopf}); it allows in particular an explicit representation of $v^{(0)}$ in
probabilistic terms. The scheme of proof
 extends without any significant modification if the covariance kernel decreases heat-kernel-like in space, e.g. if $\xi_{n+\half}\overset{(d)}{=}e^{c\Del}
\xi$ where $\xi$ is a  standard space white noise,
and $c>0$ is some constant.

\bigskip
\noindent Furthermore, it follows from the proof (see section 5) that the value of $v^{(0)}$ may be obtained by equating it to the 
constant $\tilde{v}^{(0)}$ such that $\langle w(t,0)\rangle_{\tilde{v}^{(0)}}=O(1)$ independently of $t$, in coherence 
with the value obtained in Carmona-Hu \cite{CarHu} in a discrete setting for a random directed polymer
measure  (see
section 2.1), where $w$ is the Cole-Hopf transform of $h$ (see below). Let us note that the equality between $v^{(0)}$ and
$\tilde{v}^{(0)}$ points out to the fact that we are in a weak disorder regime
in which the annealed and quenched free energies coincide. However, our proof
is independent of that of Carmona and Hu (see \cite{CarHu}, Theorem 1.5), based
on Gaussian concentration inequalities.

\bigskip
**************************************************************************

\medskip\noindent
The  {\em proof} follows closely the article by  Iagolnitzer-Magnen \cite{IagMag}
on weakly self-avoiding polymers in four dimensions, which is the main reference
for the present work. Namely, up to the change of function
$h\mapsto w:=e^{\frac{\lambda}{\nu}h}$ (called {\em Cole-Hopf transform}) and of coupling constant, $g:=\frac{\lambda}{\nu}\sqrt{D}$, the KPZ equation
is equivalent to the linear equation  $(\partial_t-\nu\Del)w=g\eta w$, solved as
$w(t,x):=\int dy\, G_{\eta}((t,x),(0,y))w_0(y)$, where $G_{\eta}\equiv
\left(\partial_t-\nu\Del-g\eta\right)^{-1}$  is a {\em random resolvent}. Formally then,
our problem is a parabolic counterpart to the large-scale analysis of polymers
in a weak random potential solved in \cite{IagMag} by studying the equilibrium
resolvent $\left(\Del+\II g\eta\right)^{-1}$, where the "$\II$"-coefficient
is the Edwards model representation of the self-avoiding condition (the model
is solved for $g\ll 1$ but the self-avoiding condition is recovered for $g=1$). Though the
two models are physically unrelated, one must analyze similar mathematical objects. As is often the case, the model with a time evolution (i.e. the parabolic one) turns out
to be easier than the equilibrium model (i.e. the elliptic one), because of the
{\em causality constraint}.

\medskip\noindent  The {\em general scheme of proof}, following, as mentioned
above,  the philosophy of constructive
field theory,  is to introduce a multi-scale expansion and define a renormalization mapping, $\nu=\nu^{(0)}\longrightarrow \nu^{(1)}\longrightarrow \ldots \longrightarrow \nu^{(\infty)}:=\nu_{eff}$, 
$D=D^{(0)}\longrightarrow D^{(1)}\longrightarrow \ldots \longrightarrow D^{(\infty)}:=D_{eff}$ or equivalently $g^{(0)}:=\frac{\lambda}{\nu_0}
\sqrt{D^{(0)}}\to g^{(1)}\longrightarrow \ldots\longrightarrow g^{(\infty)}\equiv
g_{eff}=\frac{\lambda}{\nu_{eff}}
\sqrt{D_{eff}}$ (later on interpreted as the flow of the  {\em coupling constant}
through the Cole-Hopf transform), $v=v^{(0)}\longrightarrow v^{(1)}\longrightarrow\ldots\longrightarrow v^{(\infty)}\equiv v_{eff}:=0$
ensuring the convergence of the expansion at each scale and allowing to control error terms.
The average interface velocity $v^{(0)}$ is fixed by requiring that the asymptotic
velocity $v_{eff}$ vanishes. The original parameters $\nu^{(0)},D^{(0)},v^{(0)}$, called {\em bare parameters}, 
describe the theory at scale $O(1)$, while the Edwards-Wilkinson model with scale $j$
parameters $\nu^{(j)},D^{(j)}$ and drift velocity $v^{(j)}$ give a good approximation
of the theory at time distances of order $\eps^{-1}=2^{-j}$, which becomes asymptotically
exact in the infra-red limit, when $j\to\infty$.  This goal is achieved in
general by using  a  phase-space expansion, i.e.  a {\em horizontal cluster expansion}  
casting into the form of a series   the interactions at a given energy-momentum level between the degrees of freedom,
and a {\em vertical cluster} or {\em momentum-decoupling expansion} separating the
different energy-momentum levels. Energy, resp. momentum, are the Fourier
conjugate variables of time and space; here a given energy-momentum level $j$ is
adequately defined by considering heat-kernel {\em propagators}
$$ G_{\nu}((t,x),(t',x'))=e^{\nu(t-t')\Del}(x-x')=p_{\nu(t-t')}(x-x') $$
 with $t-t'\approx 2^j$. Then the above series (roughly speaking, a
truncated power series in the coupling constants with a bounded integral,
Taylor-like remainder)  converge if
the bare coupling constant $g^{(0)}$ is small enough.

\medskip\noindent With our choice of covariance
function for $\eta$,  however, the flow of the parameters $\nu,v$ is actually trivial starting from $j=1$, i.e. $\nu^{(j)}=\nu_{eff},\, v^{(j)}=v_{eff}=0$ for
$j\ge 1$, and
the noise strength $D$, defined by resumming connected diagrams  with four external legs, though scale-dependent, requires no renormalization at all, because the equation is 
{\em infra-red
super-renormalizable}, and the total correction (obtained by summing over scales) is {\em finite}. This, and also the causality condition preventing the
so-called {\em low-momentum field accumulation} problem \cite{GJa,FMRS,Unt-mode},  leads to a much simplified
framework, from which the phase space analysis has almost disappeared. Only scale 0,
 two-point diagrams need to be renormalized, with a contribution at
near zero momentum $\vec{k}$
 $$v^{(0)}+(\nu_{eff}-\nu^{(0)})|\vec{k}|^2\equiv v^{(0)}-(\nu_{eff}-\nu^{(0)})\Del,$$ leaving a remainder of parabolic order three in the momenta, i.e.  $O(\nabla^3)$ or  $O(\nabla\partial_t)$.  Scale 0 diagrams are connected by "low-momentum" heat-kernel
propagators $G((t,\cdot),(t',\cdot))$ with $t-t'\approx 2^j$, $j\ge 1$.  A crucial point
in the proof is that, thanks to the $\nabla^3$,  remainders  integrated over space-time
cost a factor $O(1)$, namely (see (\ref{eq:PW1}) and (\ref{eq:PW1bis}))
\BEA 
&& \int_{t''}^t dt'\, \int dx'\,  G((t,x),(t',x')) \, |\nabla^3 G((t',x'),(t'',x'')| 
\nonumber\\
&&\lesssim 
\left(\int dt'\, (1+|t'-t''|)^{-3/2}\right) \, p_{\nu(t-t'')}(c|x-x''|) = O(1)\  p_{\nu(t-t')}(c|x-x''|) \nonumber\\
\EEA
or, simply said, $G\nabla^3 G\lesssim G$.
 What is left of the cluster expansions is adequately resummed
as in \cite{IagMag} into the random resolvent in the form of localized
"vertex insertions" (see section 5), thereby suppressing combinatorial factors
which make the series divergent. Then the contribution of all  vertex
insertions is bounded by some contour integral of a modified
resolvent  through the use of
Cauchy's formula. 

\medskip\noindent An extra complication comes however from the inverse Cole-Hopf transform. Applying cluster expansions -- which is done in practice by
differentiation with respect to some additional parameters --  to $\log(w)$ leads to rational expressions of the form $\frac{"D_1 w" \cdots "D_n w"}{w^n}$, where
the $D_i$'s are differential operators, acting on "replicas" of $w$. Then the
scale 0 diagrams requiring renormalization can be factorized, hence averaged
with respect to the measure $\langle \, \cdot\, \rangle$. Remaining terms are shown
to yield a convergent series in the form of a sum over "polymers"  for $\lambda$ small enough.

\medskip
The $\lambda$ and $\eps$-pre-factors contained in Theorem \ref{th:main} may
be guessed from the following guiding principles, put into light by the cluster expansion. 
 
\medskip\noindent (1) First, the two-point function of the renormalized Edwards-Wilkinson equation, 
\BEA &&  \langle h(\eps^{-1} t,\eps^{-1/2} x)h(\eps^{-1} t',\eps^{-\half} x'
)\rangle_{0;\nu_{eff},D_{eff}}  \nonumber\\
&&\qquad = D_{eff}\int_{0}^{\eps^{-1} t'} ds\, \int dy\,
 p_{\nu_{eff}(\eps^{-1}t-s)}(\eps^{-\half}x-y) p_{\nu_{eff}(\eps^{-1}t'-s)}(y-\eps^{-\half}x') \nonumber\\
&&\qquad =  D_{eff} \int_{0}^{\eps^{-1} t'} ds\, p_{\nu_{eff}(\eps^{-1}(t+t')-2s)}(\eps^{-\half}(x-x'))  \label{eq:E-W-scaling}
\EEA
 ($t\ge t'>0$), scales like $\eps^{\frac{d}{2}-1}$, as can be seen
by simply rescaling variables $t'\to \eps t', (x,x')\to (\eps^{1/2}x,\eps^{1/2}x')$ in
the integral. There
are two regimes: the {\em equilibrium regime} ($t-t'\lesssim |x-x'|^2$), in
which
$\langle h(t,x)h(t',x')\rangle_{0;\nu_{eff},D_{eff}}\approx \int_{0}^t ds\, p_{2(t-s)}(x-x')\approx |x-x'|^{-(d-2)}$ is essentially the equilibrium Green function of the Laplacian;
the {\em dynamical regime}  ($t-t'\gtrsim |x-x'|^2$), in which
$\langle h(t,x)h(t',x')\rangle_{0;\nu_{eff},D_{eff}}\approx \int_{0}^{t'} ds\, p_{t+t'-2s}(0)\approx |t-t'|^{-(\frac{d}{2}-1)}$. 

\medskip\noindent (2) The {\em connected} quantities $\Big\langle \prod_{i=1}^{2N}  h_{\eps^{-1} t_i}(\eps^{-\half} x_i) \Big\rangle^{{\mathrm{connected}}}_{\lambda;v^{(0)},\nu^{(0)},D^{(0)}}$  (also
called {\em truncated $2N$-point functions}) are
$O\left(\left[\lambda^{2} \eps^{\frac{d}{2}-1}\right]^{2N-1}\right)$. Namely, Gaussian pairwise
contractions yield the expected
scaling in $O(\eps^{N(\frac{d}{2}-1)})$, i.e. $O(\eps^{\frac{d}{2}-1})$ per
link, as expected from (1); whereas the connected expectation requires 
$N-1$ supplementary links and twice as much vertices (since these are not
present in the linear theory) in the expansion, contributing an extra small
$O\left( \left[(\lambda^2 \eps^{\frac{d}{2}-1})\right]^{N-1}\right)$  
prefactor. The cluster expansion makes it possible to develop those links
explicitly.

%%%%%%%%%%%%%%%%%%%%%%%%%%%%%%%%%%%%%%%%%

\bigskip\noindent The plan of the article is as follows. We start by recalling the
Cole-Hopf transform  in section  1, and make the bridge to previous results on the
subject stated in terms of the associated directed polymer measure. We then introduce
in section 2 
a multi-scale expansion for the propagators, together with multi-scale estimates (also
called "power-counting"), which are the building blocks of our approach.  Sections 3, 4, and 5 are the heart of the article.
The dressed equation, and the cluster expansion thereof, is presented in section 3.
Section 4 is dedicated to renormalization; the scale 0 counterterms obtained by
factorizing two-point functions through
a supplementary Mayer expansion  are bounded. Then we show in section 5 how to bound the
sum of all terms produced by the expansion, and obtain final bounds for $N$-point
functions, proving thus our main result, Theorem 0.1. Finally, there are two appendices. In the first one, we provide detailed combinatorial formulas for the horizontal and Mayer  cluster
expansions.  The second one is merely dedicated to a technical result. Pictures are provided, which are there to help the reader
visualize the outcome of the various expansions.

%%%%%%%%%%%%%%%%%%%%%%%%%

\bigskip\noindent
{\em Notations.}
\begin{enumerate}
\item {\em (parabolic distance)} Let $d((t,x),(t',x')):=\sqrt{|t-t'|+|x-x'|^2}$
$(t,t'\in\R_+, \, x,x'\in\R^d)$. Similarly, for $U,U'\subset\R_+\times\R^d$, 
$d((t,x),U):=\inf_{(t',x')\in U} d((t,x),(t',x'))$, $d(U,U'):=\max\left(\sup_{(t,x)\in U}
 d((t,x),U'), \sup_{(t',x')\in U'} d(U,(t',x')) \right)$ (Hausdorff distance). Then $\bar{d}$ is the space projection of
the distance $d$, i.e. $\bar{d}(x,x'):=d((0,x),(0,x'))=|x-x'|$, etc.
\item Let $f,g:E\to\R$ be two functions on some set $E$. We write $|f(z)|\lesssim |g(z)|$  if there exists some 
inessential constant $C$ (possibly depending  on the  parameters
$D,\nu$ and on the space dimension $d$), uniform in $\lambda$
for $\lambda$ small enough, such that
$|f(z)|\le C|g(z)|$. Then, by definition, $|g(z)|\gtrsim |f(z)|$. If $|f(z)|\lesssim
|g(z)|$ and $|g(z)|\lesssim |f(z)|$, we write $|f(z)|\approx |g(z)|$.
\item In many situations, one obtains $(t,x)$-dependent functions $f(t,x)$ such that
$f$ decays Gaussian-like, $f(t,x)\le e^{-c|x|^2/t}$ for some positive constant $c$
bounded away from $0$. We then write $f(t,x)\le e^{-c|x|^2/t}$ without further
specifying the value of $c$, which may change from line to line. For instance, if 
$p_{\nu t}(x)=e^{\nu t\Del}(x)$ is the heat kernel, then we may write
  $p_{\nu t}(x)\lesssim t^{-d/2} e^{-c|x|^2/\nu t}\lesssim t^{-d/2} e^{-c'|x|^2/t}$,
leaving out the dependence in the parameter $\nu$ as explained in 2.  Note however that,
if $\nu'\le \nu$, $p_{\nu' t}(x)\lesssim p_{\nu t}(x)$, whereas the inequality
$p_{\nu t}(x)\lesssim p_{\nu' t}(x)$ does not hold uniformly in $x$ because the space
decay of $p_{\nu t}(\cdot)$ is slower than that of $p_{\nu' t}(\cdot)$.

\end{enumerate}

\medskip\noindent{\bf Acknowledgements.} We wish to thank H. Spohn, F. Toninelli and
the referee for numerous discussions, suggestions and corrections, which have hopefully contributed in particular to the readability of the paper.

%%%%%%%%%%%%%%%%%%%%%%%%%
%%%%%%%%%%%%%%%%%%%%%%%%

\section{Cole-Hopf transform}  \label{section:Cole-Hopf}

%%%%%%%%%%%%%%%%%%%%%%%%%%%%%%%ש
%%%%%%%%%%%%%%%%%%%%%%%%%%%%%%%

\medskip\noindent It is well-known that $w:=e^{\frac{\lambda}{\nu^{(0)}} h}$ is a solution of
the linear equation with multiplicative noise,

\BEQ (\partial_t-\nu^{(0)}\Del)w(t,a)=g^{(0)} \left(\eta(t,a) - v^{(0)} \right) w(t,a) \label{eq:w}
\EEQ
where
\BEQ g^{(0)}:=  \frac{\lambda}{\nu^{(0)}}\sqrt{D^{(0)}}=O(\lambda) \label{eq:g0}
\EEQ 
plays the r\^ole of a {\em bare coupling constant},
  from which (representing the solution as a Wiener integral by
Feynman-Kac's formula)
\BEQ h(T,a)=\frac{\nu^{(0)}}{\lambda} \log w(T,a), \qquad w(T,a)=   \esper^a\left[ e^{g^{(0)}\int_0^T dt\, 
\left( \eta(T-t,B_{t}) - v^{(0)} \right) } e^{\frac{\lambda}{\nu^{(0)}} h_0(B_T)} \right],  \label{eq:Wiener-integral}
\EEQ
 where the expectation $\esper^a$ is relative to the Wiener measure 
on  $d$-dimensional Brownian paths $(B_t)_{0\le t\le T}$ issued
from $a\in\R^d$ with $\nu^{(0)}$-normalization, i.e. $\esper^a[ (B_t^i-a)^2]=2\nu^{(0)}t$, $i=1,\ldots,d$. Thus $w(T,a)$ may be interpreted as the
partition function of a directed polymer, see e.g. \cite{CarHu} and references within, but we shall not need this interpretation in the article. Note that $(B_t)_{t\ge 0} \overset{(d)}{=} (W_{2\nu^{(0)} t})_{t\ge 0}$,
where $W$ is now a standard Brownian motion, from which -- forgetting about the
regularization and using the variable $2\nu^{(0)} t$ instead of $t$ --  
$$\int_0^T dt\, \eta(T-t,B_t)\sim \frac{1}{2\nu^{(0)}} \int_0^{2\nu^{(0)} T} du\, \eta(\frac{u}{2\nu^{(0)}},W_u) \overset{(d)}{=} \frac{1}{\sqrt{2\nu^{(0)}}} \int_0^{2\nu^{(0)} T} du\,  \eta(u,W_u).$$
 Thus $w(T,a)$
may be expanded in a series in the parameter $g:=\frac{g^{(0)}}{\sqrt{2\nu^{(0)}}}=
\frac{1}{\sqrt{2}}\frac{\lambda}{(\nu^{(0)})^{3/2}} \sqrt{D^{(0)}}$.
\medskip\noindent  Similarly, $\nabla w=\frac{\lambda}{\nu^{(0)}} e^{\frac{\lambda}{\nu^{(0)}} h}\nabla h$, or conversely
$\nabla h=\frac{\nu^{(0)}}{\lambda} \frac{\nabla w}{w}$, from which

\BEA && \nabla h(T,a)= e^{-\frac{\lambda}{\nu^{(0)}} h(T,a)}  \left(   \esper^a\left[ e^{g^{(0)}\int_0^T dt\, 
\left(  \eta(T-t,B_{t}) - v^{(0)} \right) }\  e^{\frac{\lambda}{\nu^{(0)}} h_0(B_T)} \nabla h_0(B_T) \right]  \right. \nonumber\\
&& \left. \qquad + \sqrt{D^{(0)}}\  \esper^a\left[  \int_0^T dt\,  e^{g^{(0)}\int_0^t ds\, 
\left(  \eta(T-s,B_{s}) - v^{(0)} \right) } \ 
\, \cdot\,  \nabla \eta(T-t,B_{t}) \, \cdot\,  e^{\frac{\lambda}{\nu^{(0)}} h_{T-t}(B_{t})} \right] \right) \nonumber\\
  \label{eq:4.3}
\EEA

\medskip\noindent  Without using the general theory developed in \cite{Unt-KPZ1,Unt-KPZ2}, eq. (\ref{eq:Wiener-integral}) and (\ref{eq:4.3}) show that a.s. $h$,$\nabla h$ exist and are $C^1$  for $h_0$, say, $C^1$ and compactly supported. The Cole-Hopf solution coincides with
the solution defined for more general Hamilton-Jacobi equations in \cite{Unt-KPZ1,Unt-KPZ2}.

\bigskip\noindent
{\em For the rest of the subsection only, we assume that $\eta$ is a discretized "kick force", i.e. $\eta\big|_{[n,n+1)}=:\xi_{n+\half}$ are independent and constant in time, in order to compare with the existing literature.} Since $(\eta|_{[n-1,n)})_{n\ge 0}$ are independent fields,
letting $v^{(0)}:=\tilde{v}^{(0)}$, where
\BEQ \tilde{v}^{(0)}:=\frac{1}{g^{(0)}} \log\,  \langle \esper^0\left[ e^{g^{(0)} \int_0^1
dt\, \eta(0,B_t)}\right]\rangle  \label{eq:Cole-Hopf-v0} \EEQ
leads to $\langle w(n,a)\rangle_{\tilde{v}^{(0)}}=1$ for any $n\in\N$ and $a\in\R^d$ if $w_0=1$, whence more generally
\BEQ \langle w(n,a)\rangle_{\tilde{v}^{(0)}}=O(1).\EEQ
Expanding the exponential  in (\ref{eq:Cole-Hopf-v0}) and using
\BEQ \Big\langle\Big|\int_0^1 dt\, \eta(0,B_t)\Big|^p \Big\rangle
\le \int_0^1 dt \, \langle |\eta^p(0,B_t)| \rangle \lesssim C^p \Gamma(p/2) \langle\eta^2(0,B_t) \rangle^{p/2}=O( (C')^p \Gamma(p/2)),
\EEQ  
one gets:  $\langle \esper^0\left[ e^{g^{(0)} \int_0^1
dt\, \eta(0,B_t)} \right]\rangle=e^{O(\lambda^2)}$, whence $\tilde{v}^{(0)}=O(\lambda)$.

\bigskip\noindent
Let us state an easy preliminary result, adapted from Carmona and Hu \cite{CarHu}.

\begin{Lemma} \label{lem:superadditive1}
There exists some positive constant $v^{(0)}$ such that the solution of
the KPZ equation with zero bare velocity, 
\BEQ (\partial_t-\nu^{(0)}\Del)h(t,x)=\sqrt{D^{(0)}}\eta(t,x)+ \lambda
|\nabla h(t,x)|^2 \EEQ
verifies
\BEQ \frac{1}{T}\langle h(T,x)\rangle\to_{T\to\infty} 
\liminf_{T\to\infty}  \frac{1}{T}\langle h(T,x)\rangle=:v^{(0)}.
\label{eq:Fekete}
\EEQ
Furthermore, $0\le v^{(0)}\le \tilde{v}^{(0)}$.

\end{Lemma}

\medskip\noindent {\bf Proof} (see \cite{CarHu}, Lemma 3.1)  Let, for $f$ general
forcing term, 
\BEQ w_T(a|f):= \esper^0 \left[ e^{g^{(0)}\int_0^T dt\, 
f(t,a+B_{T-t})  }  \right]
\EEQ

and 

\BEQ w_T(a,b|f):= \esper^0 \left[ e^{g^{(0)}\int_0^T dt\, 
 f(t,a+B_{T-t}) } 
 \ \big{|}\ a+B_T=b \right]  \label{eq:wab}
\EEQ
Conditioning with respect to the terminal condition, $a+B_T=b$, means that
we average with respect to the law of the {\em Brownian bridge from $(0,a)$ to $(T,b)$} (see e.g. \cite{KS}).
Then, for $T,T'\in\N$,
\BEA w_{T+T'}(x|\eta) &=& \int p_T(x,dy) w_T(x,y|\eta(\cdot+T'))
w_{T'}(y|\eta) \nonumber\\
&=& w_T(x| \eta(\cdot+T')) \int p_T(x,dy) \pi_{T,T'}(x,y|\eta(\cdot+T')) w_{T'}(y|\eta)
\EEA
where $\pi_{T,T'}(x,y|\eta(\cdot+T')):=\frac{w_T(x,y|\eta(\cdot+T'))}{w_T(x|\eta(\cdot+T'))}$. By construction, 
$\int p_T(x,dy) \pi_{T,T'}(x,y|\eta(\cdot+T'))=1$. Hence (by concavity of the log)
\BEQ h_{T+T'}(x)\ge  h(T,x) +  \int p_T(x,dy)\pi_{T,T'}(x,y|\eta(\cdot+T'))
h_{T'}(y).\EEQ
Taking the expectation with respect to the noise and using independence of 
$\eta(\cdot+T')$ from $\eta\big|_{[0,T']}$, together with space translation invariance, one gets the {\em superadditive inequality},
\BEQ \langle h_{T+T'}(x)\rangle=\langle h_{T+T'}(0)\rangle \ge 
\langle h_T(0)\rangle + \langle h_{T'}(0)\rangle.\EEQ
On the other hand, by convexity of exp, $\langle h_T(0)\rangle\ge 0$.
Fekete's superadditive lemma allows us to conclude to the existence
of some constant $v^{(0)}$ verifying (\ref{eq:Fekete}). This
is the constant whose existence is asserted in Main Theorem
(see (\ref{eq:0.5})).  Furthermore, by Jensen's inequality, $v^{(0)}\le \tilde{v}^{(0)}$, as observed already in \cite{CarHu}, Prop. 1.4.  
 \hfill \eop

\medskip\noindent As mentioned in the Introduction, Carmona and Hu \cite{CarHu} actually prove
the existence of a limit random variable $h_{\infty}(0):=$a.s.-lim$_{t\to\infty}
h(t,0)$ for the solution of the KPZ equation with velocity $\tilde{v}^{(0)}$, and 
a Gaussian lower large deviation theorem (Theorem 1.5 in \cite{CarHu}) for $h_{\infty}(0)$ of the form
\BEQ \proba[h_{\infty}(0)\le -A]\lesssim e^{-cA^2}, \qquad A>0 \EEQ
from which it is clear in particular that $v^{(0)}=\tilde{v}^{(0)}$.

\bigskip\noindent
Because the equation for $w$ is {\em linear}, there exists a random kernel
$G_{\eta}=G_{\eta}((t,x),(t',x'))$ ($t>t'$) such that 
\BEQ w_t(x|\eta)=\int dx'\, G_{\eta}((t,x),(t',x'))w_{t'}(x'|\eta).\EEQ
From the above formulas one sees that 
\BEQ G_{\eta}((T,a),(0,b))\equiv w_T(a,b|\eta).\EEQ
The kernel $G_{\eta}$, called {\em random propagator}, is the matter of the
next subsection.

%%%%%%%%%%%%%%%%%%%%%%%%%%%%%%%
%%%%%%%%%%%%%%%%%%%%%%%%%%%%%%%%%%

\section{Multi-scale expansion and vertex representation}

%%%%%%%%%%%%%%%%%%%%%%%%

 We discuss in this section two
different points of view on the KPZ equation (\ref{eq:KPZ}):

\begin{enumerate}
\item  First (see section 1), due to our specific choice of quadratic nonlinearity $V(\nabla h)=|\nabla h|^2$, the {\em Cole-Hopf transform}
maps (\ref{eq:KPZ}) into a {\em linear equation} for a {\em Cole-Hopf field} $w$ with multiplicative noise, which is
explicitly solved in terms of an average over Brownian paths, giving rise
to {\em Cole-Hopf solutions}. Conjugating with respect to the Cole-Hopf transform, 
these may be seen to coincide with the  ${\cal W}$-solutions introduced elsewhere
\cite{Unt-KPZ1}. This point of view,
in combination with
martingale theorems and Gaussian concentration inequalities,
is extensively used in the literature \cite{Bol,ImbSpe,CarHu,ComYos},  where people have been at least
as much interested in the resulting weighted measure on paths, interpreted as a
directed polymer measure.  A lot of properties of this measure have been derived
in all dimensions, in the small ($\lambda\ll 1$) or large ($\lambda\gg 1$) disorder regime, with attention focused on asymptotic theorems, large-deviation properties, scaling exponent, etc. However, not much can be derived therefrom concerning the asymptotic behavior of {\em $N$-point functions of the original
KPZ field $h$} for $N\ge 2$, because they are not directly accessible from the directed polymer
measure due to necessity of taking the {\em inverse} Cole-Hopf transform.  

\item Second (see \S \ref{subsection:vertex}) -- and this our approach here --, starting either directly from the KPZ equation or
the Cole-Hopf transformed linear equation, one may try to expand the solution in powers of $\lambda$ for $\lambda$ small enough.  In the first case, the idea is more or less to apply iteratively Duhamel expansion. In the second case, one is led to a {\em vertex representation}
based on an expansion of the random resolvent.

\end{enumerate}

\medskip\noindent The second point of view may look very naive  to mathematicians at first sight -- though physicists have long known how to build predictions out of perturbative
expansions --; such
approaches in PDE theory lead in general only to existence "in the small", i.e. for
a small enough initial condition. Because here we have a SPDE with a right-hand side,
one may expect to get only short-time existence.  However, it turns out that combining it to
very basic finite-time bounds for the solution in a finite box, and to the 
apparatus of cluster expansions and renormalization,  yields exact asymptotics for
$N$-point functions in the large-scale limit!  Thus this semi-perturbative approach
for $\lambda\ll 1$ is much more successful than previous approaches 1. and 2.,  whose results are not required, and actually can be  rederived directly up to some point. The
key point  is  to assess the precise amount of expansion needed to get the
leading large-scale behavior without producing at the same time diverging series.

%%%%%%%%%%%%%%%%%%%%%%%%%%%%

\subsection{Multi-scale decompositions and power-counting}

%%%%%%%%%%%%%%%%%%%%%%%%%
%%%%%%%%%%%%%%%%%%%%%%%%%%%

In the following somewhat technical section, we cut propagators into {\em scales}, and
space-time into {\em scaled boxes}, paving the way for the cluster expansions of
section 3. The more PDE-minded reader may find it more reassuring to read section 2
first, and then navigate between sections 1 and 3.

\medskip

\begin{Definition}[phase space]  \label{def:phase-space}
\begin{itemize}
\item[(i)] {\em (boxes)}  Let
$$\D^j:=\cup_{(k_0,\vec{k})\in\N\times\Z^d} [k_0 2^j,(k_0+1) 2^j)\times
[k_1 2^{j/2},(k_1+1)2^{j/2}]\times\cdots\times [k_d 2^{j/2},(k_d+1)2^{j/2}]$$
$(j\ge 0)$ and 
 $\D=\cup_{j=0}^{+\infty} \D^j$.
If $(t,x)\in\Del$ with $\Del\in\D^j$, we write $\Del^j_{(t,x)}:=\Del$.

\item[(ii)] {\em (momentum-decoupling $\tau$-parameters)}
If $\tau^0:\D^0\to[0,1]$, we write $\tau^0_t:=\tau(\Del^0_t)$.
% Similarly,
%$\tau_t:=(\tau^j_t)_j$ is the collection of all $\tau$-coefficients localized
% at $\tau$.

\item[(iii)] (space projection) If $\Del\in\Del^j$, $\Del:=[k_0 2^j,(k_0+1) 2^j)\times [k_1 2^j,(k_1+1)2^{j/2}]\times\ldots\times [k_d 2^{j/2},(k_d+1)2^{j/2}]$, we let
$\bar{\Del}:= [k_1 2^j,(k_1+1)2^{j/2}]\times\ldots\times [k_d 2^{j/2},(k_d+1)2^{j/2}]$.
Then $\bar{\D}^j$ is the union of all such cubes in $\R^d$.

\end{itemize}
\end{Definition}

Let $\nu>0$. We let  $G_{\nu}:=(\partial_t-\nu\Del)^{-1}$ be the heat kernel with diffusion coefficient $\nu$,
\BEQ G_{\nu}(t,x;t',x'):=p_{\nu (t-t')}(x-x') \ {\mathrm{if}} 
\ t,t'\ge 0 \ {\mathrm{and}}\ t-t'>0, \qquad 0 \ {\mathrm{else}}
\EEQ
where $p_{\nu(t-t')}(x-x')=\frac{e^{-|x-x'|^2/4\nu(t-t')}}{(4\pi\nu(t-t'))^{d/2}}$ is the kernel of the heat operator $e^{\nu(t-t')\Del}$. {\em When $\nu:=\nu^{(0)}$ is the bare  viscosity, we write simply $G_{\nu^{(0)}}=:G$. }

\medskip\noindent In the following definition, if $f:\R_+\to\R$, we let: $f^j(t):=f(2^{-j}t)$ ($j\ge 1$).

\begin{Definition}[multi-scale decompositions]  \label{def:Aj}
\label{def:one-step}

Choose a smooth partition of unity $1=\chi^{0}\ast\chi^{0} + \sum_{j=1}^{+\infty} (\chi\ast\chi)^j$ of $\R_+$  for some
smooth functions $\chi:\R_+\to [0,1]$ with compact support $\subset [\half,2]$, and
$\chi^{0}:\R_+\to [0,1]$ with compact support $\subset [0,1]$.
Let $A^j(t,t'),B^j(t,t')$ ($j\ge 0$, $t>t'>0$) be the operator-valued, time-convolution
kernels defined by
\BEQ A_{\nu}^0(t,t')\equiv B_{\nu}^0(t,t'):= \chi^{ 0}(t-t') e^{\nu (t-t')\Del}  \EEQ
and, for $j\ge 1$, 
\BEQ  A_{\nu}^j(t,t')\equiv B_{\nu}^j(t,t')  := 2^{-j/2}     \chi^j(t-t') e^{\nu (t-t')\Del}.\EEQ

They define operators $A^j_{\nu},B^j_{\nu}:L^2(\R_+\times\R^d)\to L^2(\R_+\times\R^d)$
through $(A^j f)(t):=\int_0^t dt'\, A^j(t,t')f(t'),\ (B^j f)(t):=\int_0^t dt'\, B^j(t,t')f(t')$.

%\item[(4)] Let $a^j\equiv b^j:\R_+\times\R^d\to\R$, $j\ge 0$  be the
% centered Gaussian process
%with covariance function $C^{j,j'}_a=C^{j,j'}_b$ defined by $C^{j,j'}_{a}
%(t,x;t',x'):=(-\Del)^{1/2} 
%\langle \nabla\phi^j(t,x) \nabla \phi^{j'}(t',x')\rangle$.

\end{Definition}

\medskip\noindent {\em Remark.}  {\em If $(t,x)$ is connected to $(t,'x')$ by some $A^j$ or $B^j$ with
$j\ge 1$, then $t-t'>1$, hence $\langle \eta(t,x)\eta(t',x')\rangle=0$.} This
property (due to an adequate choice of cut-offs) is convenient since it implies that
{\em two-point functions require only a scale 0 renormalization} (see \S \ref{subsection:two-point}).

\bigskip\noindent
Note that $(\chi\ast\chi)^j=(2^{-j/2}\chi^j)\ast(2^{-j/2}\chi^j)$ $(j\ge 1)$. Hence, by construction,
\begin{itemize}
\item {\em The $A^j_{\nu}$'s provide a  decomposition of the kernel $G_{\nu}$ into a sum
of positive kernels}: namely, 
\BEQ \sum_{j\ge 0} A_{\nu}^j B_{\nu}^j(t,t')=  (\chi^{0}\ast\chi^{0})(t-t')  e^{\nu (t-t')\Del}\, dt \, +\, \sum_{j\ge 1}  ( (2^{-j/2}\chi^j)\ast(2^{-j/2}\chi^j))(t-t')
e^{\nu (t-t')\Del}\, dt=G_{\nu}(t,t'). \label{eq:AjBj-dec} \EEQ
Furthermore, letting 
\BEQ G^j_{\nu}:= A^j_{\nu}B^j_{\nu}, \qquad j\ge 0\EEQ
we have  $\sum_{j\ge 0} G^j_{\nu}=G_{\nu}$, and  $G^j_{\nu}$ is "roughly" $2^{j/2} A^j_{\nu}$ (we say "roughly", because $(\chi\ast\chi)^j$ and $\chi^j$ do not have
exactly the same time support -- a more precise statement may be e.g. that 
$2^{j/2} A^j_{c\nu}(\cdot,\cdot)\lesssim G^j_{\nu}\lesssim  2^{j/2} A^j_{\nu/c}(\cdot,\cdot)$ for some $0<c<1$).
\end{itemize}

\begin{Definition} \label{def:A}
\begin{enumerate}
\item Let $A_{\nu}(\cdot;\cdot,\cdot)$ be the following kernel on $(\R_+\times\R^d)\times
(\N\times\R_+\times\R^d)$,
\BEQ A_{\nu}((t,x);j,(t',x')):= A_{\nu}^j((t,x),(t',x')).\EEQ
\item Let $B_{\nu}(\cdot,\cdot;\cdot)$ be the following kernel in $(\N\times\R_+\times\R^d)\times
(\R_+\times\R^d)$,
\BEQ B_{\nu}(j,(t,x);(t',x')):= B_{\nu}^j((t,x),(t',x')).\EEQ
\end{enumerate}
\end{Definition}

In other words, letting ${\cal H}$ be an auxiliary separable Hilbert space with
orthonormal basis denoted by $\vec{e}^j$, $j\ge 0$, or equivalently, 
$|j\rangle$ (in quantum mechanical notation),
 $A_{\nu}(\cdot,\cdot)$ is the kernel of the operator 
 \BEQ A_{\nu}: {\cal H}\otimes L^2(\R_+\times\R^d) \to L^2(\R_+\times\R^d) \EEQ
  defined by $A_{\nu}(\vec{e}^j \otimes f)=A^j_{\nu}(f)$; equivalently, $A_{\nu}:=\sum_{j\ge 0} A^j_{\nu} \langle j|$ has a linear form-valued kernel on $(\R_+\times\R^d)\times(\R_+\times\R^d)$, 
 \BEQ A_{\nu}(\cdot,\cdot)\equiv \sum_{j\ge 0} A^j(\cdot,\cdot) \langle j|. 
 \label{eq:A-vec} \EEQ
  Dualizing,
 $B_{\nu}(\cdot,\cdot)$ is  the kernel of the operator 
 \BEQ B_{\nu}: L^2(\R_+\times\R^d) \to  {\cal H}\otimes L^2(\R_+\times\R^d) \EEQ
  defined by $B_{\nu}(f)=\sum_{j\ge 0} B^j_{\nu}(f)
 \vec{e}^j$; in other words, $B_{\nu}:=\sum_{j\ge 0} B^j_{\nu}  |j\rangle$, with
 associated vector-valued kernel
 \BEQ  B_{\nu}(\cdot,\cdot)\equiv \sum_{j\ge 0} B^j(\cdot,\cdot) |j\rangle. 
 \label{eq:B-vec} \EEQ
 Thus
 the decomposition of $G_{\nu}$, see (\ref{eq:AjBj-dec}), is equivalent to the identity
\BEQ A_{\nu}B_{\nu}=\sum_{j,j'\ge 0} A^j_{\nu} B^{j'}_{\nu} \langle j|j'\rangle=
\sum_{j\ge 0} A^j_{\nu}B^j_{\nu}=G_{\nu} \label{eq:AB-dec}
\EEQ 
which lies at the core of the vertex representation in \S \ref{subsection:vertex}. 

\medskip\noindent {\em As in the case of  $G_{\nu}$, we write simply $A_{\nu^{(0)}}=:A$, $B_{\nu^{(0)}}=:B$.}

\medskip\noindent  The following
estimates for the kernel $A^{j}(t,x;t',x')=B^{j}(t,x;t',x')$ of
$A^j=B^j$ are easily shown:

\begin{Lemma}[multi-scale estimates for $A$ and $B$]  \label{lem:multi-scale-estimates-A-B} Let $j\ge 1$.

\begin{itemize}
\item[(i)] (single-scale estimates) 
\BEQ |\partial^{\kappa'}_t \nabla^{\vec{\kappa}} A^{j}(t,x;t',x')|\lesssim
(2^{-j/2})^{2\kappa'+|\vec{\kappa}|} (2^{-j/2})^{d+1} e^{-c2^{-j}|x-x'|^2} {\bf 1}_{t-t'\approx 2^j};\EEQ

\BEQ \int dt' \, dx'\, A^{j}(t,x;t',x') \approx 2^{j/2}; \EEQ

\BEQ ||A^j f||_{L^2} \lesssim  (2^{-j/2})^{d/2} ||f||_{L^2}.\EEQ

\item[(ii)] (two-scale estimates) let $1\le j$ and $\kappa,\kappa'\ge 0$, then
\BEQ |(\nabla^{\vec{\kappa}}A^j \ \nabla^{\vec{\kappa}'}B^{j})((t,x),(t',x'))| \lesssim
 (2^{-j/2})^{d+|\vec{\kappa}|+|\vec{\kappa'}|} e^{-c2^{-j}|x-x'|^2}
 {\bf 1}_{t-t'\approx 2^{j}}.
\label{eq:two-scale-estimate} \EEQ

\end{itemize}
\end{Lemma}

From (ii) it results that $(\nabla^{\vec{\kappa}} A^j \ \nabla^{\vec{\kappa}'} B^{j})(\cdot,\cdot)$ scales 
like $\nabla^{\vec{\kappa}+\vec{\kappa}'} G^{j'}(\cdot,\cdot)$ -- or, more precisely, like
$\nabla^{\kappa+\kappa'} G^{j}_{\nu}(\cdot,\cdot)$, with $\nu\approx \nu^{(0)}$, or equivalently, like
$2^{j/2} \nabla^{\vec{\kappa}+\vec{\kappa}'} B^{j}_{\nu}(\cdot,\cdot)$. Also, it is clear that
$G(\cdot,\cdot)\lesssim \sum_{k\ge 0} 2^{k/2} A^k_{\nu}(\cdot,\cdot)$. As  immediate corollary, expanding  $G$ over scales, it comes out
\BEQ  (B^j G)(\cdot,\cdot) \lesssim 2^{j/2} G^{\to j}_{\nu}(\cdot,\cdot).
\label{eq:GAB} \EEQ
\BEQ |(\nabla^3 G^j \, \cdot\,  G)(\cdot,\cdot)|,\qquad  |(\partial_t\nabla G^j\,  \cdot\, G)(\cdot,\cdot)| \lesssim
2^{-j/2} G^{\to j}_{\nu}(\cdot,\cdot) \EEQ
and finally {\em the first of our two key power-counting estimates},
\BEQ |(\nabla^3 G\,  \cdot\, G)(\cdot,\cdot)|,\qquad |(\partial_t\nabla G \, \cdot\, G)(\cdot,\cdot)|\lesssim G_{\nu}(\cdot,\cdot),  \label{eq:PW1}
\EEQ 
whereas $\partial_t^{\kappa_0}\nabla^{\vec{\kappa}} G\, \cdot\, G$, $|\kappa|:=2\kappa_0+|\vec{\kappa}|\equiv 2\kappa_0+\kappa_1+\ldots+\kappa_d$  {\em diverges in
the stationary limit when $|\kappa|\le 2$}, i.e. $ (\partial_t^{\kappa_0}\nabla^{\vec{\kappa}} G\, \cdot\, G)((t,x),(0,x))\approx t^{1-\kappa/2}$ 
$(|\kappa|<2)$, $\log(t)$ ($|\kappa|=2$), therefore $\overset{t\to +\infty}{\longrightarrow} +\infty$. In all these estimates it is intended that $\nu\approx
\nu^{(0)}$.
\medskip

{\em Proof.} 
\begin{itemize}
\item[(i)] 
Immediate consequence of the elementary heat kernel estimates,
$|\partial_t^{\kappa'}\nabla^{\vec{\kappa}}  p_{\nu(t-t')}(x)| \lesssim (t-t')^{-\kappa'-|\vec{\kappa}/2} p_{2\nu(t-t')}(x).$  Note that the time support and scaled exponential space decay leave an {\em effective space-time integration volume}
$O(2^{j(1+d/2)})$. The $L^2$-norm estimate is also a
consequence of: $||A^j f||_{L^2}^2\lesssim \int_{t-t'\approx 2^j} \frac{dt}{t-t'}
||e^{(t-t')\nu\Del}f_{t'}||_{L^2}^2$ and the easy inequality $||e^{(t-t')\nu\Del}f||_{L^2}^2\lesssim (t-t')^{-d/2}
||f||_{L^2}^2$ (standard parabolic estimate).

\item[(ii)]  Integrating $\int dt''\, \int dx''\, \nabla^{\vec{\kappa}}A^j((t,x),(t'',x''))
\nabla^{\vec{\kappa}'}B^{j}((t'',x''),(t',x'))$ by parts with respect to $t''$, and
remarking that $t''$ ranges in a time-interval of size $O(2^{j})$, we obtain 
\BEA  &&\Big| (\nabla^{\vec{\kappa}} A^j \ \nabla^{\vec{\kappa}'}B^{j})((t,x),(t',x')) \Big| =
\Big| (A^j\  \nabla^{\vec{\kappa}+\vec{\kappa}'}B^{j})((t,x),(t',x')) \Big| \nonumber\\
&&\qquad\lesssim  \left( 2^{j/2} e^{2^j \frac{\nu^{(0)}}{c} \Del}\, \cdot\,
(2^{-j/2})^{1+|\vec{\kappa}|+|\vec{\kappa}'|} e^{2^{j'} \frac{\nu^{(0)}}{c} \Del}  \right) ((t,x),(t',x')) \nonumber\\
&&\qquad \lesssim   (2^{-j/2})^{|\vec{\kappa}|+|\vec{\kappa}'|} e^{2^{j'}\frac{\nu^{(0)}}{c'}\Del}((t,x),(t',x')).
 \EEA
\end{itemize}
\hfill\eop

One gets similarly 
\BEQ  G^j(t,x;t',x') \lesssim  2^{-jd/2} e^{-c2^{-j}|x-x'|^2} \ 
{\bf 1}_{t-t'\approx 2^j}  \label{eq:power-counting-G}
\EEQ

%%%

%%%%%%%%%%%%%%%%%%%%%ש

\bigskip\noindent
At this point we introduce a very useful

\medskip
\noindent{\bf Universal notation:} let $f=\sum_{j=0}^{+\infty} f^{(j)}$ be a function/random field/multi-scale diagram/... decomposed into its scale components, then 
\BEQ f^{\to j}:=\sum_{k\ge j} f^{(k)}=\ldots+f^{(j+2)}+f^{(j+1)}+f^{(j)} \EEQ
is the {\em scale $j$ low-momentum part} of $f$, while
\BEQ f^{j\to}:=\sum_{k\le j} f^{(k)}=f^{(j)}+f^{(j-1)}+\ldots+f^{(1)}+f^{(0)} \EEQ
is the {\em scale $j$ high-momentum part} of $f$.

\medskip
\noindent
In the particular case of the kernels $A$ and $B$, the following is intended,
\BEQ A^{\to j}(\cdot,\cdot):=\sum_{k\ge j} A^k(\cdot,\cdot) \langle k|, \qquad A^{j\to}
(\cdot,\cdot):=\sum_{k\le j} A^k(\cdot,\cdot)
\langle k| \EEQ
\BEQ B^{\to j}(\cdot,\cdot):=\sum_{k\ge j} B^k(\cdot,\cdot) |k\rangle, \qquad B^{j\to}(\cdot,\cdot):=\sum_{k\le j} B^k(\cdot,\cdot)
|k\rangle.\EEQ

 %%%%%%%%%%%%%%%%%

%%%%%%%%%%%%%%%%%%%%

%%%%%%%%%%%%%%%%%%%%%%%%
%%%%%%%%%%%%%%%%%%%%%%%%%%%%שש

%%%%%%%%%%%%%%%%%%%%%%%%%%%%%%%%%

\subsection{The vertex representation}  \label{subsection:vertex}

%%%%%%%%%%%%%%%%%%%%%%%%%%%%

Consider the KPZ equation (\ref{eq:model}). 
Recall  $(\nu^{(0)},D^{(0)},v^{(0)})$ are the bare parameters. Expanding 
blindly the exponential in Feynman-Kac's formula (\ref{eq:Wiener-integral}) would
yield a series in  the bare
coupling constant $g^{(0)}=\frac{\lambda}{\nu^{(0)}}\sqrt{D^{(0)}}$.  This is the
starting point for our expansion. In the end (see section \ref{sec:bounds}), we
shall see that it is possible to make partial resummations, and obtain thus 
expressions bounded by products of {\em short-time} kernels $G_{\eta}((t,x),(t',x'))$ with
$t-t'=O(1)$, which are 
in turn bounded using
(\ref{eq:Wiener-integral}).

\bigskip

Let us start with some general considerations. Let $f=f(t,x)$ be any
right-hand side, and $\nu>0$. The integral version of the equation
\BEQ (\partial_t-\nu\Del)w(t,x)=f(t,x)w(t,x), \EEQ
coinciding -- up to the replacement of $\nu^{(0)}$ by $\nu$ -- with (\ref{eq:w}) when $f(t,x):=g^{(0)} ( \eta(t,x) - v^{(0)}),$
is
\BEQ w(t,x)=G_{\nu}((t,x),(0,\cdot))w_0(\cdot)+G_{\nu}((t,x),\cdot)(f w)(\cdot).\EEQ
Iterating yields 
\BEQ w(t,x)=\left( G_{\nu}+G_{\nu}f G_{\nu}+G_{\nu}fG_{\nu}f G_{\nu}+\cdots\right)
((t,x),(0,\cdot))w_0(\cdot).\EEQ
The series converges under suitable hypotheses on $f$, and the general term
in the series has the form of a chronological sequence, or {\em string} of propagators $G_{\nu}$
with $g$'s sandwiched in-between, namely,
\BEA && \left(G_{\nu}f\cdots f G_{\nu}\right)(t,x;0,y)  =\int_0^t dt_1 \int dx_1\, G_{\nu}(t,x;t_1,x_1)
f(t_1,x_1) \nonumber\\
&&\qquad  \int_0^{t_1} dt_2\int dx_2\, G_{\nu}(t_1,x_1;t_2,x_2) f(t_2,x_2)
\int_0^{t_2} dt_3\int dx_3\, \cdots \nonumber\\
 \EEA

We now turn to a representation in terms of the operators $A_{\nu},B_{\nu}$ defined in Definition \ref{def:A}
by means of the auxiliary space ${\cal H}$ indexing the scales.

\medskip\noindent 
{\em To an arbitrary function $f$, we associate the following {\em general vertex}}
  \BEQ V_{\nu}(f)(t,x):=B_{\nu}(\cdot,(t,x))f(t,x) A_{\nu}((t,x),\cdot). \label{def:general-vertex} \EEQ

\smallskip\noindent Since $A_{\nu}B_{\nu}=G_{\nu}$, one sees immediately by expanding
$(1-X)^{-1}=1+X+X^2+\cdots$ that
\BEQ w=A_{\nu} \left(1-\int dt\, dx\, V_{\nu}(f)(t,x)\right)^{-1}B_{\nu}w_0.\EEQ
Here $(1-\int dt\, dx\, V_{\nu}(f)(t,x))^{-1}$ plays manifestly the r\^ole of a {\em resolvent}. 

\medskip\noindent {\em Remark.} Other choices of vertices and scale decompositions are possible; for instance, letting instead $B_{\nu}\equiv A_{\nu}:=\sqrt{G_{\nu}}=\int_0^{+\infty} e^{\nu t\Del}
\frac{dt}{\sqrt{2t}}$, and decomposing $B_{\nu},A_{\nu}$ into scales in a similar 
way as we did in Definition \ref{def:Aj},  eq. \ref{def:general-vertex}
defines a {\em scalar} vertex. However, the orthogonal projection structure of
(\ref{eq:A-vec},\ref{eq:B-vec}) yields significant simplifications, see (\ref{eq:BAB})
and section 7.

\medskip\noindent  Let $\nu=\nu^{(0)}$. {\em Recall that we write for short in
this case $G\equiv G_{\nu^{(0)}},A\equiv A_{\nu^{(0)}},B\equiv B_{\nu^{(0)}}$.} Choosing  $f=g^{(0)}(\eta-v^{(0)})$, we obtain the

\begin{Definition}[Cole-Hopf vertex]  \label{def:Cole-Hopf-vertex}
\BEQ
V_{\eta}(t,x):=B(\cdot,(t,x)) \left(g^{(0)}(\eta(t,x)-v^{(0)}) \right) A((t,x),\cdot) \label{eq:V}
\EEQ
\end{Definition}

Then the solution of (\ref{eq:w}) is
\BEQ w=A \left(1-\int dt\, dx\, V_{\eta}(t,x)\right)^{-1}Bw_0. \label{A1-VB} \EEQ

In other words, letting

\begin{Definition}[random resolvent/propagator]
\BEQ R_{\eta}:=\left(1-\int dt\, dx\, V_{\eta}(t,x)\right)^{-1}, \qquad G_{\eta}:=AR_{\eta}B \EEQ
\end{Definition}
we have
\BEQ w(t,x)=(A R_{\eta} B)((t,x),(0,\cdot))w_0(\cdot)=G_{\eta}((t,x),(0,\cdot))w_0(\cdot).  \label{eq:Reta} \EEQ

%%%%%%%%%%%%%%%%%%%%%%%

\section{Cluster expansions}  \label{sec:cluster}

%%%%%%%%%%%%%%%%%%%%%%%%%%%%%%

The general principle of multi-scale expansions is that each field has one degree of freedom per box in $\D=\cup_{j\ge 0}\D^j$ (an idea made precise by wavelet expansions). In order to understand
the effect of the weak coupling between the degrees of freedom belonging to 
different boxes, one interpolates
between the totally decoupled theory and the coupled theory by introducing parameters.
These are of two kinds. {\em Horizontal parameters} (denoted by the letter {\em s})
test the coupling between two boxes of the same scale. {\em Vertical parameters}
(denoted by the letter {\em $\tau$}) test the coupling between a given box $\Del\in\D^j$, $j\ge 0$ and
the boxes {\em below} it, i.e. the boxes $\Del^k\in\D^k$, $k\ge j$ (one per scale) such
that $\Del^k\supset\Del^j$. (In the case of  the KPZ equation in its Cole-Hopf formulation, the only essential counterterms for
renormalization are produced at scale $0$, so we shall only test the coupling between
a box in $\D^0$ and the boxes below it). For the coupled theory,  these parameters are equal to $1$; for the totally decoupled theory, on the other hand, they are equal to $0$.
Taylor expanding to some order around $0$ with respect to the $s$- and $\tau$-parameters produces
in general  a combinatorial sum over products of so-called {\em multi-scale polymers} (unions of boxes). Any polymer is connected by links between boxes for which 
the relevant parameter, $s$ or $\tau$, is $>0$; such terms are written in terms of  Taylor integral remainders. In equilibrium statistical field theory, there
appear pieces totally isolated from remaining boxes; they
correspond to vacuum diagrams, and -- as well-known -- disappear when one computes connected expectations.
In our context, these do not appear ($Z=1$ automatically for dynamical theories,
because the noise measure is normalized from the beginning). On the other hand,
renormalization is in general a necessity
 in either setting, due to the following reason. Differentiating with respect to
 a $\tau$-link originated from a box $\Del^j\in\D^j$ produces low-momentum fields
 in some box $\Del^k\supset\Del^j$, $k> j$. Imagine one applies $\ge 1$ differentiations with respect to some of the vertical parameters located in boxes at
 the bottom of the polymer,  in total $N_{ext}$ of them, and then sets all of these
 vertical parameters to $0$. Thus this polymer
 "floats" at a certain height with respect to its external legs, measured by
 the difference $j_{ext,min}-j_{int,max}=$(min of scales $k$ of the $N_{ext}$ low-momentum fields) - (max of scales $j$ of
 bottom boxes). Then {\em the  quantity integrated in volume obtained by summing over
 all possible locations of the polymer with respect to its external legs}   is not
 a vacuum diagram; it is to be seen rather as some {\em insertion} contributing to
 the evaluation of the polymers located {\em below}. Computations show that, for
 $N_{ext}$ small enough (in our case, for $N_{ext}=2$ only), this contribution diverges in the limit when $j_{ext,min}-j_{int,max}\to\infty$. Thus such insertions
 contribute to the large-scale limit. The idea of Wilson's renormalization scheme
 is to {\em absorb the diverging part of these insertions into a scale by scale redefinition of the parameters of the theory}.
 
 \noindent
Here an essential simplification comes through the fact
that only scale 0 diagrams need to be renormalized, but the
general philosophy remains the same. 
 
\bigskip
\noindent In most theories, $N$-point functions are of the form $\langle P(h)\rangle$, where
$P(h)$ is a polynomial in the random field $h=h(t,x)$; however, here $h$ is the
{\em logarithm} of $w$. This a feature specific to this particular model.   Let us write down here explicitly the effect of successive differentiations on
an expression of the form $P(\log h)$. Incorporating the interpolating parameters transforms $w(t,x)$ 
into $w(\tau^0,\vec{s};t,x)$, where $\vec{s}$ and $\tau^0$ are scale 0 parameters.  Now, we need to differentiate with respect
to $s$- and $\tau^0-$parameters the $N$-point function $\langle \log(w_1(\tau^0,\vec{s})) \cdots \log(w_N(\tau^0,\vec{s}))\rangle$, where we have let
$w_k(\cdot):=w(\cdot;t_k,x_k)$. Then (letting $D_1,D_2,\cdots$ denote the
derivative with respect to various $s$- or $\tau^0$-parameters)
\BEQ D_1 \log(w_k(\cdot))=\frac{D_1 w_k(\cdot)}{w_k(\cdot)},\ 
D_2 D_1 \log(w_k(\cdot))=\frac{D_2 D_1 w_k(\cdot)}{w_k(\cdot)} -
\frac{D_1 w_k(\cdot) D_2 w_i(\cdot)}{(w_k(\cdot))^2}, \cdots
\EEQ 
\BEQ D_n\cdots D_1 \log(w_k(\cdot))=\sum_{m=1}^n (-1)^{m-1} (m-1)! \sum_{i_1+\ldots+i_m=n} \sum_{I_1,\ldots,I_m} \frac{ \left[\prod_{i\in I_1} D_i \right]w_k(\cdot) \cdots
\left[\prod_{i\in I_m} D_m\right] w_k(\cdot)}{(w_k(\cdot))^m},
\label{eq:DDD}
\EEQ
where the last sum ranges over all partitions of $\{1,\ldots,n\}$ into $m$ disjoint
subsets $I_1,\ldots,I_m$ with $|I_1|=i_1,\ldots,|I_m|=i_m$. Thus the derivatives
apply to a product $w_{i,1}\cdots w_{i,k}$ of $k$ "replicas" of $w_i$. The latter expression
generalizes  easily to some  combinatorial expression of the same type for
$D_n\cdots D_1 \left\{ \log w_1(\cdot))\cdots \log w_N(\cdot)) \right\}$ which
is of the general form
\BEQ  \sum_{(I_{k,j})} c_{\vec{I}}
\frac{ \prod_{k\le N} \prod_{j\le m_k}\Big( \left[\prod_{i\in I_{k,j}} D_i\right]w_k(\cdot) \Big) }{\prod_k
(w_k(\cdot))^{m_k}}     \label{eq:DDDD}
 \EEQ
with  $\uplus_k \uplus_{j\le m_k} I_{k,j}=\{1,\ldots,n\}$, for some coefficients
$c_{\vec{I}}$ depending on the choice of the sets $(I_{k,j})$, plus terms involving
one or several $\log w$ which have not been differentiated.
 The conclusion of this discussion is that we need only
 evaluate the $D_i$'s on so-called {\em replicated products} $\prod_{k=1}^{N'} w_k(\vec{s},\tau^0;t,x)$ taking
 into account the "replicas", or equivalently
   (by (\ref{eq:Reta})) on a product of 
  $N':=\sum_{k\le N} m_k$ noisy resolvents $R_{\eta}$. This is what we do in the next paragraphs.

\bigskip\noindent Let us finally mention our {\em implicit integration convention: whenever
a formula contains more space-time variables in the r.-h.s. than in the l.-h.s.,
supplementary variables are implicitly integrated over}. 

%%%%%%%%%%%%%%%%%%

\subsection{The dressed equation}

%%%%%%%%%%%%%%%%%%%%%%%%%%%%%

We now proceed -- as a  preparation to the renormalization
step --  to separate the $0$-th scale from the others. The outcome is a
"dressed" vertex $V(\tau^0)$. Let $\tau^0:\D^0\to[0,1]$. First we need to
dress the operators $A,B$.

\begin{Definition}[dressed fields]
\begin{enumerate}
\item ($A$-field)

Let 
\BEQ  A(\tau^0;(t,x),\cdot) := A^{0}((t,x),\cdot) \langle 0| \, +\, \tau^0_{(t,x)} A^{\to 1}((t,x),\cdot)
 \EEQ

\item  ($B$-field) The dressing procedure is the same, except that it acts
on the second set of variables, namely,
\BEQ B(\tau^0;\cdot,(t',x')) := B^{0}(\cdot,(t',x')) |0\rangle\, +\, \tau^0_{(t',x')} B^{\to 1}(\cdot,(t',x'))
 \EEQ

\end{enumerate}
\end{Definition}

The idea is the following. Start from a space-time dependent  field, say,
 $\phi(t,x)$, and make it $\tau$-dependent
as indicated.  Then Taylor's formula,
$\phi(\tau_{t,x}^0;t,x)=\phi(0;t,x)+\tau^0_{t,x} \partial_{\tau^0_{t,x}} \phi(0;t,x)+\ldots$ reads simply
$\phi(\tau_{t,x}^0;t,x)=\phi^{(0)}(t,x)+\tau_{t,x}^0 \phi^{\to 1}(t,x)$.  In other words, by differentiating
$\phi(\tau_{t,x}^0;t,x)$ with respect to $\tau_{t,x}^0$, one separates the zeroth scale component $\phi^{(0)}(t,x)$ from the low-momentum  field $\phi^{\to 1}(t,x)$.

\medskip

Renormalization involves a priori  the introduction of scale counterterms 
$\del g^{(j)}:=g^{(j+1)}-g^{(j)}$ (recall  $g^{(j)}:=\frac{\lambda}{\nu_{eff}} \sqrt{D^{(j)}}$ by definition), $\del v^{(j)}:=v^{(j+1)}-v^{(j)}$, 
$\del\nu^{(j)}:=\nu^{(j+1)}-\nu^{(j)}$. Due to our hypotheses on the
covariance kernel $\langle \eta(t,x)\eta(t',x')\rangle$, 
 it actually happens (as proved in section
\ref{sec:renormalization}) that only two-point {\em scale $0$} diagrams
absolutely
need renormalization; thus we choose to take  $\del g^{(j)}=0$ for all $j\ge 0$, and $\del \nu^{(j)},\del v^{(j)}\equiv 0$
for every $j\ge 1$. Since we want $\nu^{(j)}\to_{j\to\infty} \nu_{eff}$,
 $v^{(j)}\to_{j\to\infty} 0$, this implies simply that $g^{(j)}=g^{(0)},\, 
 \nu^{(j)}=\nu_{eff}, \, v^{(j)}=0$ for all $j\ge 1$.
Thus dressing the vertex is a very simple matter. First 
 (in order to avoid having to differentiate characteristic functions of 
scale $0$ boxes coming out of the horizontal cluster, see \S \ref{subsection:hor-cluster}),
we introduce
\BEQ
\Del^{\to 0}:=\bar{\chi}^{(0)}\ast\Del,  \label{eq:Del->0}
\EEQ
where $\bar{\chi}^{(0)}:\R^d\to \R$ is any normalized smooth "bump" function,
such that e.g. $\supp(\bar{\chi}^{(0)})\subset B(0,1)$, $\int dx\, \bar{\chi}^{(0)}(x)=1$; $\Del^{\to 0}$ is  a regularized version of $\Del$. It is useful
to assume that $\bar{\chi}^{(0)}$ is isotropic though (see section 7), which
improves the precision of the asymptotics in Theorem 0.1.

\begin{Definition}[dressed vertex and effective propagators]
\label{def:dressed-vertex}
Let, for $\tau^0:\D^0\to[0,1]$,
\begin{itemize}
\item[(i)]
\BEA && V_{\eta}(\tau^0;t,x) := B(\tau^0;\cdot,(t,x)) \left( g^{(0)}(\eta(t,x)-v^{(0)}) \right) A(\tau^0;(t,x),\cdot) \nonumber\\
&& \qquad
+  B^{\to 1}(\cdot,(t,x)) \left( (1-(\tau^0_{t,x})^2) (\nu_{eff}-\nu^{(0)}) \Del^{\to 0} \right) A^{\to 1}((t,x),\cdot)  \label{eq:dressed-vertex}
\EEA

\item[(ii)] 
\BEQ R_{\eta}(\tau^0):=(1-\int dt\, dx\, V_{\eta}(\tau^0;t,x))^{-1}.\EEQ
\end{itemize}
\end{Definition}

\medskip\noindent Let us comment formula (\ref{eq:dressed-vertex}), which is the starting point of all
subsequent computations.

\medskip\noindent The {\em first} line of (\ref{eq:dressed-vertex}),
\BEQ V^{(0)}_{\eta}(\tau^0;t,x) := B(\tau^0;\cdot,(t,x)) \left( g^{(0)}(\eta(t,x)-v^{(0)}) \right) A(\tau^0;(t,x),\cdot) \EEQ
is simply a dressed version of the Cole-Hopf vertex (\ref{eq:V}).

\medskip\noindent The {\em second} line,
\BEQ \del V_{\eta}(\tau^0;t,x):=B^{\to 1}(\cdot,(t,x)) \left( (1-(\tau^0_{t,x})^2) (\nu_{eff}-\nu^{(0)}) \Del^{\to 0} \right) A^{\to 1}((t,x),\cdot) 
\label{eq:del-V-eta} \EEQ
{\em vanishes when $\tau^0\equiv 1$}, which ensures that one recovers the original
Cole-Hopf vertex, i.e. $V_{\eta}(\tau^0\equiv 1;\cdot)=V_{\eta}(\cdot).$ 
It may be decomposed into two pieces, which are proportional but play a very different
r\^ole. The first one, $-(\tau^0_{t,x})^2   B^{\to 1}(\cdot,(t,x)) \left( (\nu_{eff}-\nu^{(0)}) \Del^{\to 0} \right) A^{\to 1}((t,x),\cdot)$, is a low-momentum counterterm
which resums the corresponding zero-momentum contribution of scale $0$ two-point functions (see \S \ref{subsection:two-point}). The second one,\\ $+B^{\to 1}(\cdot,(t,x)) \left( (\nu_{eff}-\nu^{(0)}) \Del^{\to 0} \right) A^{\to 1}((t,x),\cdot)$, leads to an
{\em effective propagator} 
\BEA
 \tilde{G}_{eff} &:=& A^{\to 1} \ \cdot\  \sum_{n\ge 0} \Big(\del V_{\eta}(\tau^0\equiv 0)\Big)^n \ \cdot\  B^{\to 1} \nonumber\\
&=& A^{\to 1} \left(1-(\nu_{eff}-\nu^{(0)})  B^{\to 1} \Del^{\to 0} A^{\to 1} \right)^{-1} B^{\to 1}
\EEA
which plays an essential r\^ole in the large-scale limit discussed in section 5. As
proved in Lemma \ref{lem:7}, $\tilde{G}_{eff}$ may be replaced in that limit by
$G_{eff}:=(\partial_t-\nu_{eff}\Del)^{-1}$ with an excellent approximation. Thus
$\nu_{eff}$ is, indeed, an effective viscosity. Namely, {\em it is shown in \S 7 that} 
\BEQ \tilde{G}_{eff}((\eps^{-1}t,\eps^{-1/2}x),(\eps^{-1}t',\eps^{-1/2}x')=G_{eff}((\eps^{-1}t,\eps^{-1/2}x),(\eps^{-1}t',\eps^{-1/2}x')) +"O(\eps)", 
\label{eq:O(eps)}
\EEQ
meaning the following (see Lemma \ref{lem:7}). Assume $t-t'\approx 1$ and $\eps\approx 2^{-j}\ll 1$, so that
$\eps^{-1}(t-t')\approx 2^j$. {\em Then the error term $"O(\eps)"$ is  equal to 
$O(\eps)$ times an exponentially decreasing kernel which is bounded by
$G_{\nu^{(0)}+O(\lambda^2)}((\eps^{-1}t,\eps^{-1/2}x),(\eps^{-1}t',\eps^{-1/2}x')) =
\eps^{d/2} G_{\nu^{(0)}+O(\lambda^2)}((t,x),(t',x'))$ in a very large space-time
region including the "normal regime"  $\frac{|x-x'|^2}{t-t'}\lesssim 1$.}

%%%%%%%%%%%%%%%%%%ש

\subsection{Horizontal cluster expansion} \label{subsection:hor-cluster}

%%%%%%%%%%%%%%%%%%%%%

The general principle  is outlined in  section 6.
We only need   a  scale 0 cluster expansion, which  we apply using (\ref{eq:BKAR}) to
\BEQ F\equiv F(A^0,B^0|\eta;A^{\to 1},B^{\to 1}):=\log(w_1(\tau^0,\vec{s}=1)\cdots
\log(w_N(\tau^0,\vec{s}=1)) \label{eq:specific-F}, \EEQ
 where
the $A^{\to 1}$'s and $B^{\to 1}$'s are only spectators. To be specific, $w(\tau^0,\vec{s})$ in the above expression is
defined as follows: 
\BEQ w(\tau^0,\vec{s};t,x)=(AR_{\eta}(\tau^0)(\vec{s}) B)((t,x),(0,\cdot))w_0(\cdot); 
\label{eq:wtaus} \EEQ

\BEA && R_{\eta}(\tau^0)(\vec{s})(t,x;t',x'):=\del(t-t')\del(x-x')+  \sum_{n=1}^{\infty}\int dx_{1}...dx_{n} \nonumber\\
&&\ \ 
\prod_{i=1}^n \Big[  \ \  \left( s_{\Del^0_{t_{i-1},x_{i-1}},\Del^0_{t_i,x_i}} B^0((t_{i-1},x_{i-1}),
(t_i,x_i))\, |0\rangle\, + \tau^0_{t_i,x_i} B^{\to 1} ((t_{i-1},x_{i-1}),(t_i,x_i)) \right) 
  \nonumber\\
&& \qquad\ \  \ \cdot\  (g^{(0)}(\eta(t_i,x_i)-v^{(0)})) \nonumber\\
&&\qquad\ \ \cdot\   \left( s_{\Del^0_{t_{i},x_{i}},\Del^0_{t_{i+1},x_{i+1}}} A^0((t_{i},x_{i}),
(t_{i+1},x_{i+1}))\, \langle 0|\, + \tau^0_{t_i,x_i} A^{\to 1} ((t_{i},x_{i}),(t_{i+1},x_{i+1})) \right)  \nonumber\\ 
&& \qquad
\  +\   B^{\to 1}(t_{i-1},x_{i-1}),(t_i,x_i)) \left( (1-(\tau^0_{t_i,x_i})^2) (\nu_{eff}-\nu^{(0)}) \Del^{\to 0} \right) A^{\to 1}((t_i,x_i),(t_{i+1},x_{i+1}))\  \Big]\nonumber\\  \label{eq:Rtaus} 
\EEA
where (by convention) $(t_0,x_0)\equiv(t,x),(t_{n+1},x_{n+1})\equiv (t',x')$. This way,
$F$ appears as a functional of $A^0,B^0$, to which the BKAR cluster expansion
formula (\ref{eq:BKAR}) applies. 

\medskip\noindent The outcome is an
{\em expression of $F$ in terms of a sum over scale $0$ forests $\F^0$},

\BEQ  \langle F(A^0,B^0|\eta) \rangle =\sum_{\F^0\in{\cal F}^0}  \left(
\prod_{\ell\in L(\F^0)}\int_0^1   dw_{\ell}\right) \left(  \left(\prod_{\ell\in L(\F^0)} \frac{d}{d s_{\ell}}\right) 
\langle F(A^0(\vec{s}(\vec{w})),B^0(\vec{s}(w)))|\eta \rangle_{\vec{s}(\vec{w})} \right),   \label{eq:appl-BKAR-F} \EEQ
see section 6 for detailed notations.

\medskip\noindent Let $\ell=(\Del,\Del')$, $\Del,\Del'\in\D^0$ be a pair of linked boxes.
We use the shortened notation $V_{\eta}(\tau^0)(\vec{s}(\vec{w})):=V_{\eta}(\tau^0)(A^0(\vec{s}(w)),B^0(\vec{s}(w)))$ and $R_{\eta}(\tau^0)(\vec{s}(\vec{w})):=\frac{1}{1-\int dt\, dx\, V_{\eta}(\tau^0;t,x)(A^0(\vec{s}(\vec{w})),B^0(\vec{s}(w)))}$.  A direct computation yields
\BEQ \frac{\partial}{\partial s_{\ell}} R_{\eta}(\tau^0)(\vec{s}(\vec{w}))= R_{\eta}(\tau^0)(\vec{s}(\vec{w}))
\left( \frac{d}{ds_{\ell}} \int dt\, dx\, V(\tau^0;t,x)(\vec{s}(\vec{w})) \right)   R_{\eta}(\tau^0)(\vec{s}(\vec{w}))  \label{eq:3.8}
\EEQ
Then
\BEA &&  \frac{\partial}{\partial s_{\ell}} \int dt\, dx\,  V_{\eta}(\tau^0;t,x)(\vec{s}(\vec{w})) =
\int dt\, dx\,
B(\tau^{0},\vec{s}(\vec{w}))(\cdot,(t,x))\ \cdot \nonumber\\
&&\qquad \cdot\  \left( g^{(0)}(\eta(t,x)-v^{(0)}) \right) \left(\frac{d}{ds_{\ell}} A(\tau^{(0)},\vec{s}(\vec{w}))((t,x),\cdot)\right)  \nonumber\\
&& \ +\int dt'\, dx'\, \left(\frac{d}{ds_{\ell}} B(\tau^{(0)},\vec{s}(\vec{w}))(\cdot,(t',x')) \right) \ \cdot \nonumber\\
&&\qquad\cdot\  \left( g^{(0)}(\eta(t',x')-v^{(0)}) \right) \left( A(\tau^{(0)},\vec{s}(\vec{w}))((t',x'),\cdot)\right)  \nonumber\\
\label{eq:3.9}
\EEA

Finally, if $(t,x)\in\Del$, $(t',x')\in\Del'$, $\Del,\Del'\in\D^0$,
\BEA \frac{\partial }{\partial s_{\ell}}A(\tau^0,\vec{s}(\vec{w})((t,x),(t',x')) &=& \frac{\partial}{\partial s_{\ell}}A^0(\vec{s}(\vec{w}))((t,x),(t',x')) \langle 0| \nonumber\\
&=&  A^0((t,x),(t',x')) \langle 0| \ \cdot\ 
{\bf 1}_{\ell=\{\Del,\Del'\}} \label{eq:3.10} \EEA
and similarly for $B$. Hence $(t,x)$, resp. $(t',x')$ in (\ref{eq:3.9}) is integrated
over $\Del$, resp. $\Del'$.

\medskip\noindent On the other hand (see (\ref{eq:d/ds-eta})), {\em $\frac{d}{ds_{\ell}}$
also acts on the covariance kernel of $\eta$, according to the rules:}
\BEQ \frac{d}{ds_{\ell}} \Big\langle \big(\, \cdot\, \big) \Big\rangle_{\vec{s}(\vec{w})} \equiv \int_{\Del_{\ell}}dz_{\ell} \int_{\Del'_{\ell}}dz'_{\ell} \, \langle \eta(z_{\ell})\eta(z'_{\ell})\rangle_{\vec{s}=1} \ \cdot\  \Big\langle
 \frac{\del}{\del\eta(z_{\ell})} \frac{\del}{\del \eta(z'_{\ell})}
\big(\, \cdot\, \big) \Big\rangle_{\vec{s}(\vec{w})}  \label{eq:3.11}
\EEQ

\BEQ \frac{\del}{\del \eta(z_{\ell})} R_{\eta}(\tau^0)(\vec{s}(\vec{w}))= R_{\eta}(\tau^0)(\vec{s}(\vec{w}))
\left( \frac{\del}{\del \eta(z_{\ell})} \int dz\,  V(\tau^0)(\vec{s}(\vec{w}))(z) \right)   R_{\eta}(\tau^0)(\vec{s}(\vec{w}))  \label{eq:3.12}
\EEQ

\BEQ  \frac{\del}{\del \eta(z_{\ell})}  \int dz\, V(\tau^0)(\vec{s}(\vec{w}))(z)= B(\tau^0,\vec{s}(\vec{w}))(\cdot,z_{\ell})  g^{(0)} A(\tau^0,\vec{s}(\vec{w}))(z_{\ell},\cdot), \label{eq:3.13} \EEQ

with now averages defined with respect to the $\vec{s}$-dependent Gaussian measure
$\langle \, \cdot\, \rangle_{\vec{s}(\vec{w})}.$ 

\medskip\noindent Clearly, $\frac{d}{ds_{\ell}}$ (or $\frac{\del}{\del \eta(z)}$,
$z=z_{\ell}$ or $z'_{\ell}$) can also act directly on one of the $A(\tau^0,\vec{s}(\vec{w}))(\cdot,\cdot)$, $B(\tau^0,\vec{s}(\vec{w}))(\cdot,\cdot)$ or $\eta$'s produced by previous differentiations. 

\bigskip\noindent Summarizing, turning to  the specific case (\ref{eq:specific-F}),  the result of the expansion (\ref{eq:appl-BKAR-F})
may be rewritten, using the notations of   (\ref{eq:DDDD}), and separating
 the action of the $s$-derivatives on the covariance kernel of $\eta$ from the action on the propagators $A^0,B^0$, and splitting the
 $s$-derivatives according to the index $(k,j)$ of the string
on which they act -- or possibly the {\em pair} of indices $(k,j),(k',j')$ for $\eta$-pairings between two
different strings -- 
\BEA &&  \sum_{\F^0\in{\cal F}^0}  \left(
\prod_{\ell\in L(\F^0)}\int_0^1   dw_{\ell}\right)  
\sum_{(L_G)_{k,j},(L_{\eta})_{k,j},(L_{\eta})_{(k,j),(k',j')}}   c_{\vec{I}} \ \  \Big(
\prod_{\ell\in L_{\eta}
} \int_{\Del_{\ell}}dz_{\ell}\ \int_{\Del'_{\ell}}dz'_{\ell}  \langle \eta(z_{\ell})\eta(z'_{\ell})\rangle_{\vec{s}=1} \Big)  \nonumber\\
&&  \ \Big\langle  
\frac{ \prod_{k\le N} \prod_{j\le m_k} \Big( \Big[  (D_G)_{k,j} (D_{\eta})_{k,j} (D_{\eta})_{(k,j),\cdot} (D_{\eta})_{\cdot,(k,j)}  
  \Big]
w_k(\tau^0,\vec{s};\cdot)  \Big)}{\prod_{k\le N}
(w_k(\tau^0,\vec{s};\cdot))^{m_k}}    \Big\rangle_{\vec{s}(\vec{w})} , 
\nonumber\\
  \label{eq:appl-BKAR-F2} \EEA
  where: $c_{\vec{I}}$ is as in (\ref{eq:DDDD});\\
  $L(\F^0)=L_G\uplus L_{\eta}$;\\
   $L_G=\uplus_{(k,j)} (L_G)_{k,j}$ ({\em propagator links}); \\
$L_{\eta}=\uplus_{k,j} (L_{\eta})_{k,j}\uplus_{(k,j),(k',j')} 
(L_{\eta})_{(k,j),(k',j')}$ ({\em noise links}); \\
$(L_{\eta})_{(k,j),\cdot}:=\uplus_{(k',j')\not=(k,j)} (L_{\eta})_{(k,j),(k',j')}, (L_{\eta})_{\cdot,(k,j)}:=\uplus_{(k',j')\not=(k,j)} (L_{\eta})_{(k',j'),(k,j)}$ (noise links
between two strings); \\
$(D_G)_{k,j}:= \prod_{\ell\in (L_G)_{k,j}} \frac{\partial}{\partial s_{\ell}}$ (derivatives acting on propagators $A^0$ or $B^0$);\\
 $(D_{\eta})_{k,j}:= \prod_{\ell\in (L_{\eta})_{k,j}} \frac{\del^2}{\del \eta(z_{\ell})\del \eta(z'_{\ell})} $ (double derivatives acting
 on two noise fields located on the same string);\\
$(D_{\eta})_{(k,j),\cdot}:=\prod_{\ell\in 
(L_{\eta})_{(k,j),\cdot}} \frac{\del}{\del\eta(z_{\ell})\del\eta(z'_{\ell})}, \ 
(D_{\eta})_{\cdot,(k,j)}:= \prod_{\ell'\in (L_{\eta})_{\cdot,(k,j)}} \frac{\del}{\del \eta(z_{\ell})\del\eta(z'_{\ell})}$ (resp. on two different
strings, including that of index $(k,j)$);\\
$m_k=$Card$\Big\{j | (L_G)_{k,j}\cup (L_{\eta})_{k,j} \cup (L_{\eta})_{(k,j),\cdot} \cup (L_{\eta})_{\cdot,(k,j)} \not=\emptyset\Big\}$

\medskip
  with $w_k(\tau^0;\vec{s})$ defined as in (\ref{eq:wtaus},\ref{eq:Rtaus}).

\bigskip\noindent  In other words:
\begin{itemize}
\item[(i)]  (see (\ref{eq:3.8},\ref{eq:3.9},\ref{eq:3.10})), each $s$-derivative along a link acting on a random resolvent  (i) singles out
a {\em localized} $A^0$- or $B^0$-propagator between the two boxes connected by the link,
and produces (ii) a supplementary  $B-$, resp. $A-$
propagator  ending, resp. starting in one of the two boxes; (iii)
a "renormalized" noise field 
\BEQ \tilde{\eta}(t,x):=g^{(0)}(\eta(t,x)-v^{(0)}) \label{eq:etatilde} \EEQ
 sandwiched between the  localized scale $0$ propagator, and another propagator with unspecified scale; (iv) and supplementary resolvents $R_{\eta}(\tau^0)(\vec{s}(\vec{w}))$,  whose scale 0 components $R_{\eta}^{(0)}$ will later
on be produced explicitly by the vertical expansion. Because all these scale 0 operators are causal, they may be seen as {\bf beads} stringed on an (open) {\bf string} propagating causally, with {\bf dangling $\tilde{\eta}$-ends} on each bead. See Fig. 1 below.

\smallskip\noindent  {\em Sequences $\int_{\Del} dt\, dx\, B^{\bullet}(\cdot,(t,x))
\tilde{\eta}(t,x)A^{\bullet}((t,x),\cdot)$ integrated in a box $\Del\in\D^0$,   are called}  {\bf vertices} by
reference to Definition \ref{def:dressed-vertex}.  
\item[(ii)] an $s$-derivative acting directly on some $A$ or $B$ turns into an
$A^0$ or $B^0$ linking two specified boxes;
\item[(iii)] the cluster in $\eta$ (see in particular
 (\ref{eq:3.11},\ref{eq:3.12},\ref{eq:3.13})) produces from $0$ to $2$ vertices (depending on whether the $\frac{\del}{\del \eta(z_{\ell})},\frac{\del}{\del \eta(z'_{\ell})}$ act on a resolvent or directly on some dangling
$\tilde{\eta}$), and a {\em local link between two vertices}, by which we mean that one gets
some {\bf pairing of} (old or new) {\bf  vertices} $\int_{\Del} dz\, B^{\bullet}(\cdot,z) A^{\bullet}(z,\cdot)$, 
$\int_{\Del'} dz'\, B(\cdot,z')A(z',\cdot)$, multiplied with the finite-range kernel 
$\langle \eta(z)\eta(z')\rangle_{\vec{s}=1}$, which forces $d(\Del,\Del')=O(1)$.

\end{itemize}

\bigskip
A general term in (\ref{eq:appl-BKAR-F2}) is in the form of a product of  $N'$  strings with  {\em beads} or  {\em inserted vertices} and {\em dangling $\tilde{\eta}$-ends}, schematically,
letting $z_i^j:=(t_i^j,x_i^j)$  ($1\le i\le N'$, $j\ge 1$) be intermediate
coordinates implicitly integrated over with $t_i>t_i^1>t_i^2>\ldots>t_i^{3n_i}
\equiv 0$,
\BEA 
&& \left( \prod_{j=1}^{n_1-1} \tilde{\eta}(z_1^{3j})  \right) \  A^{\bullet}((t_1,x_1),z_1^1) \nonumber\\
&&\ \ \left( R_{\eta}(z^1_1,z^2_1) \prod_{j=1}^{n_1-1} B^{\bullet}(z^{3j-1}_1,
z_1^{3j})A^{\bullet}(z_1^{3j},z_1^{3j+1}) R_{\eta}(z_1^{3j+1},z_1^{3j+2})  \right) B^{\bullet}
(z_1^{3n_1-1},z_1^{3n_1}) w_0(z_1^{3n_1})
\nonumber\\
&& \nonumber\\
&&\qquad\qquad\qquad\qquad\qquad\qquad \qquad\qquad\vdots  \nonumber\\
&&\nonumber\\
&& \left( \prod_{j=1}^{n_{N'}-1} \tilde{\eta}(z_{N'}^{3j})  \right) \ 
A^{\bullet}((t_{N'},x_{N'}),z_{N'}^1) \nonumber\\
&&\ \ \left( R_{\eta}(z^1_{N'},z^2_{N'}) \prod_{j=1}^{n_{N'}-1} B^{\bullet}(z^{3j-1}_{N'},
z_{N'}^{3j})A^{\bullet}(z_{N'}^{3j},z_{N'}^{3j+1}) R_{\eta}(z_{N'}^{3j+1},z_{N'}^{3j+2})  \right) B^{\bullet}
(z_{N'}^{3n_{N'}-1},z_{N'}^{3n_{N'}}) w_0(z_{N'}^{3n_{N'}})\nonumber\\
\label{eq:3.24}
\EEA
averaged w.r. to the measure $\langle \ \cdot\ \rangle_{\vec{s}(\vec{w})}$, where some of the $B$'s and $A$'s are localized, $0$-th scale  propagators, others being
"grey" for the moment (i.e. of unspecified scale), and $\tilde{\eta}(\cdot)=g^{(0)}(\eta(\cdot)-v^{(0)})$, see eq. (\ref{eq:etatilde}).
As seen from the previous formulas in this very subsection, such terms should be
summed over forests, integrated w.r. to interpolation coefficients $\vec{w}$. 
Intermediate 
coordinates $z^j_i$ are integrated over $0$-scale boxes $\Del_{\ell},\Del'_{\ell}$. 
Also missing are coefficients
$c_{\vec{I}}(\vec{z})$ now depending on  $\vec{z}:=
(z_{\ell},z'_{\ell})_{\ell\in 
L({\cal F}^0)}$  through the pairing factors  $\langle \eta(z_{\ell})\eta(z'_{\ell})\rangle$ due to the cluster
expansion in $\eta$. A more explicit expression shall be given at the very
end of section 3, after we have completed the vertical cluster expansion.

%%%%%%%%%%%%%%%%%%%%%%%%%%%%%

\subsection{Vertical cluster or momentum-decoupling expansion}  \label{subsection:vert-cluster}

%%%%%%%%%%%%%%%%%%%%%%%%%ש

{\em After} performing the scale $0$ horizontal cluster expansion, one must still perform on the contribution associated 
to a given forest $\F^0$  another expansion called
{\em vertical cluster} or {\em momentum-decoupling expansion}. 
This consists simply in applying the operator

\BEQ {\mathrm{Vert}}^0=  \prod_{\Del\in\D^0} \left( \sum_{\mu_{\Del}=0}^{2} \partial^{\mu_{\Del}}_{\tau_{\Del}}\big|_{\tau_{\Del}=0} +  \int_0^1  d\tau_{\Del}
\frac{(1-\tau_{\Del})^{2}}{
2!} \partial_{\tau_{\Del}}^{3} \right).  \label{eq:Vert0} \EEQ

Fix  a box $\Del^0\in\D^0$. A derivative $\partial_{\tau^0}$, acting on a dressed field $A(\tau;\cdot)$, simply beheads $A^0$ -- the highest-momentum component of $A$
--, and yields $A^{\to 1}$. On the other hand, if $(t,x)\in\Del^0$,
\BEQ
\partial_{\tau_{\Del^0}} R_{\eta}(\tau)=R_{\eta}(\tau) \left(\partial_{\tau_{\Del^0}} V_{\eta}(\tau)\right)R_{\eta}(\tau) 
\EEQ
\BEA && \partial_{\tau_{\Del^0}} V_{\eta}(\tau;t,x)= \nonumber\\
&&  B^{\to 1}(\cdot,(t,x))
 \left( g^{(0)}  (\eta(t,x)-v^{(0)}) \right) A(\tau;(t,x),\cdot) +   B(\tau;\cdot,(t,x)) \left( g^{(0)}(\eta(t,x)-v^{(0)}) \right) A^{\to 1}((t,x),\cdot)  \nonumber\\
&& + B^{\to 1}(\cdot,(t,x)) \left(-2\tau^0_{t,x} (\nu^{(0)}-\nu_{eff})\Del^{\to 0}
\right) A^{\to 1}((t,x),\cdot)
\EEA
Therefore, {\em the vertical cluster expansion acts by inserting vertices, just as
the horizontal cluster expansion does}.  On the other hand, these vertices comprise
at least one low-momentum field. {\em $0$-scale boxes in which these low-momentum fields are integrated}
(here $\Del^0$)  constitute the {\bf external boxes} or (looking
more precisely at the nature -- $A$ or $B$ - and the scale of the
low-momentum fields) the 
{\bf external structure} of the associated polymers. Such
low-momentum fields are called {\bf external legs} of the polymer. The
order of differentiation in $\tau_{\Del}$ is denoted by $\mu_{\Del}$;  for the
Taylor remainder in (\ref{eq:Vert0}) one has $\mu_{\Del}=3$.  Since
each $\tau$-derivative contributes an external leg, the number
of external legs of a polymer is equal to the number of $\tau$-derivatives that have been applied to it.  Thus $\mu_{\Del}$ can be interpreted as a {\em multiplicity}, by which we mean that
a polymer containing $\Del$ has $\mu_{\Del}$ external legs starting
from the box $\Del$.

\bigskip\noindent Now that we have completed the cluster expansion, {\em a fundamental observation to be made is the following.}  Let $\Del:=\Del^0_{t,x},\Del':=\Del^0_{t',x'}$. If $\Del,\Del'$ belong to different
 components of $\F^0$, then $R_{\eta}(\tau)(\vec{s}(\vec{w}))((t,x),(t',x'))=0$.
 In the contrary case, letting $\T^0$ be the tree containing $\Del$ and $\Del'$, 
 {\em  $R_{\eta}(\tau)(\vec{s}(\vec{w}))((t,x),(t',x'))$ depends
 only on the values of $\eta$ in the image $|\T^0|:=\{\Del\in\D^0\ |\ \Del\in\T^0\}$
 of the polymer.}

\bigskip\noindent We illustrate the double horizontal/vertical cluster expansion by
Fig. 1, where the following pictural conventions are used. {\em 
Wavy lines} are pairings $\langle \eta(t,x)\eta(t',x')\rangle$ produced by the cluster expansion in $\eta$; the attached
$d/ds$ is a reminder of the action of the cluster operator
$d/ds$ which produced the pairing.
{\em Wavy half-lines} with added symbol $\tilde{\eta}$ stand for
dangling  $\tilde{\eta}$-ends; when evaluating averaged $N$-point
functions, they are contracted {\em inside} their connected component
(polymer).  {\em Scale 0 thick lines} are space-time convolutions
$A^0 R^{(0)}(\tau^0=0)B^0$; an attached $d/ds$ signals the fact
that either $A^0$ or $B^0$ has been produced by the propagator
cluster. {\em Scale $j$ thick lines} ($j\ge 1$) are either 
$A^j$ or $B^j$ or $G^j=A^j B^j$. 
%\vskip -10cm
\begin{figure}[H]  \label{fig:cluster}
  \centering
  %\vskip 3cm
  \caption{Cluster expansions.}
  \vskip -2cm
   \includegraphics[scale=0.75]{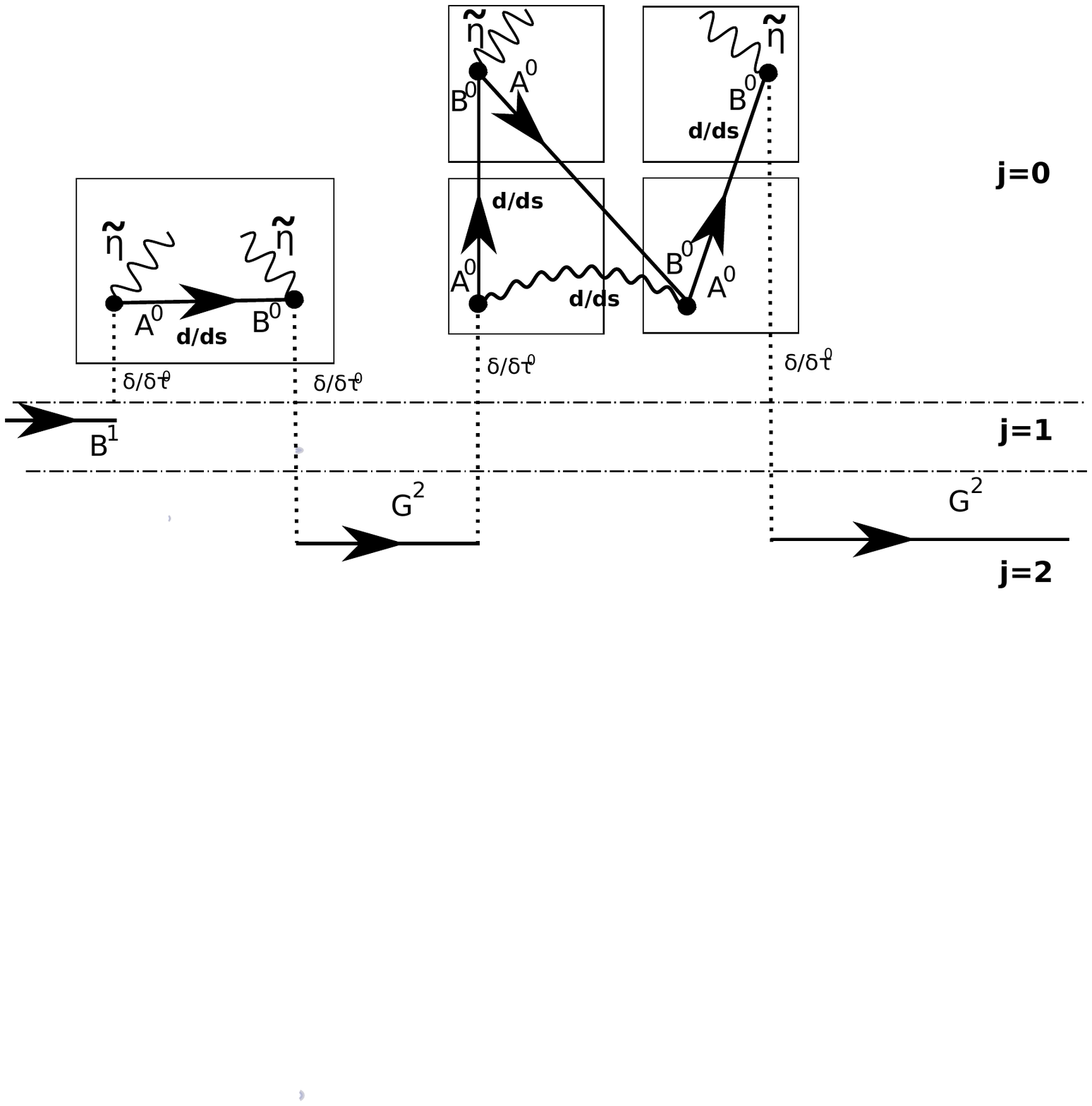}
\end{figure}

%%%%%%%%%%%%%%%%%%%%%%%%%%%%%%ש

\vskip -10cm\noindent
 The final outcome of this section is the following compact
expression, where $V(\F^0)$ is the set of vertices connected by a forest $\F^0$, and
$n_{\Del}$, $\Del\in V(\F^0)$ is the coordination number of a given vertex of the forest:
\BEA &&  \Big\langle \log(w_1(\tau^0,\vec{s}=1;t_1,x_1))\cdots
\log(w_N(\tau^0,\vec{s}=1;t_N,x_N)) \Big\rangle=\nonumber\\ 
 &&  \sum_{\F^0\in{\cal F}^0}  \left(
\prod_{\ell\in L(\F^0)}\int_0^1   dw_{\ell}\right)  
\sum_{L_G, L_{\eta},L_{vert},\vec{\mu}}   c_{\vec{I}} \ \  \Big(
\prod_{\ell\in L_{\eta}
} \int_{\Del_{\ell}}dz_{\ell}\ \int_{\Del'_{\ell}}dz'_{\ell}  \langle \eta(z_{\ell})\eta(z'_{\ell})\rangle_{\vec{s}=1} \Big)  \nonumber\\
&&  \ \Big\langle  
\frac{ \prod_{k\le N} \prod_{j\le m_k} \Big( \Big[  (D_G)_{k,j} (D_{\eta})_{k,j} (D_{\eta})_{(k,j),\cdot} (D_{\eta})_{\cdot,(k,j)}
(D_{\tau})_{k,j}  
  \Big]
w_k(\tau^0,\vec{s};\cdot)  \Big)}{\prod_{k\le N}
(w_k(\tau^0,\vec{s};\cdot))^{m_k}}    \Big\rangle_{\vec{s}(\vec{w})} , 
\nonumber\\
  \label{eq:appl-BKAR-F3} \EEA
where; $\vec{\mu}:=(\mu_{\Del})_{\Del\in\D^0}$;  $(D_{\tau})_{k,j}=\prod_{\Del\in L_{k,j}} D_{\tau_{\Del}}(\mu_{\Del})$, and
\BEQ D_{\tau_{\Del}}(\mu_{\Del})=\partial_{\tau_{\Del}}^{\mu_{\Del}} \Big|_{\tau_{\Del}=0}
\qquad (\mu_{\Del}=0,1,2),\qquad \int_0^1 d\tau_{\Del} \, 
\frac{(1-\tau_{\Del})^2}{2!} \partial^3_{\tau_{\Del}} \qquad
(\mu_{\Del}=3),\EEQ

with $\tilde{\eta}=g^{(0)}(\eta-v^{(0)})$, 
 featuring a product of strings indexed by $k,j$,
where, for each box $\Del\in\F^0$: 
\begin{itemize}
\item[(i)] the horizontal cluster expansion has produced $0\le n'_{\Del}\le n_{\Del}$
vertices integrated over $z'_n(\Del)\in\Del$, $n=0,\ldots,n'_{\Del}$;
\item[(ii)] the vertical cluster expansion has produced $0\le n''_{\Del}\le \mu_{\Del}
\le 3$   vertices integrated over $z''_n(\Del)\in\Del$, $n=0,\ldots,n'''_{\Del}$;
\end{itemize}
and $\{z'_n(\Del),z''_n(\Del)\}_{\Del,n} = \{z_{k,j}^i\}$.

%%%%%%%%%%%%%%%%%%%%%%%%%%שש

%%%%%%%%%%%%%%%%%%%%%%%%%%%%%%ש
%%%%%%%%%%%%%%%%%%%%%%%%%%ש

\section{Renormalization}  \label{sec:renormalization}

%%%%%%%%%%%%%%%%%%%%%%%%%%%%שש
%%%%%%%%%%%%%%%%%%%%%%%%%%%%%%%%

We now proceed to the renormalization stage. As explained in the introduction
to section 3, renormalization consists in general in computing, and compensating
by equal counterterms, the "diverging part" of the sum of 
diagrams with a given external structure. In a multi-scale setting,
one considers instead the  so-called "local part" of the {\em sum of all polymers} with internal legs of scale $\le j$
and given external structure, made up of a product of  external legs of scale $>j$;
such local parts are compensated by counterterms of scale $j$.

\medskip\noindent Given the simplicity of this stage in the present model, we spare the reader
a full-length explanation of these ideas (that can be found e.g. in \cite{MagUnt2}
or \cite{Unt-mode}), and describe instead what we do in
simple terms.

\medskip\noindent  The main step is the estimation of the {\em two-point function}. The idea is 
roughly the following. Low-momentum propagators $G^{\to 1}((t_i,x_i),(t_f,x_f))=\sum_{j\ge 1} A^j\langle j|\, B^j|j\rangle\, ((t_i,x_i),(t_f,x_f))$,
$j\ge 1$ occupying on a string the time-section between initial time $t_i$ and final
time $t_f$, may be cut anywhere into two parts by a scale $0$ vertex
insertion, according to the rule
\BEA && G^{\to 1}((t_i,x_i),(t_f,x_f))\rightsquigarrow \sum_{j,k\ge 1} A^j((t_i,x_i),\cdot)\langle j|\, \Big[ B^j|j\rangle  \Big( g^{(0)}\eta A^0 \langle 0|\,  B^0 |0\rangle \, g^{(0)}\eta +\cdots\Big)  A^k \langle k|\,  \Big](\cdot,\cdot) \ \cdot
\nonumber\\ &&\qquad \cdot\   B^k(\cdot,(t_f,x_f))\,  |k\rangle \ = G^{\to 1}((t,x),\cdot)\  K_{\eta}(\cdot,\cdot)\  G^{\to 1}(\cdot,(t',x')) \nonumber\\
\EEA

\medskip\noindent The random kernel between parentheses, 
\BEA
&&
K_{\eta}((t,x),(t',x')):=
\Big(g^{(0)}\eta A^0 \langle 0| \frac{1}{1- B^0 |0\rangle \, g^{(0)}\eta A^0 \langle 0| }  B^0 |0\rangle \, g^{(0)}  \eta\Big) ((t,x),(t',x'))  \nonumber\\
&&=
\Big( g^{(0)}\eta A^0 \langle 0|\,  B^0 |0\rangle \, g^{(0)}\eta +\cdots\Big)((t,x),(t',x'))
\end{eqnarray}

 containing only
$A^0$- and $B^0$-components, is (as can be shown) 
$O(1)$ in average, and decreases exponentially fast when $d((t,x),(t',x'))$ is large, while
\BEQ 
G^j((t',x'),\cdot)\simeq G^j((t,x),\cdot) \label{eq:Gjt'x'tx}
\EEQ if $d((t,x),(t',x'))=O(1)$ and 
$j\gg 1$. Thus it makes
sense to assume that its main contribution to the string is the averaged zero-momentum
quantity $v(t):=\int dt'\, dx'\, \langle K_{\eta}((t,x),(t',x'))\rangle$ (later on identified
as $g^{(0)}v^{(0)}$, up to some small correction).   Assuming for simplicity that 
$v(t)\equiv v$ is a constant, we must consider the sum of the geometric series
$G^{\to 1}+ G^{\to 1} v G^{\to 1} + G^{\to 1} vG^{\to 1}vG^{\to 1}\cdots$. Since
now $(G \ast G)((t,x),(t',x'))=\int_{t'}^t dt''\, \int dx''\,  p_{t-t''}(x-x'')p_{t''-t'}(x''-x')=(t-t')G((t,x),(t',x'))$, one sees that the large-scale (i.e. $t-t'\to
+\infty$) correction to $G^{\to 1}$ is infinite. On the other hand, the geometric
series may be resummed exactly, $G+GvG+GvGvG+\cdots= (\partial_t-\nu^{(0)}\Del-v)^{-1}$. This explains why we incorporated $v^{(0)}$ into the equation. Considering instead
a second-order Taylor expansion in $x-x'$ in (\ref{eq:Gjt'x'tx}) yields (see
 similarly (\ref{eq:delnu})) a
contribution $\del\nu$, compensated by $\del V$ (see (\ref{eq:dressed-vertex},
\ref{eq:del-V-eta})), creating a geometric series $\simeq G+G\del\nu\Del G+G\del\nu\Del G\del\nu\Del G+\cdots=(\partial_t-\nu^{(0)}\Del-
\del\nu \Del)^{-1}$; thus $\nu_{eff}:=\nu^{(0)}+\del\nu$ may be interpreted
as an effective viscosity. Now, further corrections, of the type $G^{\to 1}\rightsquigarrow G^{\to 1} \partial^{\kappa}G^{\to 1}$ with $|\kappa|\ge 3$, 
see our {\em first key power-counting estimate} (\ref{eq:PW1}),
 {\em finite} in the large-scale limit, need not be considered.

\medskip\noindent 
In a general renormalizable theory, only a finite number of  $N$-point functions yield infinite contributions in the large-scale
limit. It turns out here, however, that only $N=2$ point functions yield an infinite contribution,  because of our
{\em second key power-counting estimate} (\ref{eq:PW2}). We content ourselves with briefly discussing
diagrammatics for $N=4$ in \S \ref{subsection:four-point}.

%%%%%%%%%%%%%%%%%%%%%%%%%
%%%%%%%%%%%%%%%%%%%%%%%%ש

%%%%%%%%%%%%%%%%%%%%%%

\subsection{Two-point function}  \label{subsection:two-point}

%%%%%%%%%%%%%%%%%%%%%%%%ש

  Consider a {\em piece ${\cal S}$} of a string $A((t_{init},x_{init}),\cdot) (1-V_{\eta})^{-1}(\cdot,\cdot) B(\cdot,(0,y)) e^{\frac{\lambda}{\nu^{(0)}} h_0(y)}$ running from
 initial position $(t_{init},x_{init})$ to final position $(0,y)$, {\em connected by the horizontal cluster alone
 (i.e. obtained by letting $\tau^0\equiv 0$).} By construction, it has  two external legs, one at each temporal end. Then (letting
$\tilde{\eta}(t,x):=g^{(0)}(\eta(t,x)-v^{(0)})$ -- see Definition \ref{def:dressed-vertex} --,
and $L_{\eta }(\F^{0})$ be the set of cluster links coming from the perturbation of the measure on  $\eta $ -- compare with 
eq. (\ref{eq:appl-BKAR-F2}), while now $k=j=1$ since there is only
one string, and $I_{1,1}=L_{\eta}(\F^0)$ --)

\BEA 
&& {\cal S}:= \Big(\prod _{\ell \in L(\F^{0})} \int dw_{\ell} \Big) \sum_{L_G,L_{\eta}} \Big(\prod _{\ell \in L_{\eta }} \int_{\Del_{\ell}} dz_{\ell}\,  \int_{\Del'_{\ell}} dz'_{\ell}\,    
\langle  \eta (z_{\ell }) \eta (z'_{\ell }) \rangle_{\vec{s}=\vec{1}} \Big)\ \cdot\  \nonumber\\ &&
\  D_G D_{\eta} \Big\{ B^{\to 1}(\cdot,(t,x)) \frac{1}{1-\int dt"\, dx"\,  V^{(0)}_{\eta}(\tau^0=0)(\vec{s}(\vec{w}))(t",x")}((t,x),(t',x')) A^{\to 1}((t',x'),\cdot) \Big\} \nonumber\\
 &&=
B^{\to 1}(\cdot,(t,x)) \ \cdot\ \tilde{\eta}(t,x)\ \cdot\ A^{\to 1}((t,x),\cdot) \nonumber\\
&&\ \  +\sum_{n\ge 0} \int dt_1\, dx_1\cdots \int dt_n\, dx_n\, 
\Big(\prod _{\ell \in L(\F^{0})} \int dw_{\ell} \Big)\sum_{L_G,L_{\eta}}  \Big(\prod _{\ell \in L_{\eta}}  \int_{\Del_{\ell}} dz_{\ell } \int_{\Del'_{\ell}} dz'_{\ell}    
\, \langle  \eta (z_{\ell }) \eta (z'_{\ell }) \rangle  \Big)
\nonumber\\
&&\qquad D_G D_{\eta} \Big\{  B^{\to 1}(\cdot, (t,x)) \ \cdot\ \Big[ \tilde{\eta}(t,x) A^0(\vec{s}(\vec{w}))((t,x),\cdot) \langle 0| \ \cdot\ R_{\eta}^{(0)}(\tau^0=0)(\vec{s}(\vec{w}))(\cdot,\cdot) \ \cdot\  \nonumber\\
&&\qquad \big( B^0(\vec{s}(\vec{w}))(\cdot,(t_1,x_1)) |0\rangle \,  \tilde{\eta}(t_1,x_1))A^0(\vec{s}(\vec{w}))((t_1,x_1),\cdot) \langle 0| \big)
\ \cdot\ R_{\eta}^{(0)}(\tau^0=0)(\cdot,\cdot) \nonumber\\
&&\qquad \cdots\   
\big( B^0(\vec{s}(\vec{w}))(\cdot,(t_i,x_i)) |0\rangle\,   \tilde{\eta}(t_i,x_i))A^0(\vec{s}(\vec{w}))((t_i,x_i),\cdot) \langle 0| \big)\ \cdot\ 
 R_{\eta}^{(0)}(\tau^0=0)(\cdot,\cdot) \nonumber\\
&&\qquad \cdots\   
\big( B^0(\vec{s}(\vec{w}))(\cdot,(t_n,x_n)) |0\rangle\,   \tilde{\eta}(t_n,x_n))A^0(\vec{s}(\vec{w}))((t_n,x_n),\cdot) \langle 0| \big)\ \cdot\ 
 R_{\eta}^{(0)}(\tau^0=0)(\cdot,\cdot) \nonumber\\
&&\qquad B^0(\vec{s}(\vec{w}))(\cdot,(t',x')) |0\rangle \, 
\tilde{\eta}(t',x') \Big]  \cdot\  A^{\to 1}((t',x'),\cdot) \Big\} \label{eq:two-point1}
 \EEA

where $n$ is the number of internal vertices, and $R_{\eta}^{(0)}(\tau^0=0)$ 
are "scale $0$ resolvents",
\BEQ R^{(0)}_{\eta}(\tau^0=0)(\vec{s}(\vec{w}))((t,x),(t',x'))=\Big(1-\int dt\, dx\,  V^{(0)}_{\eta}(\tau^0=0)(\vec{s}(\vec{w}))(t,x)\Big)^{-1}.\EEQ
The operators $D_G$, $D_{\eta}$ are as in (\ref{eq:appl-BKAR-F3}),
with only one string involved, say, $D_G\equiv (D_G)_{1,1}, D_{\eta}\equiv
(D_{\eta})_{1,1}$. Each $\partial/\partial s_{\ell}$ appearing in $D_G$ suppresses
one of the $s$-factors in front of the propagators; each 
$\frac{\del}{\del \eta(z_{\ell})\del\eta(z'_{\ell})}$ appearing
in $D_{\eta}$ takes out the corresponding pair of $\tilde{\eta}$'s.
How this is done is specified by the choice of $(\F^0,L_G,L_{\eta})$. Thus the action of $D_G,D_{\eta}$ is extremely simple and produces no extra combinatorial factors.

\medskip\noindent It is convenient to describe the lonely term in the first line of (\ref{eq:two-point1}), 
obtained simply by differentiating twice with respect to  $\frac{d}{d\tau^0_{\Del^0}}$ in a box $\Del^0=\Del^0_{t,x}$ untouched by the horizontal cluster,  as an "$n=-1$" contribution; note that it contains implicitly a Dirac function $\del((t,x),(t',x'))$.

\medskip\noindent {\em Assume that the
$\eta$'s inside the brackets $\Big[ \ \cdot\ \Big]$ contract pairwise}, or
equivalently, that no $\eta$-field on $\cal S$ pairs to an $\eta$-field on another
string. 
By the first property below  Definition \ref{def:Aj}, namely,   since $\langle \eta(t_i,x_i)\eta(t_{i'},x_{i'})\rangle=0$ if $(t_i,x_i)$ is connected to $(t_{i'},x_{i'})$ by
some low-momentum propagator $A^{\to 1}$ or $B^{\to 1}$, {\bf only scale $0$ diagrams  contribute}; which explains why {\em we need not consider} generalizations of (\ref{eq:two-point1}) with brackets $\Big[ \ \cdot \ \Big]$ including {\em lower-momentum} $A$'s and
$B$'s. Note that, since  $R_{\eta}^{(0)}(\tau^0=0)=\Id+B^0(\vec{s}(\vec{w}))(\cdot,\cdot) |0\rangle\,  \tilde{\eta}(\cdot,\cdot) A^0  \langle 0| +\cdots$, other choices of
external legs are not allowed, for instance,
\BEQ \Big[\ \cdot\ \Big]\ \cdot\ B^{\to 1}= \Big[ \cdots A^0 \langle 0| \Big] 
\Big(B^1 |1\rangle + B^2 |2\rangle+\cdots\Big)\equiv 0   \label{eq:BAB} \EEQ
because the basis $(|j\rangle)_{j\ge 0}$ is orthonormal.

\bigskip\noindent Let $\Sigma_0((t,x),(t',x'))$ be the average with respect to
the measure in $\eta$ of the sum of all contributions like
the one in $\left[ \ \cdot\ \right]$ in (\ref{eq:two-point1}); 
the kernel $\Sigma_0((t,x),(t',x'))$ must be seen as a deterministic insertion 
on the string between $(t,x)$ and $(t',x')$. {\em For  reasons explained in
{\bf C.} below, we
symmetrize the kernel $\Sigma_0$ by 
letting $\Sigma_0((t',x'),(t,x)):=\Sigma_0((t,x),(t',x'))$ if $t'<t$.}  We split the discussion into a number 
of steps.

\bigskip\noindent {\bf A.} {\em A first step} consists in {\em displacing the final external leg
$A^{\to 1}((t',x'),\cdot)$ to the location $(t,x)$  of the initial external
leg $B^{\to 1}(\cdot,(t,x))$} (or conversely, see below). Namely,
\begin{eqnarray}\label{loc1}
&&  B^{\to 1}(.(t,x))\Sigma_0((t,x),(t',x'))A^{\to 1}((t',x'),.)  =
  B^{\to 1}(.(t,x)) A^{\to 1}((t,x),.)\Sigma_0((t,x),(t',x')) \nonumber\\
  &&\qquad +B^{\to 1}(.(t,x))\Sigma_0((t,x),(t',x'))\  \big[ A^{\to 1}((t',x'),.)-A^{\to 1}((t,x),.)\big].
\end{eqnarray}
 Then we Taylor expand  $ A^{\to 1}((t',x'),.)-A^{\to 1}((t,x),.)  $   to parabolic order three:
  \begin{eqnarray}\label{loc2}
&& A^{\to 1}((t',x'),.)-A^{\to 1}((t,x),.) =\nonumber\\
 &&\qquad =\Big((t'-t)\partial _{t}+(x'-x)\cdot \nabla_{x} 
  +\half \sum _{i,j}(x'-x)_{i} (x'-x)_{j}\partial_{x_i}\partial_{x_j}
  \Big)A^{\to 1}((t,x),.) \nonumber\\
  &&\qquad +\int _{0}^{1} du\, \frac{(1-u)^{2}}{2} \frac{d^{3}}{du^{3}} \big\{ A^{\to 1}(((1-u^2)t+u^2 t',(1-u)x+ux'),.)  \big\}  
\end{eqnarray}
See Fig. 2 for an illustration.  

\begin{figure}[htpb]  \label{fig:local-part}
  \centering
  \caption{Displacement of external legs.}
  \vskip 1cm
   \includegraphics[scale=1.3]{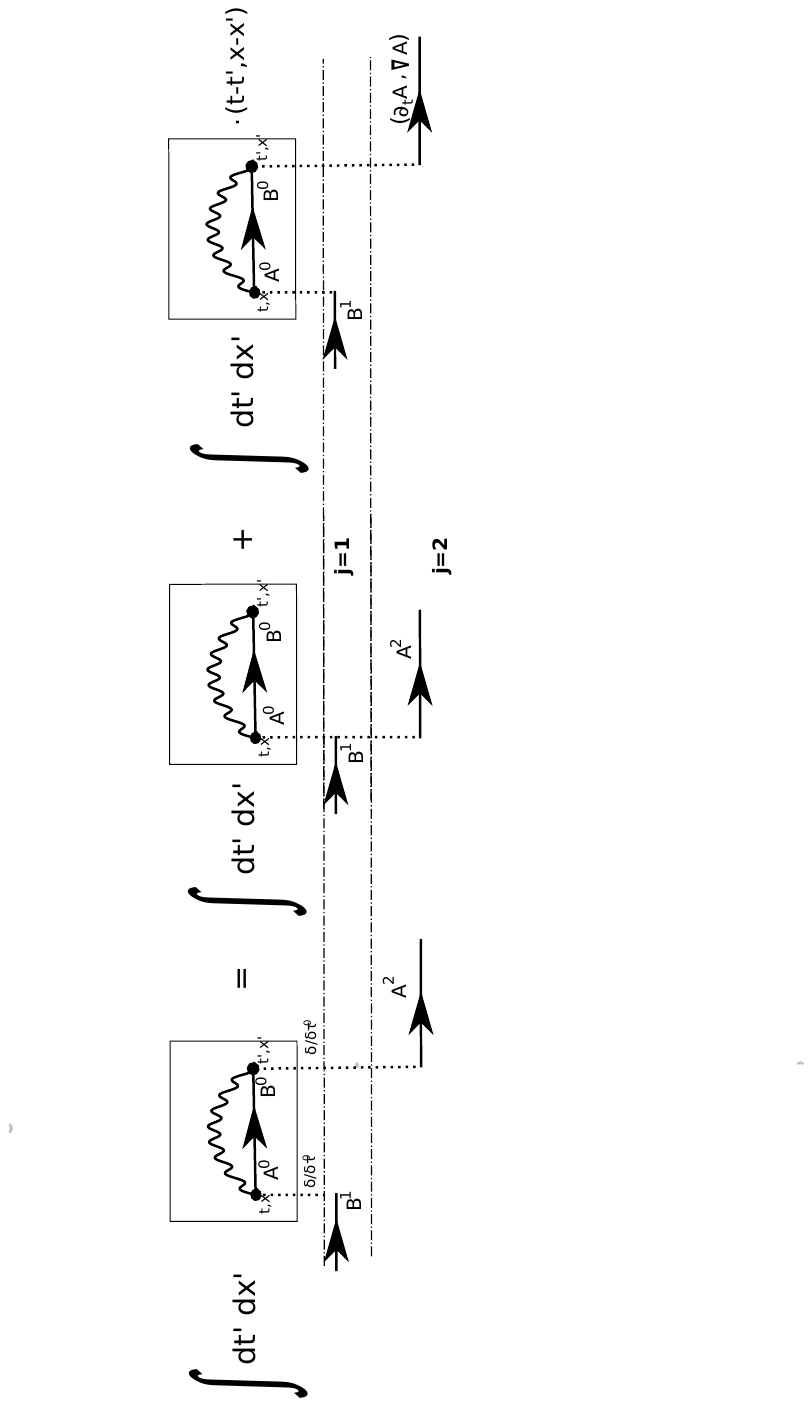}
\end{figure}

 The integral remainder term in (\ref{loc2})  is a sum of  derivatives of
parabolic order $\ge 3$  (more precisely, ranging in $\{3,\ldots,6\}$), 
 \begin{eqnarray}
\label{reg2}
&& \left|\frac{d^{3}}{du^{3}} A^{\to 1}(\cdot,\cdot)\right| \lesssim 
|x-x'|^3 |\nabla^3 A^{\to 1}(\cdot,\cdot)|+ (t-t')  |x-x'|\  |\partial _{t}\nabla A^{\to 1}(\cdot,\cdot)|  \nonumber\\
&& \qquad +
  (t-t')^{2} |\partial _{t}^{2}  A^{\to 1}(\cdot,\cdot)|+   (t-t')^{3} \, 
| \partial_t^3  A^{\to 1}(\cdot,\cdot)| + (t-t')^2 |x-x'| \ |\partial^2_t\nabla A^{\to 1}
(\cdot,\cdot)|  \nonumber\\
&&\qquad\qquad\qquad\qquad\qquad\qquad\qquad
+(t-t') |x-x'|^{2} |\partial _{t}\nabla^{2} A^{\to 1}(\cdot,\cdot)|   \label{eq:4.8}
\end{eqnarray}

\medskip\noindent The main terms in (\ref{reg2}) are those on the first line; splitting $A^{\to 1}$
into its constituent scales $\sum_{j'\ge 1} A^{j'}\langle j'|$, we known from section 1 that
$\nabla^3 A^{j'}, \partial_t\nabla A^{j'}\sim 2^{-3j'/2} A^{j'}$, whereas $|x-x'|^n |t-t'|^m \Sigma_0((t,x),(t',x'))=O(\Sigma_0((t,x),(t',x')))$ for all $n,m\ge 0$ (due
to the exponential decrease in $d((t,x),(t',x'))$, see below), all together a gain of  $O(2^{-3j'/2})$.

\medskip\noindent
The other terms are dealt with below, namely the first term in the r.-h.s. of 
(\ref{loc1}) and the first line in the r.-h.s. of (\ref{loc2}) -- more
precisely, only the second-order, traced term $\frac{1}{6} |x'-x|^2 \Del_x
A^{\to 1}((t,x),\cdot)$, the other ones vanishing by symmetry --; they contribute to the renormalization of the two-point function.

\medskip\noindent In order to get the smaller of  two factors, we displace instead
the initial external leg $B^{\cdot}(\cdot,(t,x))$ to the location of the final
external leg $A^{\cdot}((t',x'),\cdot)$ if the scale $j$ of the $B$-leg is strictly lower
than the scale $j'$ of the $A$-leg, i.e. if $j> j'$, yielding a small factor
$O(2^{-3j/2})$. 
\medskip\noindent
Summarizing: were it not for (i) the {\em boundary conditions at initial time $t_{init}$
and final time $0$}, and (ii) the  {\em non-overlapping condition between the scale $0$ boxes
chosen by the horizontal cluster expansion},  the contribution would be  (taking
into account the symmetrization of the kernel $\Sigma_0$, and considering -- as an
intermediate step only -- the natural extension of the model to negative times)
\BEA &&  \half\sum_{1\le j\le j'} \Big( \Big\{\int_{-\infty}^{+\infty} dt'\, dx'\,  \Sigma_0((t,x),(t',x')) \Big\}\,  B^{j}(\cdot,(t,x))A^{j'}((t,x),\cdot) 
 \nonumber\\
&& + \Big\{\frac{1}{2d} \int_{-\infty}^{+\infty} dt'\, dx'\, |x-x'|^2 \Sigma_0((t,x),(t',x')) 
\Big\} \, B^{j}(\cdot,(t,x)) \Del A^{j'}((t,x),\cdot) \, \Big),  \label{eq:4.9}
\EEA
 plus the same expression up to the exchange of $A$, $B$ and $(t,x),(t',x')$
summed over  $j>j'\ge 1$. Choosing $v^{(0)}$ such that $\int dt'\, dx'\, 
\Sigma_0((t,x),(t',x'))=0$, and letting  
\BEQ \del\nu:=\frac{1}{4d} \int_{-\infty}^{+\infty} dt'\, dx'\, |x-x'|^2\,  \Sigma_0((t,x),(t',x')), \label{eq:delnu0}
\EEQ
 this is 
   equivalent to the addition to the vertex $V_{\eta}$   of\\ $
  \int dt\, dx\, \half\frac{d^2}{d\tau^0_{t,x}} B^{\to 1}(\cdot,(t,x)) ((\tau^0_{t,x})^2\del\nu \Del) A^{\to 1}((t,x),\cdot)$, compensating the term proportional to 
  $(\tau^0_{t,x})^2$ in $\del V_{\eta}$, see (\ref{eq:dressed-vertex}). 

\medskip\noindent Let us consider objections (i) and (ii) separately. First, because of the boundary conditions,  the integral $\half \int_{-\infty}^{+\infty} dt'\, dx' \, (\cdots)=\int_{-\infty}^t dt'\, (\cdots) $ in (4.9) must be replaced by $\int_0^t dt'\, (\cdots)$.
Similarly, if $j>j'$, the integral $\half \int_{-\infty}^{+\infty} dt\, dx\, (\cdots)=\int_{t'}^{+\infty} dt\, (\cdots)$ must be replaced by $\int_{t'}^{t_{init}} \, (\cdots)$.
Differences $\big(\int_{-\infty}^t - \int_0^t\big) dt'\, (\cdots)=\int_{-\infty}^0
dt'\, (\cdots)$, resp. $\big(\int_{t'}^{+\infty}-\int_{t'}^{t_{init}}\big) dt\, (\cdots)=\int_{t_{init}}^{+\infty} dt\, (\cdots)$, are shown in {\bf D.} to be
exponentially small in the distance to the boundary, $t-0$, resp. $t_{init}-t'$. Thus
one may equivalently define $\del\nu$  by an integral over {\em positive times}, which is more natural given that we are considering an initial-value problem,
\BEQ \del\nu:=\frac{1}{2d} \lim_{t\to +\infty}  \int_{0}^t dt'\, dx'\, |x-x'|^2\,  \Sigma_0((t,x),(t',x')). \label{eq:delnu}
\EEQ

\medskip\noindent Next, due to the non-overlapping condition, the factorization  of $\Sigma_0$ fails. The solution to this well-known problem is 
through a {\em Mayer expansion}.

\bigskip\noindent {\bf B (Mayer expansion).} Namely, we shall now   apply
the restricted cluster expansion, see Proposition \ref{prop:BK2}, to
the result of our expansion.  Cluster expansions have produced a scale 0 forest 
$\F^0$ of boxes, whose  tree components, together with their
external structure made up of low-momentum $A$'s and $B$'s, are called
{\em polymers}, and denoted by $\P_1,\ldots,\P_N$.  The {\em objects} are now
scale $0$ polymers $\P$ in    ${\cal O}=\{\P_1,\ldots,\P_N\}$ ; a link
$\ell\in L({\cal O})$ is a pair of polymers $\{\P_n,\P_{n'}\}$, $n\not=n'$.
Objects of type 2 are  polymers with $>2$ external legs, whose non-overlap conditions we shall not remove at this stage, because these polymers are already convergent, hence
do not need to be renormalized.
Then objects of type $1$  are   polymers with two external legs; note that -- due to the displacement of externel legs operated in
 {\bf A.} --
the two external legs are located in the {\em same} scale $0$ box.

\medskip\noindent Implicit in the outcome of the cluster expansions is  the {\bf non-overlapping condition},
 \BEA {\mathrm{NonOverlap}}(\P_1,\ldots,\P_N) &:=& \prod_{(\P_n,\P_{n'}) } {\bf 1}_{\P_n,\P_{n'}\ {\mathrm{non}}-{\mathrm{overlapping}}} \nonumber\\
 &=&
\prod_{(\P_n,\P_{n'}) } \prod_{\Del\in\vec{\Del}(\P_n),\Del'\in
\vec{\Del}(\P_{n'})}
 \left( 1 + \left(
 {\bf 1}_{\Del\not=\Del'}-1 \right)\right) \nonumber\\ \EEA
 stating that a box $\Del$ belonging to the image of $\P_n$ and
 a box $\Del'$ belonging to the image of $\P_{n'}$ are necessarily
 distinct.
As in the proof of BKAR formula (see Proposition \ref{prop:BKAR}), we
choose some polymer, say $\P_1$, with $2$ external legs, and 
weaken the non-overlap condition between $\P_1$ and all the other polymers by introducing a parameter $S_1$,

 \BEA && {\mathrm{NonOverlap}}(\P_1,\ldots,\P_N) (S_1)=
  \prod_{\{\P_n,\P_{n'}\}_{n,n'\not=1} } \prod_{\Del\in{\bf\Del}(\P_n),\Del'\in
{\bf\Del}(\P_{n'})} {\bf 1}_{\Del\not=\Del'} \ \cdot\nonumber\\
 &&\qquad 
\prod_{ (\Del,\Del')\in\vec{\Del}(\P_1)\times \vec{\Del}(\P_{n'})\setminus \vec{\Del}_{ext}(\P_1)\times \vec{\Del}_{ext}(\P_{n'}) } \left( 1 +  S_1 \left( {\bf 1}_{\Del\not=\Del'}-1 \right)\right),  \nonumber\\ \label{eq:Mayer1}\\ \EEA
where $\vec{\Del}_{ext}(\P)\subset\vec{\Del}(\P)$ is the subset of boxes $\Del$ with external legs - i.e. that have been differentiation
with respect to $\tau_{\Del}$ -, 
 and Taylor expand in $S_1$ to order $1$; each factor 
 \BEQ {\bf 1}_{\Del\not=\Del'}-1 =-{\bf 1}_{\Del=\Del'} \EEQ
  produced by
 differentiation is a {\bf Mayer link} between $\P_1$ and
 some $\P_{n'}, n'\not=1$, or more precisely, some  box $\Del\in\vec{\Del}(\P_1)$ and some
 box $\Del'\in\vec{\Del}(\P_{n'})$, implying an explicit overlap between $\P_1$ and
 $\P_{n'}$, and adding a link to the forest $\F^0$.
 Iterating the procedure and applying Proposition \ref{prop:BK2}
 to the weakened non-overlap condition
 \BEA &&{\mathrm{NonOverlap}}(\P_1,\ldots,\P_N) (\vec{S}):
=  \prod_{\{\P_n,\P_{n'}\} } \prod_{\Del\in{\bf\Del}_{ext}(\P_n),\Del'\in
{\bf\Del}_{ext}(\P_{n'})} {\bf 1}_{\Del\not=\Del'} \ \cdot\nonumber\\
&&
 \qquad 
\prod_{(\Del,\Del')\in\vec{\Del}(\P_n)\times\vec{\Del}(\P_{n'})\setminus
{\bf\Del}_{ext}(\P_n)\times
{\bf\Del}_{ext}(\P_{n'})} \left( 1 +  S_{\P_n ,\P_{n'}} \left( {\bf 1}_{\Del\not=\Del'}-1 \right)\right),  \nonumber\\ \label{eq:Mayer2}
\EEA
 
 The outcome is a sum 
 \BEA &&
 \sum_{\G^0\in {\cal F}_{res}({\cal O})} \Big( \prod_{\ell \in 
 L(\G^0)} \int_0^1 dW_{\ell}\Big)\ \   {\mathrm{NonOverlap}}(\vec{S}(\vec{W})), \nonumber\\
 &&\qquad  {\mathrm{NonOverlap}}(\vec{S}(\vec{W})):= \Big[ \Big(\prod_{\ell\in L(\G^0)} \frac{\partial}{\partial S_{\ell}} \Big) {\mathrm{ NonOverlap}}(\P_1,\ldots,\P_N)\Big] (\vec{S}(\vec{W}))  \label{eq:Mayer3}
\nonumber\\
 \EEA
 Links $\ell=\ell_{\P_n,\P_{n'}}\in L(\G^0)$ are obtained as links between {\em polymers}, however the corresponding differentiation 
$\frac{\partial}{\partial S_{\ell}}$ is immediately rewritten as
a sum over pairs over boxes $(\Del,\Del')\in\vec{\Del}(\P_n)\times 
\vec{\Del}(\P_{n'})$. Thus we see Mayer links as {\em links between boxes}.  As such they add up to the set of links $L(\F^0)$ produced
by the horizontal cluster expansion, producing a forest $\bar{\F}^0$
with same vertices as $\F^0$ but larger set of links $L(\bar{\F}^0)\equiv 
L(\F^0)\uplus L_{{\mathrm{Mayer}}}$, where $L_{{\mathrm{Mayer}}}$
(in bijection with $L(\G^0)$ is the set of Mayer links. Since
a forest is characterized by its set of links, we  rewrite in practice 
(\ref{eq:Mayer3}) as
\BEA && \sum_{L_{{\mathrm{Mayer}}}} \Big( \prod_{\ell \in 
 L_{{\mathrm{Mayer}}}} \int_0^1 dW_{\ell}\Big)\ \   {\mathrm{Mayer}}(\vec{S}(\vec{W})), \nonumber\\
 &&\qquad  {\mathrm{Mayer}}(\vec{S}(\vec{W})):= \Big[ \Big(\prod_{\ell\in L_{{\mathrm{Mayer}}}} \frac{\partial}{\partial S_{\ell}} \Big) {\mathrm{ NonOverlap}}(\P_1,\ldots,\P_N)\Big] (\vec{S}(\vec{W}))  \label{eq:Mayer4}.
 \EEA

\medskip\noindent The number of external legs of a set of polymers connected by Mayer links is the
sum of the number of external legs of each of the polymers. In particular, any
Mayer connected component containing at least two polymers has $\ge 4$ external legs;
it has become convergent.

\medskip\noindent   Let us now give some necessary precisions. Since the Mayer expansion is really applied to the non-overlap function NonOverlap and {\em not} to
the outcome of the expansion,  one must
still extend the outcome of the expansion  to the case when the $\P_n$, $n=1,\ldots,N$ have some overlap. The natural way to do this is to assume that
 the random variables $(\eta\big|_{\P_n})_{n=1,\ldots,N}$ remain independent even when they overlap. This may be understood in the following way. Choose a different color for each polymer $\P_n=\P_1,\ldots,\P_N$, and paint with that color
 {\em all} intervals $\Del\in \P_n\cap\D^0$. If $\Del\in{\bf \Del}_{ext}(\P_n)$, then its
 external links to the  $A^{\to 1}, B^{\to 1}$ below it are left in black. The
 previous discussion implies that intervals with different colors may superpose; on the other hand, {\em external inclusion links} may {\em not}, so that
 {\em low-momentum fields} $B^{\to1}((\cdot),(t,x)), A^{\to1}((t,x),\cdot)$, $(t,x)\in\Del^{0}$ with $\Del^{0}\in{\bf\Del}_{ext}(\P_n)$,
do not superpose and may be left in black.

\medskip\noindent Hence one must see $\eta$ as living on a two-dimensional set, $\D^{0}\times\{ {\mathrm{colors}}\}$, so that copies of $\eta$ with different colors are independent of each other. This defines a new, {\em extended and restricted to the zeroth scale} resolvent $\tilde{R}_{\eta}^{(0)}(\tau^0=0)$ 
associated to an {\em extended} field $\eta :\R_+\times\R^d\times \{{\mathrm{colors}}\}\to\R$, and {\em Mayer-extended polymers}.
 By abuse of notation, we shall skip the tilde in the sequel, and always implicitly extend the fields and the measures of  scale $0$ by taking into account colors.

\bigskip\noindent {\bf C (counterterms).} We now {\em define $\Sigma((t,x),(t',x'))$} to be the Mayerization of the sum of all contributions like
the one in $\left[ \ \cdot\ \right]$ in (\ref{eq:two-point1}), in which the two external legs have been displaced into the same box as in {\bf A.}, so that
there is no non-overlapping restriction on the support but for the box containing $(t,x)$. 
Note that  Mayer links between 
polymers with two external legs produce Mayer polymers with $\ge 4$ external legs, which
are therefore convergent (see \S 4.2).

\medskip\noindent Then (provided
that the limit does exist)
\BEQ g^{(0)}v^{(0)}:=\lim_{T\to +\infty} \ \int_0^T dt' \, \int  dx'\,  \Sigma((T,x),(t',x')) =\lim_{T\to +\infty}
\int^T_{t'} dt\, \int dx\, \Sigma((t,x),(t',x')). \label{eq:g0v0} \EEQ
The result does not depend on $x$. Furthermore, as shown below, letting
\BEQ g^{(0)}v^{(0)}(T):= \int_0^T dt' \, \int dx'\,  \Sigma((T,x),(t',x'))=\int_{t'}^{T+t'} dt \, \int dx\, \Sigma((t,x),(t',x')),
\EEQ
with $T=t$, resp. $t_{init}-t'$, the {\em boundary correction $\del v^{(0)}(T)$
to $v^{(0)}$ decreases exponentially with $T$}, namely,
\BEQ \del v^{(0)}(T):= v^{(0)}-v^{(0)}(T)=O((Cg^{(0)})^{cT})\to_{T\to +\infty} 0
\label{eq:v0-v0m}
 \EEQ
for some constants $C,c>0$.

\bigskip\noindent
Consider once again the first term in the r.-h.s. of (\ref{loc1}) and the first line in the
r.-h.s. of (\ref{loc2}), but this time {\em after} the Mayer expansion; summing, we get
if $j'\ge j$
(with a factor $\half$ due to the symmetrization of $\Sigma_0$)
\BEA \label{eq:Sigma1}
&&  \half \int dt'\, dx'\, B^{j}(\cdot, (t,x)) \Sigma((t,x),(t',x')) \nonumber\\
&&\qquad
 \Big(1\ +\ (t'-t)\partial _{t}+(x'-x)\cdot \nabla _{x}
  +\half \sum _{i,j}(x'-x)_{i} (x'-x)_{j}\partial_{x_i}\partial_{x_j}
  \Big)
 A^{j'}((t,x),\cdot) \nonumber\\
\EEA

  The first term in  (\ref{eq:Sigma1}) vanishes for an adequate choice of $v^{(0)}$,
  as shown below. Then  the second (thanks to the symmetrization) and third terms vanish by parity, and 
 the fourth one  vanishes for $i \ne j $ by isotropy.  The remaining term in
 (\ref{eq:Sigma1}) may be  absorbed into  a redefinition  of $\nu$. Namely,
 we define for any $i=1,\ldots,d$, 
\BEQ \nu_{eff}-\nu^{(0)}:=\frac{1}{4}\int dt' \, dx' (x'_i-x_i)^2 
\Sigma((t,x),(t',x')).  \label{eq:delnu1/4} \EEQ
Thus 

\BEA &&
 \half \int dt'\, dx'\, B^{j}(\cdot, (t,x)) \Sigma((t,x),(t',x'))
A^{j'}((t',x'),\cdot) = v^{(0)} B^{j}(\cdot,(t,x)) A^{j'}((t,x),\cdot)  \nonumber\\
&&\ + (\nu_{eff}-\nu^{(0)})  B^{j}(\cdot,(t,x)) \Del^{\to 0}_x A^{j'}((t,x),\cdot) +  \mbox{remainders}
\EEA

Remainders include the previously discussed integral remainder term in (\ref{loc2}),
and the cut-off difference
\BEQ (\nu_{eff}-\nu^{(0)}) B^{j'\to}(\cdot,(t,x)) (\Del_x-\Del^{\to 0}_x)A^{j'}((t,x),\cdot),  \label{cut-off-diff} \EEQ
which  is bounded in absolute value by $O(|\nu_{eff}-\nu^{(0)}|\, B^{j'\to}(\cdot,(t,x)))$ times
\BEQ  \int dx'\, \bar{\chi}^0(x') |\nabla^2 A^{j'}((t,x),\cdot)-\nabla^2 A^{j'}((t,x+x'),\cdot)| 
\sim 2^{-3j/2} A^{j'}((t,x),\cdot),\EEQ
of the same order as the integral remainder term.

\bigskip\noindent The leading-order contribution in the coupling constant of 
$\nu_{eff}-\nu^{(0)}$ is obtained (as seen from (\ref{eq:delnu1/4}), letting $(t',x')=(0,0)$ and integrating in $(t,x)$ instead) by contracting the $\eta$'s in the expression
\BEQ \half (g^{(0)})^2 \int_0^{\infty} dt \, \int dx\, x_1^2\,  \eta(t,x)(A^0 B^0)((t,x),(0,0))
\eta(t',x'). \label{eq:leading-nu-eff0}
 \EEQ
 This is the $n=2$ term in   (\ref{eq:two-point1}) with $R_{\eta}^{(0)}(\cdot,\cdot)$ substituted by its leading order term $\del(\cdot\, -_, \cdot)$. By (\ref{eq:cov-eta}), one gets
\BEQ \nu_{eff}-\nu^{(0)}=\half \frac{\lambda^2 D^{(0)}}{(\nu^{(0)})^2} \int_0^{\infty} dt\,
\int dx\, x_1^2  (\omega\ast\omega)(t,x) (A^0 B^0)(t,x).  \label{eq:leading-nu-eff}
\EEQ

\bigskip\noindent
The simplest contributions to $v^{(0)}$ are obtained by taking $n=-1,0$
in (\ref{eq:two-point1}) and replacing
\BEQ R_{\eta}^{(0)}(\tau^0=0)(\vec{s}(\vec{w}))=\frac{1}{1-\int dt\, dx\, V^{(0)}(\tau^0=0)(\vec{s}(\vec{w}))(t,x)}  \label{eq:Reta0}
\EEQ by its lowest-order term $1$. Demanding that the $"n=-1"$-term compensates
exactly the sum for $n\ge 0$, we get an implicit equation for $v^{(0)}$,
\BEA
g^{(0)} v^{(0)}&=&  (g^{(0)})^2 \int dt' \, dx'\, 
G^0((t,x),(t',x'))\langle (\eta(t,x)-v^{(0)})(\eta(t',x')-v^{(0)})\rangle \nonumber\\
&&\qquad\qquad + O((g^{(0)}+g^{(0)}v^{(0)})^2 g^{(0)}v^{(0)})+ O((g^{(0)}+g^{(0)}v^{(0)})^4)   \label{eq:4.21}
\EEA
The implicit function theorem yields a unique solution 
\BEQ v^{(0)}=g^{(0)} \int dt'\, dx'\, G^0((t,x),(t',x')) \, \langle \eta(t,x)
\eta(t',x')\rangle + O ((g^{(0)})^3),  \label{eq:4.22}\EEQ
{\em provided} one can show that the series in $n$ converges, and that subleading
terms are indeed bounded as suggested in (\ref{eq:4.21}) and (\ref{eq:4.22}). This
is our next task.

\bigskip\noindent
{\bf D (bounds).}  {\em We now proceed to bound $v^{(0)}$ and 
$\nu_{eff}-\nu^{(0)}$.} 

\medskip\noindent  Let us first bound {\em scale $0$ resolvents}. They are of the
form (\ref{eq:Reta0}), where 
\BEQ V_{\eta}^{(0)}(\tau^0=0)(\vec{s}(\vec{w}))(t,x):=B^0(\vec{s}(\vec{w}))(\cdot,(t,x))\tilde{\eta}(t,x)
A^0(\vec{s}(\vec{w}))((t,x),\cdot)
\EEQ
where $\tilde{\eta}(t,x):=g^{(0)}(\eta(t,x)-v^{(0)})$.
Now, as explained in \S \ref{subsection:hor-cluster}, only the $\tilde{\eta}$'s belonging
to the image $|\T|$ of the connected component $\T$ (i.e. polymer) of $\F^0$ containing $(t,x)$ 
contribute. Denote then $\tilde{\eta}_{|\T|}(t,x):={\bf 1}_{(t,x)\in|\T|}\tilde{\eta}(t,x)$
the restriction of $\tilde{\eta}$ to $|\T|$. 
Expanding each $R^{(0)}_{\eta}(\tau^0=0)(\vec{s}(\vec{w}))((t_i,x_i),(t_{i+1},x_{i+1}))$ yields  $\del((t_i,x_i),(t_{i+1},x_{i+1}))+\Big(R^{(0)}_{\eta}(\tau^0=0)(\vec{s}(\vec{w}))((t_i,x_i),(t_{i+1},x_{i+1})) - \del((t_i,x_i),(t_{i+1},x_{i+1})) \Big)$,
with (expanding (\ref{eq:Reta0}))

\BEA && \Big| R^{(0)}_{\eta}(\tau^0=0)(\vec{s}(\vec{w}))((t_i,x_i),(t_{i+1},x_{i+1})) - \del((t_i,x_i),(t_{i+1},x_{i+1})) \Big| \nonumber\\
&&= \Big|
  B^0(\vec{s}(\vec{w}))((t_i,x_i),\cdot)  |0\rangle\,  \tilde{\eta}(\cdot) A^0(\vec{s}(\vec{w}))(\cdot,(t_{i+1},x_{i+1}))  \langle 0|\,   \nonumber\\
&& \qquad
+  B^0(\vec{s}(\vec{w}))((t_i,x_i),\cdot) |0\rangle\  \tilde{\eta}(\cdot)  A^0(\vec{s}(\vec{w}))(\cdot,\cdot) \langle 0|\,  B^0(\vec{s}(\vec{w}))(\cdot,\cdot) |0\rangle\ 
\nonumber\\
&&\qquad\qquad \tilde{\eta}(\cdot) A^0(\vec{s}(\vec{w}))(\cdot,(t_{i+1},x_{i+1})) \langle 0|\, +\ldots \Big| \nonumber\\
&&\le    B^0((t_i,x_i),\cdot) |0\rangle \ |\tilde{\eta}_{|\T|}(\cdot)|\,  A^0(\cdot,(t_{i+1},x_{i+1}))  \langle 0| \nonumber\\
&&\qquad +  B^0((t_i,x_i),(t'_i,x'_i)) |0\rangle\  |\tilde{\eta}_{|\T|}(t'_i,x'_i)| \ \cdot\ G_{|\eta_{|\T|}|}((t'_i,x'_i),(t'_{i+1},x'_{i+1})) \ \cdot\  \nonumber\\
&& \qquad\qquad\qquad \cdot\ |\tilde{\eta}_{|\T|}(t'_{i+1},x'_{i+1})| \,  A^0((t'_{i+1},x'_{i+1}),(t_{i+1},x_{i+1})) \langle 0|.   \nonumber\\
 \label{eq:pre-ARB-expansion}
\EEA

\medskip\noindent {\em Remark that (as follows from causality and from the fact that
boxes  of $\cal S$ are connected through $A^0$'s and $B^0$'s)  $t'_i-t'_{i+1}\le t_i-t_{i+1}\le 2$.} \\ Thus (letting $|\eta_{\Del}|:=\sup_{(t,x)\in\Del} |\eta(t,x)|$ for $\Del\in\D^0$)
\BEA && \langle G_{|\eta_{|\T|}|}((t'_i,x'_i),(t'_{i+1},x'_{i+1})\rangle \le G^{2\to}((t'_i,x'_i),(t'_{i+1},x'_{i+1})) \max
\langle \prod_{\Del\in\T} e^{\theta_{\Del}g^{(0)}|\eta_{\Del}|} \rangle \nonumber\\
&& \qquad  \le G^{2\to}((t'_i,x'_i),(t'_{i+1},x'_{i+1}))
\max e^{c(g^{(0)})^2 \sum_{\Del}\theta_{\Del}^2} \le G^{2\to}((t'_i,x'_i),(t'_{i+1},x'_{i+1})) e^{c'(g^{(0)})^2}  \label{Geta1} \nonumber\\
\EEA
if the maximum ranges over all possible choices of  occupation times 
$\theta_{\Del}:=|\{0\le s\le t'_i-t'_{i+1}\ |\ B_s\in\Del\}|$ for the
Brownian bridge from $(0,x'_i)$ to $(t'_i-t'_{i+1},x'_{i+1})$, since
$\sum_{\Del} \theta_{\Del}^2\lesssim \sum_{\Del} \theta_{\Del}=t'_i-t'_{i+1}\lesssim 1$. The bound for $\langle \prod_{\Del\in\T} e^{\theta_{\Del}g^{(0)}|\eta_{\Del}|} \rangle$ is obtained by rewriting the product $\prod_{\Del\in\T}
(\cdots)$ as a finite product, 
$\prod_{\vec{\eps}} \left(\prod_{\Del\in \T_{\vec{\eps}}}  (\cdots) \right)$, where
$\vec{\eps}\in\{0,1\}^{d+1}$ and $\T_{\vec{\eps}}:=\{\Del=[k_0,k_0+1)\times[k_1,k_1+1]\times\cdots
\times [k_d,k_d+1]\ |\ k_i-\eps_i\equiv 0 \mod 2, i=0,\ldots,d\}$, each of these
a product of independent variables, and uses H\"older's inequality, 
$\Big|\langle \prod_{\vec{\eps}} X_{\vec{\eps}} \rangle \Big| \le \prod_{\vec{\eps}} \left(\langle (X_{\vec{\eps}})^{2^{d+1}} \rangle \right)^{2^{-(d+1)}}$.

\medskip\noindent However, because the $G_{|\eta_{|\T|}|}((t'_i,x'_i),(t'_{i+1},x'_{i+1}))$, 
$i\ge 1$ are not independent in general, one should make the following easy
adaptation of the argument around (\ref{Geta1}). Split the total time interval
$[t',t]$ in (\ref{eq:two-point1}) into a union $I_1\cup I_2\cup\cdots$, 
$I_1:=[t_{i_1},t_{i_0}]$, $I_2:=[t_{i_2},t_{i_1}],\ldots$,  $i_0<i_1<i_2<\ldots$,
in such a way that $t_{i_{k-1}}-t_{i_k-1}<1<t_{i_{k-1}}-t_{i_k}$, and bound
as in (\ref{Geta1}) the products $\langle Y_k^2\rangle:=\Big\langle
\Big( \prod_{i=i_{k-1}}^{i_k-2} G_{|\eta_{|\T|}|}((t'_i,x'_i),(t'_{i+1},x'_{i+1}))
\Big)^2 \Big\rangle$. Since $t_{i_{k-1}}-t_{i_{k+1}}>2$, the random variables
$(Y_{2k+\eps})_k$, $\eps=0,1$ are independent, hence one concludes as above
using $|\langle\prod_k Y_k\rangle|\le (\prod_k \langle Y_{2k}^2\rangle)^{1/2}
(\prod_k \langle Y_{2k+1}^2\rangle)^{1/2}$.

%\medskip\noindent {\em Consider now a general random resolvent $G_{|\eta_{|\T|}|
%}(t'_i,x'_i),(t'_{i+1},x'_{i+1}))$.}  This case is treated in the same way by remarking
% that, if\\
%$T+1>t'_i>T>\ldots>T-n+1>t'_{i+1}>T-n$,
%\BEA && G_{|\eta_{|\T|}|}(t'_i,x'_i),(t'_{i+1},x'_{i+1}))\le \int_{|\T|} dy_0\cdots 
% \int_{|\T|} dy_{n-1}\,  \nonumber\\
%&& \qquad
 %G_{|\eta_{|\T|}|}((t'_i,x'_i),(T,y_0)) G_{|\eta_{|\T|}|}((T,y_0),(T-1,y_1))\cdots 
 %G_{|\eta_{|\T|}|}((T-n+1,y_{n-1}),(t'_{i+1},x'_{i+1})). \nonumber\\
 %\EEA

%The integration $\int dy_i G((T-i+1,y_{i-1}),(T-i,y_i))$ yields simply $1$, hence

%\BEQ \langle G_{|\eta_{\T}|}(t'_i,x'_i),(t'_{i+1},x'_{i+1}))\rangle
% \le C^{t'_i-t'_{i+1}}, \label{Geta2} \EEQ
%where $C=O(e^{c(g^{(0)})^2})=O(1)$.

\medskip\noindent {\em Let us now bound in average the product of  the $\eta$-dependent terms along ${\cal S}$}, namely,
the product of the dangling $\tilde{\eta}$'s with the $R^{(0)}_{\eta}$'s. Compared
to (\ref{Geta1}), one must now face the case when $O(n(\Del))$
dangling $\tilde{\eta}$'s are produced inside a box $\Del$, where $n(\Del)$ is the
(a priori arbitrary large) coordination number of $\Del$ in $\T$. Keeping aside for
further use the small factor $O(g^{(0)})$ per vertex, this leads to
replacing $\langle e^{\theta_{\Del}g^{(0)}|\eta_{\Del}|} \rangle$ in (\ref{Geta1}) by
$\langle (|\eta_{\Del}|+O(1))^{n(\Del)}  e^{\theta_{\Del}g^{(0)}|\eta_{\Del}|} \rangle
=e^{c'(g^{(0)})^2\theta_{\Del}^2}\ \cdot\ O(C^{n(\Del)}\Gamma(n(\Del)/2))$, with $C=O(1)$. These factors, traditionally called {\em local factorials}, are easily shown to
pose no real threat to the convergence of the sum over all polymers. Namely,
if $\Del\sim_{\T} \Del_1,\ldots,\Del_{n(\Del)-1}$,  $d(\Del,\Del_1)\le\ldots\le d(\Del,\Del_{n(\Del)-1})$, then (i) $d(\Del,\Del_n)\gtrsim n^{1/d}$; (ii) for each $n=1,\ldots,n(\Del)-1$, the string $\cal S$ contains a propagator, either $A^0((t,x),(t_n,x_n))$ or
$B^0((t,x),(t_n,x_n))$, with $(t,x)\in\Del,(t_n,x_n)\in\Del_n$. Rewriting $A^0((t,x),(t_n,x_n))$ as
$e^{-\frac{c}{2} |x-x_n|^2} \ \cdot\ \tilde{A}^0((t,x),(t_n,x_n))$, where $c$ is as
in Lemma \ref{lem:multi-scale-estimates-A-B}, one sees that $\tilde{A}^0(\cdot,\cdot)$
has the same scaling properties as $A^0(\cdot,\cdot)$, and has furthermore retained
the same Gaussian type space-decay, only with different constants. Putting (i) and (ii)
together, one sees easily that
\BEQ
C^{n(\Del)} \Gamma(n(\Del)/2) \ \cdot\ \prod_n e^{-\frac{c}{2}|x-x_n|^2}
\lesssim  C^{n(\Del)} \Gamma(n(\Del)/2) e^{-c'n(\Del)^{1+2/d}}  =O(1).
\label{eq:local-factorials}
\EEQ
 Thus
(at the price of replacing the  $A^0$'s and $B^0$'s along the string by 
propagators $\tilde{A}^0$'s, $\tilde{B}^0$'s with equivalent bounds), one has got rid
of local factorials.

\medskip\noindent {\em Finally,}  $|v^{(0)}|$, $|\nu_{eff}-\nu^{(0)}|$ and more
generally 
\BEQ I_{p,q}:=\int dt'\, \int dx'\, |x-x'|^p |t-t'|^q \Sigma((t,x),(t',x')), \EEQ 
$p+q\le 3$, see (\ref{eq:4.8}), {\em are simply bounded by
a sum over the number $n$ of vertices,}
\BEA && I_{p,q}\le   \sum_{n\ge 1} 
(Cg^{(0)})^{n+1} \int dt_1\, dx_1\,  A^0((t,x),(t_1,x_1)) \nonumber\\
&&\qquad\qquad\qquad \int dt'_1\, dx'_1\,  \underline{G}^{2\to}((t_1,x_1),(t'_1,x'_1)) 
 \int dt_2\, dx_2\,  B^0((t'_1,x'_1),(t_2,x_2))
\nonumber\\
&&\qquad
\cdots \int dt_{n}\, dx_{n}\,   A^0((t_{n-1},x_{n-1}),(t_{n},x_{n}))   \int dt'_{n}\, dx'_{n}\, \underline{G}^{2\to}((t_{n},x_{n}),(t'_{n},x'_{n})) \nonumber\\
&&\qquad\qquad\qquad 
 \int dt'\, dx'\, F_3((t,x),(t_1,x_1),(t'_1,x'_1),\ldots,(t'_n,x'_n),(t',x'))\,  B^0((t'_{n},x'_{n}),(t',x')). \nonumber\\  \label{eq:v0-bd}
\EEA
where $\underline{G}^{2\to}(\cdot,\cdot):=\del(\cdot,\cdot)+G^{2\to}(\cdot,\cdot)$, $C=O(1)$ and  (by H\"older's inequality)
\BEA  && F_3(\cdot) = O(n^2) \Big[ (1+t-t_1)^3+(1+t_1-t'_1)^3+\cdots+(1+t_n-t'_n)^3+(1+t'_n-t')^3  \nonumber\\
&&\qquad + (1+|x-x_1|)^3+(1+|x_1-x'_1|)^3+\cdots+(1+|x_n-x'_n|)^3+(1+|x'_n-x'|)^3
\Big].
\nonumber\\
\EEA  
Integrating space-time variables in chronological order, and using
\BEA && (1+t-t'+|x-x'|)^3 A^0((t,x),(t',x')),(1+t-t'+|x-x'|)^3 G^{2\to}((t,x),(t',x'))
\nonumber\\
&&\qquad\lesssim (t-t')^{-d/2}
e^{-c|x-x'|^2/(t-t')}
\, \cdot\, {\bf 1}_{t-t'=O(1)},
\EEA
 one gets
\BEQ g^{(0)}|v^{(0)}|,|\nu_{eff}-\nu^{(0)}|\le \sum_{n\ge 1} n^3 (C'g^{(0)})^{n+1}=O((g^{(0)})^2)
\EEQ
for another constant $C'=O(1)$.

\medskip\noindent Finally, it is clear from (\ref{eq:v0-bd}) that 
$v^{(0)}-v^{(0)}(T)$ involves only terms in the sum for which $n\gtrsim T$; thus it
is of order $O((C'g^{(0)})^{cT})$ for some constant $c>0$.

 %%%%%%%%%%%%%%%%%%%%%%%%%%%

\subsection{Four-point  function}  \label{subsection:four-point}

%%%%%%%%%%%%%%%%%%%%%%%%%%%%%%ש

Next, we discuss briefly connected {\em four-point functions}, which contribute corrections
to the noise strength $D$.  The correct way
to get an understanding \`a la Wilson of the induced flow for the parameter $D$ is
a priori to sum inductively for each fixed scale $j$ over all diagrams of lowest scale $\le j$
with four external $A$- or $B$-propagators of scale $>j$. In practice this would 
lead to introduce further  scale counterterms of the form $V^{(j)}(t,x)=B^{\to j}(\cdot,(t,x))
g^{(j)}\eta_j(t,x) A^{\to j}((t,x),\cdot)$, where $(\eta_j)_{j\ge 0}$ are independent
copies of $\eta$, with scale $\tau$-prefactors, yielding the whole machinery
of multi-scale cluster expansions. Fortunately, since the insertion of such diagrams
inside the expansion yields power-like vanishing contributions in the large-scale limit, 
such counterterms need not be introduced by hand  to make the expansion
convergent. We shall actually compute directly in \S \ref{subsection:hh} an {\em effective} value
$D_{eff}\equiv D^{(\infty)}$ for 
$D$ by considering the large-scale limit of the connected two-point function
$\langle h(\cdot)h(\cdot)\rangle$.  We shall be content here with a few indications about {\em how}
four-point functions are produced by the expansion. This subsection may be skipped
since it is not used in the proof of our Main Theorem.

\bigskip\noindent In order to obtain a four-point function, one needs two strings. Let us denote
by the index $\alpha$  vertices produced on the first string, and by the index
$\beta$ those produced on the second string. A  component  connected by the cluster is made up of a piece of string ${\cal S}_{\alpha}$ and a piece of
string ${\cal S}_{\beta}$, both of the type (\ref{eq:two-point1}).  One thus
obtains a diagram with  4 external vertices, $B^{\to(j+1)}(\cdot,(t_{\alpha},x_{\alpha})) \, \cdot\, \Big[ \, \cdot\, \Big]\, A^{\to(j+1)}((t'_{\alpha},x'_{\alpha},\cdot) \ \cdot\ B^{\to(j+1)}(\cdot,(t_{\beta},x_{\beta})) \, \cdot\, \Big[ \, \cdot\, \Big]\, A^{\to(j+1)}((t'_{\beta},x'_{\beta},\cdot)$.
To get a connected contribution, we assume that  $\eta(t_{\alpha},x_{\alpha})$ contracts with  $\eta(t_{\beta},x_{\beta})$, and similarly, $\eta(t'_{\alpha},x'_{\alpha})$ contracts with  $\eta(t'_{\beta},x'_{\beta})$. Then this means that we obtain a very simple "ladder diagram",  whose
leading term is 

\BEA && B^{\to(j+1)}(\cdot,(t_{\alpha},x_{\alpha})) \, \cdot\, \nonumber\\
&&\qquad\cdot\, 
\Big[ \eta(t_{\alpha},x_{\alpha}) A^{j\to}((t_{\alpha},x_{\alpha}),\cdot) B^{j\to} (\cdot,(t'_{\alpha},x'_{\alpha}))  \eta(t'_{\alpha},x'_{\alpha}) \Big]\,   \cdot\,  A^{\to(j+1)}((t'_{\alpha},x'_{\alpha}),\cdot)  \, \cdot  \nonumber\\
 && \cdot\   B^{\to(j+1)}(\cdot,(t_{\beta},x_{\beta})) \ \cdot\nonumber\\
&&\qquad \cdot\,   
\big[ \eta(t_{\beta},x_{\beta}) A^{j\to}((t_{\beta},x_{\beta}),\cdot)  B^{j\to}(\cdot,(t'_{\beta},x'_{\beta})) \eta(t'_{\beta},x'_{\beta}) \Big] \, \cdot A^{\to(j+1)}((t'_{\beta},x'_{\beta}),\cdot)
\label{eq:ladder} 
\EEA
with $d((t_{\alpha},x_{\alpha}),(t_{\beta},x_{\beta})),d((t'_{\alpha},x'_{\alpha}),(t'_{\beta},x'_{\beta}))=O(1)$. Renormalization corrections 
are due precisely to these (and more complicated) ladder diagrams, with $(t_{\alpha},x_{\alpha}),(t_{\beta},x_{\beta})$, resp.  $(t'_{\alpha},x'_{\alpha}),(t'_{\beta},x'_{\beta})$
belonging to $\Del$, resp. $\Del'$, where $\Del$, $\Del'$ are two distinct scale $0$
boxes where the $\eta$'s contract two-by-two. If $j=0$ then short-distance "crossed" $\eta$-contractions
are also possible.

%%%%%%%%%%%%%%%%%%%%%%%
%%%%%%%%%%%%%%%%%%%%%%%

\section{Final bounds} \label{sec:bounds}

%%%%%%%%%%%%%%%%%%%שש
%%%%%%%%%%%%%%%%%%%%%ש

We are now, at long last, ready to prove our Main Theorem. Roughly speaking, 
$N$-point functions $\langle h(t_1,x_1)\cdots h(t_N,x_N)\rangle$ have been rewritten
in terms of a series, that is, an (infinite) sum over polymers. Obviously, the first task is to ensure that this series is
convergent. This turns out to be the main point in the section; once this is understood,
the scaling behavior of $N$-point functions will be essentially obtained by looking at
the terms of lower order in $g^{(0)}$ in the series.

%%%%%%%%%%%%%%%%%%%%%%%%%%

\subsection{Small noise/large noise boxes}

%%%%%%%%%%%%%%%%%%%%%%%%%ששש

\begin{Definition}
Let $\Del\in\D^0$. Then $\Del$ is said to be a size $k$ {\em large field} box $(k\ge 0)$ if
$2^k \lambda^{-1/2}<\sup_{\Del} |\eta|\le 2^{k+1}\lambda^{-1/2}$. 

Denote by $\D^0_{LF,k}$ the set of size $k$ large field boxes,  by $\D^0_{LF}:=\uplus_{k\ge 0} \D^0_{LF,k}$ the set of all large field boxes,
and by $\D^0_{SF}:=\D^0 
\setminus \D^0_{LF}$ its complementary. The region $\D^0_{LF}$ is called
the {\em large field region}, and the region $\D^0_{SF}$ the {\em small field
region}.
\end{Definition}

By standard Gaussian deviations, if $\Del\in\D^0$, then 
\BEQ \proba[\Del\in\D^0_{\LF,k}]\le e^{-c2^{2k}/\lambda}, \qquad k\ge 0.
\label{eq:standard-G-dev}
\EEQ
The bound (\ref{eq:standard-G-dev}) also holds trivially if $\Del\in\D^0_{SF}$ by
letting formally $k=-\infty$. This trick allows to handle small noise and large noise
boxes on equal footing.

% On the other hand, each large field box comes with a small
%statistical weight $\proba[\sup_{\Del} |\eta|>\lambda^{-1/2}]=O(e^{-c/\lambda})$,
%so that no expansion is needed in the large field region.

%%%%%%%%%%%%%%%%%%%שש

\subsection{Vertex insertions and contour integrals}  \label{subsection:vertex-insertions}

%%%%%%%%%%%%%%%%%%%%

Let us recapitulate the previous steps. We start from an $N$-point function,
\BEQ \langle h(t_1,x_1)\cdots h(t_N,x_N)\rangle=\left(\frac{\nu^{(0)}}{\lambda}\right)^N F_N, \EEQ
where 
\BEA F_N &:=&  \langle \log(w(t_1,x_1))\cdots \log(w(t_N,x_N))\rangle  \nonumber\\
&=& \Big{\langle} \log\left( \int  dy_1 \,  A((t_1,x_1),\cdot) (1-V_{\eta})^{-1} (\cdot,\cdot)
 B(\cdot,(0,y_1))\,  e^{\frac{\lambda}{\nu^{(0)}} h_0(y_1)} \right) \cdots  \nonumber\\
 &&\qquad \log\left( \int  dy_N\,   A((t_N,x_N),\cdot)
  (1-V_{\eta})^{-1} (\cdot,\cdot) B(\cdot,(0,y_1))\,  e^{\frac{\lambda}{\nu^{(0)}} h_0(y_N)} \right) \Big{\rangle}
 \EEA
 and $V_{\eta}:=\int dt\, dx\, V_{\eta}(\tau=1)(t,x)=\int dt\, dx\, B(\cdot,(t,x))\left(g^{(0)}(\eta(t,x)-v^{(0)})\right)A((t,x),\cdot)$.
 Then we:

\begin{enumerate}
\item  apply to $F_N$ the horizontal and vertical cluster expansions; this results
in a sum over forests $\F^0\in{\cal F}^0$ of a rational function (see (\ref{eq:DDD}))
in {\em strings};
\item displace external legs; 
\item  contract the dangling $\eta$'s;
\item apply Mayer's expansion to scale $0$ two-point diagrams;
\item {\em factorize} the scale $0$ two-point diagram contributions. By construction
these are exactly compensated by the counterterms.
\end{enumerate} 

\medskip\noindent Consider now the various vertex insertions (in form of a kernel),
\BEQ \tilde{V}_{\alpha}(\vec{s}(\vec{w});z_{\alpha}):=c_{\alpha}(\vec{s}(\vec{w}))
\begin{cases} V_{\alpha}(z_{\alpha}) 
\qquad\ \  {\mathrm{if}}\  z_{\alpha}=z_{\ell}, \ell\in L_{\eta} \\
\tilde{\eta}(z_{\alpha})V_{\alpha}(z_{\alpha}) \qquad {\mathrm{otherwise}} 
\end{cases}  \label{eq:ctildealpha}
\EEQ
where 
\BEQ c_{\alpha}(\vec{s}(\vec{w})):=\tilde{s}_{\Del'_{\alpha},\Del_{\alpha}} \tilde{s}_{\Del_{\alpha},\Del''_{\alpha}}, \qquad 
\tilde{s}_{\Del,\Del'}=\begin{cases} 1 \qquad
{\mathrm{if}}\  \{\Del,\Del'\}\in L_G \\ s_{\Del,\Del'} \qquad {\mathrm{otherwise}}
\end{cases}
\EEQ
and
\BEQ V_{\alpha}(z_{\alpha})(z'_{\alpha},z''_{\alpha}):=g^{(0)} \partial^{\kappa'_{\alpha}} B^{j'_{\alpha}}(z'_{\alpha},z_{\alpha})
\partial^{\kappa''_{\alpha}} A^{j''_{\alpha}}(z_{\alpha},z''_{\alpha})
\label{eq:Valpha}
\EEQ 
on the strings, with $z_{\alpha}=(t_{\alpha},x_{\alpha}),z'_{\alpha}=(t'_{\alpha},x'_{\alpha})\in\Del'_{\alpha},z''_{\alpha}=(t''_{\alpha},x''_{\alpha})\in\Del''_{\alpha}$,  $\partial^{\kappa'_{\alpha}}:=\partial_{t'}^{\kappa'_{\alpha,0}} \partial_{x'}^{\vec{\kappa}'_{\alpha}}$, and similarly for
$\partial^{\kappa''_{\alpha}}$; $\alpha$ being some dummy index. Let 
$|\kappa'_{\alpha}|:=2\kappa'_{\alpha,0}+|\vec{\kappa}'_{\alpha}|$ be the
parabolic order of derivation; in particular, $|\kappa'_{\alpha}|=3$ if and only
if $\partial^{\kappa'_{\alpha}}=\partial_{t'} \nabla_{x'}$ or $\nabla_{x'}^{\vec{\kappa}'_{\alpha}}$, 
$|\vec{\kappa}'_{\alpha}|=3$; and similarly for $\kappa''_{\alpha}$.  By assumption $z_{\alpha}$ ranges over some box $\Del^0_{\alpha}$ of scale $0$, and $\Del'_{\alpha}
\in \D^{j'_{\alpha}},\Del''_{\alpha}\in\D^{j''_{\alpha}}$. 
There are four cases:
\begin{itemize}
\item[(i)] (no $\tau$-derivative, $0$-th scale vertices) $j'_{\alpha},j''_{\alpha}=0$, and  $\kappa'_{\alpha}=\kappa''_{\alpha}=0$;
\item[(ii)] (one $\tau$-derivative, beginning of $0$-th scale cluster) $j'_{\alpha}>0$, $j''_{\alpha}=0$, $\kappa''_{\alpha}=0$, $|\kappa'_{\alpha}|=0$ or $\ge 3$; 
\item[(iii)] (one $\tau$-derivative, end of $0$-th scale cluster) $j'_{\alpha}=0$, $\kappa'_{\alpha}=0$, and $j''_{\alpha}>0$, 
$|{\kappa}''_{\alpha}|=0$ or $\ge 3$.
\item[(iv)] (second $\tau$-derivative) $j'_{\alpha}>0$,  $j''_{\alpha}>0$,  and
$\kappa'_{\alpha}=\kappa''_{\alpha}=0$.
\end{itemize} 

To these, one must add insertions of a particular type, proportional to
$\del v^{(0)}$ (see (\ref{eq:v0-v0m})),
\begin{itemize}
\item[(v)] (boundary terms) $V_{\alpha}(z_{\alpha})(z'_{\alpha},z''_{\alpha}):=
\del v^{(0)}(t_{\alpha}) B^{j'_{\alpha}}(z'_{\alpha},z_{\alpha})
 A^{j''_{\alpha}}(z_{\alpha},z''_{\alpha})$ $(j''_{\alpha}\ge j'_{\alpha})$,
 resp. $\del v^{(0)}(t_{init}-t_{\alpha}) B^{j'_{\alpha}}(z'_{\alpha},z_{\alpha})
 A^{j''_{\alpha}}(z_{\alpha},z''_{\alpha})$ $(j''_{\alpha}\ge j'_{\alpha})$ 
 $(j''_{\alpha}<j'_{\alpha})$, where $t_{init}:=t_i$ if $V_{\alpha}(z_{\alpha})$
 is inserted on the $i$-th string, $i=1,\ldots,N$.
 \end{itemize}

%($\del\nu$-counterterm) $j'_{\alpha}>0$, $j''_{\alpha}>0$, and
%$\kappa'_{\alpha}=0$, "$\partial^{\kappa''_{\alpha}}$"$=\Del^{\to 0}$ ($
%|\kappa''_{\alpha}|=2$). These $\del V_{\eta}$-vertex insertions (see second line
%of (\ref{eq:dressed-vertex})) are accompagnied by a small prefactor $\frac{\nu_{eff}
%-\nu^{(0)}}{g^{(0)}}=O(g^{(0)})$. The operator $\Del^{\to 0}=\bar{\chi}^{(0)}\ast\Del$
% is
%not strictly local, so $\del V_{\eta}$ is not stricto sensu of the form
% (\ref{eq:Valpha}). Instead of delocalizing the vertex, it is simpler to replace
%  implicitly the 
%propagator $B^{j'_{\alpha}}$ by another,$\tilde{B}^{j'_{\alpha}}=B^{j'_{\alpha}}\ast
% \bar{\chi}^{(0)}$, with equivalent properties.

\medskip\noindent Vertices of type (iv) are responsible for the production
of the $v^{(0)}$ (see "$n=-1$" term in \S \ref{subsection:two-point}) and $\del\nu
$  (see second line of (\ref{eq:dressed-vertex})) counterterms; the $v^{(0)}$-counterterm is chosen in such a way as to cancel the two-point function (see
\S \ref{subsection:two-point}), while  $\del\nu$-counterterms are resummed into
the effective propagator $\tilde{G}_{eff}$ (see \S \ref{subsection:h} {\bf A.}). 
The contribution of scale $0$ vertices (i) is bounded in \S \ref{subsection:h}
{\bf A.}

\medskip\noindent 
Vertex insertions of type  (ii), (iii)
have been differentiated by the scale $0$ renormalization. More precisely, letting
$\kappa'_{\alpha}$ be the order of differentiation of a low-momentum $B^{j'_{\alpha}}$-propagator entering a
given $0$-th scale cluster (ii), and $\kappa''_{\alpha'}$ that of a low-momentum
$A^{j''_{\alpha'}}$-propagator exiting the same $0$-th scale cluster, one  
sets as in \S \ref{subsection:two-point}: $(|\kappa'_{\alpha}|\ge 3,\kappa''_{\alpha'}=0)$ if $j'_{\alpha}\ge j''_{\alpha'}$, $(\kappa'_{\alpha}=0,|\kappa''_{\alpha'}|\ge 3)$ if $j'_{\alpha}< j''_{\alpha'}$. From the point of view of power-counting (see below), we have
thus produced an essential {\em small factor} 
\BEQ O(2^{-\frac{3}{2}\max(j'_{\alpha},j''_{\alpha'})}),
\label{eq:small-factor} \EEQ
that is, $O(2^{-\frac{3}{2}j})$ per half of the low-momentum fields $A^j$ or $B^j$, having the
same effect as $\nabla^{3}$, or (considering a chronological sequence
 $A^j(\cdot,\cdot) \langle j| \ B^j(\cdot,\cdot) |j\rangle =G^j(\cdot,\cdot)$),  $O(\nabla^3)$ in average per low-momentum $G$-field. 
 
\medskip\noindent Finally, boundary vertices of type (v) enjoy an {\em exponentially small factor}. Namely, assuming e.g. that $j''_{\alpha}\ge j'_{\alpha}$, the boundary correction
$\del v^{(0)}(t_{\alpha})$ to $v^{(0)}$ is $O((Cg^{(0)})^{ct_{\alpha}}$, which
is $\lesssim t_{\alpha}^{-3/2}\lesssim 2^{-\frac{3}{2}j''_{\alpha}}$. Hence such
vertices may and will be considered -- from the power-counting point of view -- as
$O(1)$ times a vertex of type (ii) or (iii).

\medskip\noindent Recalling that  these vertices are produced
anywhere along the strings by the cluster expansion,  their contributions may be resummed 
as follows. We first need some notations. Let:

\begin{itemize}

\item  $I(\F^{0})=\{I_{\alpha}\}_{\alpha}$ be the set of vertex insertions;

\item $\Del_{\alpha}\in\D^0$ be the scale-0 box where $z_{\alpha}$
(see (\ref{eq:Valpha})) is located; 
 
\item  $L(\F^{0})$ be the set of horizontal cluster links, and  $L_{\eta}\subset L(\F^0)$ those coming specifically
from the cluster on  $\eta $ (compare with \S \ref{subsection:two-point});
\item    $L_{vert}$ be the set of   links coming from the vertical  cluster expansion;
\item $L_{Mayer}$ be the set of  Mayer links.
\end{itemize}

\medskip\noindent Then
  \BEA
&& F_N= \sum _{\F^{0}\in{\cal F}^{0}}  \sum_{L_G,L_{\eta},L_{vert},L_{Mayer},\vec{\mu}}\Big(\prod _{\ell \in L(\F^{0})} \int _{0}^{1} dw_{\ell }\Big)
 \Big(\prod _{\ell \in L_{Mayer}(\F^{0})} \int _{0}^{1} dS_{\ell }\Big)\  \  {\mathrm{Mayer}}(\vec{S}) \nonumber\\
 &&\qquad
\Big(\prod _{\ell\in L_{\eta }(\F^{0})} \langle\eta (z_{\ell }) \eta (z'_{\ell })\rangle\Big) \qquad \cdot\qquad
\Big\langle \Big(\prod _{\alpha \in I(\F^{0}) }\frac{d}{d\gamma_{\alpha }}|_{\gamma_{\alpha }=0}\Big)
 \nonumber\\
&&\  \prod _{j=1}^{N}
\log\left(\int dy_{j} \,  \tilde{A}(\vec{s}(\vec{w}))((t_j,x_j),\cdot) \frac{1}{1-V_{\eta}(\tau)-\sum_{\alpha\in I(\F^{0}) }\gamma_{\alpha}\tilde{V}_{\alpha}} (\cdot,\cdot)
 \tilde{B}(\vec{s}(\vec{w}))(\cdot,(0,y_j)) \, e^{\frac{\lambda}{\nu^{(0)}} h_0(y_j)} \right) 
  \Big\rangle_{s(w)}\nonumber\\  \label{eq:FN}
\EEA

where $\tilde{C}(\vec{s}(\vec{w})):=\tilde{C}^0(\vec{s}(\vec{w}))+C^{\to 1}$, $\tilde{C}^0(\vec{s}(\vec{w}))(z,z'):=
\begin{cases} C^0(z,z') \qquad {\mathrm{if}}\ (\Del^0_z,\Del^0_{z'})\in L_G \\ s_{\Del^0_z,
\Del^0_{z'}} C^0(z,z') \qquad {\mathrm{otherwise}} \end{cases}$
$(C=A,B)$;  $\tilde{V}_{\alpha}:=\int_{\Del_{\alpha}}  dz_{\alpha}\, \tilde{V}_{\alpha}(\vec{s}(\vec{w});z_{\alpha})$, and $V_{\eta}(\tau):=\int dt\, dx\, V_{\eta}(\tau)(\vec{s}(\vec{w}))(t,x)$ is the space-time
integration of the dressed vertex (\ref{eq:dressed-vertex}).
Note that the $s$-dependence in this expression is trivial when
it comes to {\em bounds}, since $|\tilde{C}(\vec{s}(\vec{w})(\cdot,\cdot)|\le |C(\cdot,\cdot)|$, $C=A,B$, and similarly
$|c_{\alpha}(\vec{s}(\vec{w})|\le 1$ (see (\ref{eq:ctildealpha})).

\medskip\noindent By causality, the vertex insertions may be re-expanded along
the string number $i=1,\ldots,N$ into a finite sum ${\cal S}_i$ as follows: letting
$\vec{\gamma}:=(\gamma_{\alpha})_{\alpha}$, 
\BEA && {\cal S}_i(\vec{\gamma}):= \tilde{A}(\vec{s}(\vec{w}))(t_i,x_i),\cdot) (1-V_{\eta}(\tau)-\sum_{\alpha}\gamma_{\alpha}V_{\alpha})^{-1}(\cdot,\cdot) \int dy_i\,  \tilde{B}(\vec{s}(\vec{w}))(\cdot,(0,y_i))
\, e^{\frac{\lambda}{\nu^{(0)}} h_0(y_i)} \nonumber\\
&&\ =\sum_{\alpha_1,\alpha_2,\ldots} 
\tilde{A}(\vec{s}(\vec{w}))((t_i,x_i),\cdot) \ \cdot\ \nonumber\\
&&\ \cdot\    \left(\int_{\Del'_{\alpha_1}} dz'_1\, \int_{\Del_{\alpha_1}}
dz_{\alpha_1} \int_{\Del''_{\alpha_1}}dz''_1 \right) R_{\eta}(\tau)(\cdot,z'_1)  \Big\{\gamma_{\alpha_1} \tilde{V}_{\alpha_1}(\vec{s}(\vec{w}));z_{\alpha_1})(z'_1,z''_1) \Big\} \ \cdot
\nonumber\\ 
&& \cdot \  \left(\int_{\Del'_{\alpha_2}} dz'_2\, \int_{\Del_{\alpha_2}} dz_{\alpha_2}\,  \int_{\Del''_{\alpha_2}}dz''_2 \right) R_{\eta}(\tau)(z''_1,z'_2)  \Big\{ \gamma_{\alpha_2} \tilde{V}_{\alpha_2}(\vec{s}(\vec{w}))(z_{\alpha_2})(z'_2,z''_2) \Big\} \ \cdots,
\nonumber\\  \label{eq:vertex-insertions}
\EEA
with {\em main term} (disregarding propagator renormalization, see \S \ref{subsection:h} {\bf A.})
\BEA &&  A((t_i,x_i),\cdot) \int dy_i\, B(\cdot,(0,y_i))e^{\frac{\lambda}{\nu^{(0)}} h_0(y_i)} = \int dy_i\,  G((t_i,x_i),(0,y_i)) e^{\frac{\lambda}{\nu^{(0)}} h_0(y_i)}
\nonumber\\
&&\qquad= 1+ e^{\nu^{(0)} t_i\Del} (e^{\frac{\lambda}{\nu^{(0)}} h_0}-1)(x_i) 
\le 1+\frac{\lambda}{\nu^{(0)}}\  e^{\frac{\lambda}{\nu^{(0)}} ||h_0||_{\infty}} \, \ \cdot\ (e^{\nu^{(0)}t_i\Del} |h_0|)(x_i) \nonumber\\
&&\qquad= 1+O(\lambda e^{\frac{\lambda}{\nu^{(0)}} ||h_0||_{\infty}}) \, \min(||h_0||_{L^{\infty}}, t_i^{-d/2} ||h_0||_{L^1}). 
\label{eq:main-term} 
\EEA

\medskip\noindent  {\em Since the result is analytic} in the parameters $\vec{\gamma}$ in a
neighborhood of $0$,
{\em  we may replace $\frac{d}{d\gamma_{\alpha}}\big|_{\gamma_{\alpha}=0} F(\gamma_{\alpha})$ by the Cauchy contour integral}
$$\frac{1}{2\II\pi} \oint_{\partial B(0,r_{\alpha})} \frac{d\gamma_{\alpha}}{\gamma_{\alpha}^2} F(\gamma_{\alpha}),$$ with {\em (defining $k_{\alpha}$ to be the
size of the large-field zone of $\Del_{\alpha}$ if $\Del_{\alpha}$ is large-field, i.e.  $\Del_{\alpha}\in 
\D^0_{LF,k_{\alpha}}$, $k_{\alpha}\ge 0$, and $k_{\alpha}=-\infty$ if $\Del_{\alpha}\in \D^0_{SF}$)} 
\BEQ r_{\alpha}\equiv r_{\alpha}(k_{\alpha}):=r'_{\alpha}(k_{\alpha}) r''_{\alpha}, \EEQ
 where

\BEQ (r'_{\alpha}(k_{\alpha}))^{-1}:=C (2^{k_{\alpha}+1})^{n(\Del_{\alpha})} e^{\lambda^{1/2} 2^{k_{\alpha}+1}}
\EEQ

\BEQ (r''_{\alpha})^{-1} := C g^{(0)} \int_{\Del'_{\alpha}} dz'_{\alpha}
\int_{\Del_{\alpha}} dz_{\alpha} \int_{\Del''_{\alpha}} dz''_{\alpha} \,  
|V_{\alpha}(z_{\alpha})(z'_{\alpha},z''_{\alpha})| \nonumber\\  \label{eq:r''alpha}
\EEQ
%& \approx & g^{(0)} (2^{k_{\alpha}+1})^{n(\Del_{\alpha})}  e^{\lambda^{1/2}
% 2^{k_{\alpha}+1}} 2^{-(j'_{\alpha}+j''_{\alpha})d/4}\,  2^{-j'|\kappa'_{\alpha}|/2}
%\, 2^{-j''|\kappa''_{\alpha}|/2}\, 
%e^{-c\del x^2/\del t} \EEA
for some large enough uniform constant $C$.
%, where $\del x:=\sup_{(t,x)\in\Del'_{\alpha},(t',x')\in\Del''_{\alpha}}
%|x'-x''|$, $\del t:=\sup_{(t,x)\in\Del'_{\alpha},(t',x')\in\Del''_{\alpha}}
%|t'-t''|$.
 Then 
 \BEQ |\gamma_{\alpha}|=r_{\alpha},\qquad \frac{1}{2\pi} \oint_{\partial B(0,r_{\alpha})} \frac{d|\gamma_{\alpha}|}{|\gamma_{\alpha}|^2}= 
r_{\alpha}^{-1}. \EEQ 
As we shall see, the $r_{\alpha}\equiv |\gamma_{\alpha}|$ have been chosen small enough
(depending on the order of magnitude of  the $(|\eta_{\Delta_{\alpha}}|)_{\alpha}$)
so that each  ${\cal S}_i(\vec{\gamma})$  is equal to $1+o(1)$, yielding
\BEA   F_N(\F^0,\vec{k}|\eta) &:=& \prod_{\alpha} {\bf 1}_{\Del_{\alpha}\in\D^0_{LF,k_{\alpha}}} \ \cdot\  \left[ \prod_{\alpha} \left(\frac{1}{2\II\pi} \oint_{\partial B(0,r_{\alpha})} \frac{d\gamma_{\alpha}}{\gamma_{\alpha}^2} \right)  \right]\ \left\{\prod_{i=1}^N  \log( {\cal S}_i(\vec{\gamma})) \right\}
\nonumber\\
 &=& O(1) \prod_{\alpha} r^{-1}_{\alpha}(k_{\alpha}),  \label{eq:FN-k} \EEA
 a {\em deterministic estimate} (but depending on $\vec{k}:=(k_{\alpha})_{\alpha}$).
This is step {\bf B.} in \S \ref{subsection:h}. 

\medskip\noindent  The next step  (see \S \ref{subsection:h},
step {\bf C.}) is to show that 
the averaged  infinite sum $\langle \sum_{\vec{k}} F_N(\F^0,\vec{k}|\eta)
\rangle $
is $\lesssim \prod_{\alpha} (r''_{\alpha})^{-1}$; or rather, to be precise,  $\lesssim \prod_{\alpha} (\tilde{r}''_{\alpha})^{-1}$, where (as in \S \ref{subsection:two-point}) $\tilde{r}''_{\alpha}$ is $r''_{\alpha}$ up to the replacement of $A^j,B^j$ with
equivalent kernels $\tilde{A}^j,\tilde{B}^j$.

\medskip\noindent
The final step is to show that the infinite sum
 $\sum_{\F^0}\prod_{\alpha} (\tilde{r}''_{\alpha})^{-1}$ converges; see step {\bf D.} in \S \ref{subsection:h}.

\medskip\noindent Obviously, in the course of the proof one must extract the lowest
order terms, which will give the leading behavior of the KPZ truncated functions.

%%%%%%%%%%%%%שש

\subsection{KPZ 1-point function}  \label{subsection:h}

%%%%%%%%%%%%%%%%%ש

Let us first consider the case of the $1$-point function $\langle h(t,x)\rangle=
\frac{\nu^{(0)}}{\lambda} \langle \log w(t,x)\rangle$, where there is
only one string. One must prove that $\langle h(t,x)\rangle\overset{t\to\infty}{\to}
0$. We decompose the proof into four points (see discussion at the end of \S \ref{subsection:vertex-insertions}); the
first point {\bf A.} is a preparatory step. Except that {\bf A.} must be 
supplemented with a new power-counting argument (see {\bf A'.}), the same
scheme of proof of convergence is  used  for KPZ truncated functions of higher order, see \S 5.4, 5.5, where 
details are skipped, so that one can concentrate on the asymptotic
large-scale scaling functions.

\begin{itemize}
\item[{\bf A.}] {\bf (contribution of the random resolvents)} On each string, one finds
a number of  random
resolvents $R^{(0)}_{\eta}(\tau^0=0)(\vec{s}(\vec{w}))((t_i,x_i),(t_{i+1},x_{i+1}))$. As in \S 4.1 {\bf C.}, 
such resolvents may be expanded to order two, see (\ref{eq:pre-ARB-expansion}), 
\BEA && R^{(0)}_{\eta}(\tau^0=0)((t_i,x_i),(t_{i+1},x_{i+1}))=\del((t_i,x_i),(t_{i+1},x_{i+1}))  \nonumber\\
&&\ + B^0(\vec{s}(\vec{w}))((t_i,x_i),\cdot) \tilde{\eta}(\cdot) A^0(\vec{s}(\vec{w}))(\cdot,(t_{i+1},x_{i+1})) \nonumber\\
&&\ +B^0(\vec{s}(\vec{w}))((t_i,x_i),(t'_i,x'_i)) \tilde{\eta}(t'_i,x'_i) \ \cdot\ G_{\eta}(\vec{s}(\vec{w}))((t'_i,x'_i),(t'_{i+1},x'_{i+1})) \ \cdot\ \nonumber\\
&&\qquad \cdot\ \tilde{\eta}(t'_{i+1},x'_{i+1}) ((A^0(\vec{s}(\vec{w}))((t'_{i+1},x'_{i+1}),(t_{i+1},x_{i+1})) \nonumber\\  \label{eq:pre-ARB-2}
\EEA

with $t'_i-t'_{i+1}\le 2$.
Then $A^0(\vec{s}(\vec{w})(\cdot,\cdot)\le A^0(\cdot,\cdot)$,  $B^0(\vec{s}(\vec{w})(\cdot,\cdot)\le B^0(\cdot,\cdot)$  and
\BEA && G_{\eta}(\vec{s}(\vec{w}))((t'_i,x'_i),(t'_{i+1},x'_{i+1}))\le G_{|\eta_{|\T|}|}
((t'_i,x'_i),(t'_{i+1},x'_{i+1})) \nonumber\\
&&\qquad \le G^{2\to} ((t'_i,x'_i),(t'_{i+1},x'_{i+1}))
\max \prod_{\Del\in\T} e^{\theta_{\Del}g^{(0)}|\eta_{\Del}|}. 
\EEA
Furthermore,  the expansion (\ref{eq:pre-ARB-2}) has produced new $\eta$ fields,
$O(n(\Del))$ per box $\Del\in\T$.
Thus, to each large-field box $\Del_{\alpha}\in\T$ corresponds a  factor 
$r'_{\alpha}(k_{\alpha})  |\eta_{\Del_{\alpha}}|^{O(n(\Del_{\alpha}))} e^{\theta_{\Del_{\alpha}}g^{(0)}|\eta_{\Del_{\alpha}}|}=o(1)$. {\em Concluding:
a scale $0$ resolvent $R^{(0)}_{\eta}(\tau^0=0)(\cdot,\cdot)$ may be replaced by $\del(\cdot,\cdot)+ O(g^{(0)})G^{2\to}(\cdot,\cdot)$.}

\medskip\noindent On the other hand, one also finds low-momentum resolvents\\
\BEQ \del R_{\eta}:=(1-\del\nu B^{\to 1}\Del^{\to 0}A^{\to 1})^{-1} 
 \label{eq:A-low}
\EEQ
 (see
second line of (\ref{eq:dressed-vertex})). Sandwiched between a $\partial^{\kappa''_{\alpha}} A^{j''_{\alpha}} \, \langle j''_{\alpha}|$-propagator on the left side, and a $\partial^{\kappa'_{\alpha'}}B^{j'_{\alpha'}} \, 
|j'_{\alpha'}\rangle$-propagator on the right side,
they produce, as proved in Lemma \ref{lem:7}, a  propagator $\partial^{\kappa''_{\alpha}+\kappa'_{\alpha'}} \tilde{G}_{eff}^{j,j'}(z''_{\alpha},z'_{\alpha'})$ having a priori {\em three} scales -- $j,j'$ and
$\lfloor \log_2(t''_{\alpha}-t'_{\alpha'})\rfloor$ -- which may be resummed into an effective propagator $\partial^{\kappa''_{\alpha}+\kappa'_{\alpha'}}\tilde{G}_{eff}(z''_{\alpha},z'_{\alpha'})$. Thus in the sequel these are evaluated as a constant
$O(1)$, times a  
contraction\\ $\partial^{\kappa''_{\alpha}} \tilde{A}^{\tilde{j}''_{\alpha}}
\, \langle\tilde{j}''_{\alpha}|\  \partial^{\kappa'_{\alpha'}}\tilde{B}^{\tilde{j}'_{\alpha'}}\, |\tilde{j}'_{\alpha'}\rangle\ (z''_{\alpha},z'_{\alpha'})$,
where $\tilde{A}=A_{\nu^{(0)}+O(\lambda^2)}$, $\tilde{B}=B_{\nu^{(0)}+O(\lambda^2)}$,  and $\tilde{j}''_{\alpha}=\tilde{j}'_{\alpha'}= \log_2(t''_{\alpha}-t'_{\alpha'})+O(1)$.
 
\item[{\bf B.} ] {\bf (deterministic bound for the sum (\ref{eq:FN-k}))}
 The (deterministic) product of the $V_{\alpha}$'s    is compensated by the product 
 $\prod_{\alpha} r''_{\alpha}$, leaving only a small
coefficient $C^{-1}$
per vertex. Thus the sum  ${\cal S}_1(\vec{\gamma})$ (\ref{eq:vertex-insertions}) converges to a constant
$1+O(C^{-2})$. For $C$ small enough this is in the complex disk $B(1,1/2)$, so
$ {\bf 1}_{\Del_{\alpha}\in \D^0_{LF,k_{\alpha}}} \ \cdot\ 
 \log({\cal S}_1(\vec{\gamma})$ is well-defined, and 
\BEQ {\bf 1}_{\Del_{\alpha}\in \D^0_{LF,k_{\alpha}}} \ \cdot \log {\cal S
}_1(\vec{\gamma}) \simeq {\bf 1}_{\Del_{\alpha}\in \D^0_{LF,k_{\alpha}}}
 \ \cdot \left( {\cal S}_1(\vec{\gamma})-1\right).
\EEQ

\medskip\noindent
We must now sum the scaling coefficient $\prod_{\alpha}r_{\alpha}^{-1}$ over all vertex
locations, i.e. over all
forests $\F$. Since the $(r''_{\alpha})^{-1}$'s give (up to a constant $O(1)$ per vertex) the correct order of magnitude of the vertex insertions, we  may assume that we want
to sum over all large-field indices $\vec{k}$ (see {\bf C.}), then over all forests $\F$ (see {\bf D.})  the string  ${\cal S}_1-1$, see (\ref{eq:vertex-insertions}) , where one
has set: $\gamma_{\alpha}=O(1)$ and $R_{\eta}(\tau)(\cdot,\cdot)=\del(\cdot,\cdot)+G^{2\to}(\cdot,\cdot)$.

\item[{\bf C.}] {\bf (convergence of the average in $\eta$)}  The main issue here is to
show, using standard Gaussian large deviations, that our estimates are integrable in
$\eta$.  
Proceeding as in \S \ref{subsection:two-point} {\bf C.},
we rewrite $A^j((t,x),(t',x'))$ as $e^{-\frac{c}{2} |x-x'|^2/2^j} \ \cdot\ \tilde{A}^j((t,x),(t',x'))$, where $\tilde{A}^j$ has the same scaling properties as $A^j$, and
has furthermore retained the same Gaussian type space-decay, only with different 
constants; and similarly for $B^j$, $\tilde{B}^j$. Up to a multiplicative constant
$O(1)$, this is equivalent to replacing $\nu^{(0)}$ by $\tilde{\nu}^{(0)}\approx \nu^{(0)}$.  In the process, we have gained
a small factor $\prod_{\alpha} 2^{-cn(\Del_{\alpha})^{1+2/d}}$, see (\ref{eq:local-factorials}). Then we split in two  the large-deviation factor $LF(\F,\vec{k})=\prod_{\alpha} \proba[\Del_{\alpha}\in \D^0_{LF,k}]\lesssim \left( \prod_{\alpha} e^{-\frac{c}{2} 2^{2k}/\lambda} \right)^2$.
Then
\BEQ  \Big[ LF(\F,\vec{k}) \Big]^{1/2}  \ \cdot\ \Big[\prod_{\alpha} (r'_{\alpha}(k_{\alpha}))^{-1}  2^{-cn(\Del_{\alpha})^{1+2/d}} \Big]=O(1).\EEQ
This is easily shown using the space-decay, resp. large-deviation factor when $k_{\alpha}\le n(\Del_{\alpha})$, resp. $\ge n(\Del_{\alpha})$.
The remaining factor $\left[\LF(\F,\vec{k})\right]^{1/2}$ makes the sum
over large-field indices converge to a factor $O(1)$ per vertex, $\sum_{k_{\alpha}\in \{-\infty\}\cup\N}
e^{-\frac{c}{2} 2^{2k_{\alpha}}/\lambda}=1+o(1)$. 
Thus (see (\ref{eq:FN-k}))  
$\Big| \sum_{\F} \sum_{\vec{k}} \big\langle F_1(\F,\vec{k}|\eta) \big\rangle
\Big|\lesssim \sum^{\star}_{\F} \prod_{\alpha} (\tilde{r}''_{\alpha})^{-1}
\equiv \left(\sum_{\F}\prod_{\alpha} (\tilde{r}''_{\alpha})^{-1} \right)-1,
$
where: $\sum^*_{\F} f(\F):=\sum_{\F\not=\emptyset} f(\F)$, and 
$(\tilde{r}''_{\alpha})^{-1}$ is given by the same formula as (\ref{eq:r''alpha}), but
with $V_{\alpha}(z_{\alpha})$ (see (\ref{eq:Valpha})) replaced by 
 $\tilde{V}_{\alpha}(z_{\alpha}):= g^{(0)} \partial^{\kappa'_{\alpha}} \tilde{B}^{j'_{\alpha}}(\cdot,z_{\alpha})\partial^{\kappa''_{\alpha}}\tilde{A}^{j''_{\alpha}}(z_{\alpha},\cdot)$.

\item[{\bf D.}] {\bf (convergence of the sum over forests)} Following the same technique as in \S \ref{subsection:two-point}, we
integrate space-time variables in chronological order, yielding for $n$ vertices
\BEA 
&& \int d\bar{z}_{\alpha_1} \cdots d\bar{z}_{\alpha_n} \sum_{j'_{\alpha_1},\ldots,
j'_{\alpha_n}\ge 0} \sum_{j''_{\alpha_1},\ldots,j''_{\alpha_n}\ge 0} \ \int dz''\,  A((t_1,x_1),z'') \ \cdot\ \nonumber\\
&&  \cdot\ \int dz_{\alpha_1} \partial^{\kappa'_{\alpha_1}} \tilde{B}^{j'_{\alpha_1}}
(z'',z_{\alpha_1})  |j'_{\alpha_1}\rangle  g^{(0)} \underline{G}^{2\to}(z_{\alpha_1},\bar{z}_{\alpha_1}) \int dz''_{\alpha_1} \partial^{\kappa''_{\alpha_1}}\tilde{A}^{j''_{\alpha_1}}(\bar{z}_{\alpha_1},z''_{\alpha_1})  \langle j''_{\alpha_1}| \cdot\nonumber\\
&& \cdot\  \int dz_{\alpha_2}\, \partial^{\kappa'_{\alpha_2}} \tilde{B}^{j'_{\alpha_2}}
(z''_{\alpha_1},z_{\alpha_2})  |j'_{\alpha_2}\rangle  g^{(0)}\underline{G}^{2\to}(z_{\alpha_2},\bar{z}_{\alpha_2})  \int dz''_{\alpha_2}\, \partial^{\kappa''_{\alpha_2}} \tilde{A}^{j''_{\alpha_2}}
 (z_{\alpha_2},z''_{\alpha_2}) \langle j''_{\alpha_2}| \cdots \nonumber\\
 && \cdot\ \,
 \int dz_{\alpha_n}\,  \partial^{\kappa'_{\alpha_n}} \tilde{B}^{j'_{\alpha_n}}
(z''_{\alpha_{n-1}},z_{\alpha_n}) |j'_{\alpha_n}\rangle    g^{(0)} \underline{G}^{2\to}(z_{\alpha_n},
\bar{z}_{\alpha_n}) \int dz''_{\alpha_n}\partial^{\kappa''_{\alpha_n}} \tilde{A}^{j''_{\alpha_n}}
 (z_{\alpha_n},z''_{\alpha_n}) \langle j''_{\alpha_n}|  \cdot\nonumber \\
&& \cdot \int dy_1\,  B(z''_{\alpha_n},(0,y_1))\,  e^{\frac{\lambda}{\nu^{(0)}} h_0(y_1)}\EEA
where $\underline{G}^{2\to}(\cdot,\cdot):=\del(\cdot,\cdot)+G^{2\to}(\cdot,\cdot)$
as in (\ref{eq:v0-bd}).

Since, for $t\gtrsim 1$,  $G^{2\to}e^{\tilde{\nu}^{(0)}t\Del}(\cdot,\cdot)\lesssim \int_0^{O(1)} dt'\, 
e^{\tilde{\nu}^{(0)}(t'+ct)\Del}(\cdot\,\cdot)\lesssim e^{c'\tilde{\nu}^{(0)}t\Del}(\cdot,\cdot)$, 
the contribution of the  integration in the $\bar{z}$'s  may
be absorbed into the coupling constant $g^{(0)}$ and a redefinition of 
$\tilde{\nu}^{(0)}$. Thus we may assume that
$\bar{z}_{\alpha_i}\equiv z_{\alpha_i}$. Furthermore, since $(|j\rangle)_{j\ge 0}$ is
an orthonormal basis, $j''_{\alpha_i}=j'_{\alpha_{i+1}}$, and 
$\partial^{\kappa''_{\alpha_i}}\tilde{A}^{j''_{\alpha_i}}({z}_{\alpha_i},\cdot)\langle j''_{\alpha_i}|\ \cdot\  \partial^{\kappa'_{\alpha_{i+1}}}
\tilde{B}^{j'_{\alpha_{i+1}}}(\cdot,z_{\alpha_{i+1}}) |j'_{\alpha_{i+1}}\rangle =
\partial^{\kappa''_{\alpha_i}+\kappa'_{\alpha_{i+1}}} \tilde{G}^{j''_{\alpha_i}}(z_{\alpha_i},z_{\alpha_{i+1}})$.  

\medskip\noindent  We let $z_{\alpha_i}\equiv(t_{i+1},x_{i+1})$; rewrite the derivatives
 $\partial^{\kappa''_{\alpha_i}+\kappa'_{\alpha_{i+1}}}$ as $\partial^{\kappa_{i+1}}$,
 with $|\kappa_{i+1}|=0,3$ or $6$, 
which produces an equivalent factor 
$O((1+t_{i+1}-t_{i+2})^{-|\kappa_{i+1}|/2})$; bound 
$\sum_{j\ge 0}\tilde{G}^j((t,x),(t',x'))$ by $O(1) p_{\frac{\nu^{(0)}}{c}(t-t')}(x-x')$; and use for a sequence of two {\em low-momentum} $\tilde{G}$-propagators our {\em first power-counting estimate}
(compare with  (\ref{eq:PW1})),
\BEA &&  \int_{t_{i+2}}^{t_i} dt_{i+1}\, (1+t_i-t_{i+1})^{-|\kappa_i|/2} (p_{\frac{\nu^{(0)}}{c}(t_i-t_{i+1})}\ast p_{\frac{\nu^{(0)}}{c}(t_{i+1}-t_{i+2})})(x_i-x_{i+2}) \nonumber\\
&&\qquad \le
 O(1)\, p_{\frac{\nu^{(0)}}{c}(t_i-t_{i+2})}(x_i-x_{i+2}), \label{eq:PW1bis}
 \EEA
an estimate similar to but more precise than (\ref{eq:PW1}), valid
 for $|\kappa_i|>2$.
If $|\kappa_i|=3$, resp. $6$, then we apply (\ref{eq:PW1bis}) with $|\kappa_i|$ replaced by $2^+$, keeping $(1+t_{i}-t_{i+1})^{-(\frac{1}{2})^-}$, resp. $(1+t_{i}-t_{i+1})^{-2^-}$ in
store.  If $\kappa_i=0$ then $\kappa_{i+1}\not=0$ and $t_{i+2}-t_{i+1}\gtrsim
t_{i+1}-t_i$; we obtain similarly a factor $O(1)$ and keep in store  $(1+t_{i+1}-t_{i+2})^{-(\frac{1}{2})^-}$, resp. $(1+t_{i+1}-t_{i+2})^{-2^-}$. Extra factors 
$(1+t_{i}-t_{i+1})^{-3/2}$ are used to iterate, so that there remains in store exactly
$\prod_i (1+t_{i}-t_{i+1})^{-(\half)^-}$, where the product ranges over low-momentum
propagators. Each scale $0$ propagator $\tilde{G}^0((t_i,\cdot),(t_{i+1},\cdot))$
$(t_i-t_{i+1}\le 1)$, on the other hand, has $\kappa_i=0$, but benefits from a
small factor $O(g^{(0)})$ which can be rewritten in the form $(1+t_i-t_{i+1})^{-(\half)^-} O(g^{(0)})$. 

\medskip\noindent The conclusion is the following.  Rescale the coordinates, $(t_1,x_1)\rightsquigarrow (\eps^{-1}t_1,\eps^{-1/2}x_1)$. 
The main term in (\ref{eq:PW1bis}) is 
\BEA && \int dy_1\,  \tilde{G}_{eff}((\eps^{-1}t_1,\eps^{-1/2} x_1),(0,y_1))
e^{\frac{\lambda}{\nu^{(0)}} h_0}(y_1) \nonumber\\
&& = 1+e^{\nu_{eff}\eps^{-1} t_1\Del} (e^{\frac{\lambda}{\nu^{(0)}} h_0}-1)(\eps^{-1/2}x_1) + O(\eps) e^{(\nu^{(0)}+O(\lambda^2))\eps^{-1}t_1\Del} (e^{\frac{\lambda}{\nu^{(0)}} h_0}-1)(\eps^{-1/2}x_1) \nonumber\\
&=& 1+O(\lambda e^{\frac{\lambda}{\nu^{(0)}} ||h_0||_{\infty}}) \, \eps^{d/2} ||h_0||_{L^1}
\EEA
 ($n=0$), while terms with $n\ge 1$ are bounded by $O(\lambda \eps^{d/2})$ times
a prefactor 
\BEQ \prod_{i=1}^n \big[ O(\lambda)(1+t_{i}-t_{i+1})^{-(\frac{1}{2})^-} \big] \le
O(\lambda) t_1^{- (\frac{1}{2})^-}.
\label{eq:eps1/2}
\EEQ 
yielding after rescaling an error term $O(\lambda\eps^{(\half)^-})$.
   Hence we simply get 

\BEQ  \langle h(\eps^{-1}t_1,\eps^{-1/2}x_1)\rangle \equiv \frac{\nu^{(0)}}{\lambda} \langle \log(w(\eps^{-1}t_1,\eps^{-1/2}x_1)\rangle=O( e^{\frac{\lambda}{\nu^{(0)}} ||h_0||_{\infty}}) \, \eps^{d/2} ||h_0||_{L^1}. \EEQ

\end{itemize}

%%%%%%%%%%%%%%%%%%%%%%%

\subsection{KPZ truncated 2-point function} \label{subsection:hh}

%%%%%%%%%%%%%%%%%%%%%%%%%%%

We are now interested in the large scale behavior of the connected $2$-point
function (i.e. covariance function),
\BEQ \langle h(t_1,x_1)h(t_2,x_2)\rangle_{c}=\left(\frac{\nu^{(0)}}{\lambda}\right)^2 \langle \log(w(t_1,x_1))\log(w(t_2,x_2))\rangle_{c}=: \left(\frac{\nu^{(0)}}{\lambda}\right)^2
F_{2,c}((t_1,x_1),(t_2,x_2)). \EEQ

A simple way to generate the {\em connected} two-point function is to consider
two independent replicas $\eta_1,\eta_2$ of $\eta$; then 
\BEQ \langle h(t_1,x_1)h(t_2,x_2)\rangle_c= \half
\Big\langle \Big(h(t_1,x_1|\eta_1)-h(t_1,x_1|\eta_2)\Big) \Big( h(t_2,x_2|\eta_1)-
h(t_2,x_2|\eta_2) \Big) \Big\rangle,\EEQ
where $\langle \, \cdot\, \rangle$ now refers to the expectation with respect to
the pair $(\eta_1,\eta_2)$.
We make a cluster expansion as above in the propagators and in the covariance kernels of
$\eta_1$ and $\eta_2$, and get an expression similar to (\ref{eq:FN}). 
By symmetry $(1\leftrightarrow 2)$, there is  {\em at least one four-leg vertex}, which means that there is  (at least) one pairing $\langle \eta_p(z_{\beta_1})\eta_p(z_{\beta_2})\rangle$ $(p=1,2)$, $z_{\beta_i}=(t_{\beta_i},x_{\beta_i})$
$(i=1,2)$, coming
from a  $V_{\beta_1}(z_{\beta_1})$ insertion on the 1st string, and a $V_{\beta_2}(z_{\beta_2})$ insertion
on the 2nd string; the pairing vanishes unless  $z_{\beta_1},z_{\beta_2}$ are in the
same   box
$\Del\in\D^0$ or in neighboring boxes. 

\medskip\noindent {\em Choose} among the existing such $V_{\beta_i}(z_{\beta_i})$,
$i=1,2$ {\em the earliest one anti-chronologically}, i.e. the one with the
largest $t_{\beta_i}$. Because of the finite-range nature of the kernel $\langle
\eta(\cdot)\eta(\cdot)\rangle$, there exists a pairing $\langle \eta_i(z_{\beta_1})
\eta_i(z_{\beta'_2})\rangle$ with $d(\Del_{\beta_2},\Del_{\beta'_2})=O(1)$, and
similarly a pairing $\eta_i(z_{\beta'_1})\eta_i(z_{\beta_2})\rangle$ with 
$d(\Del_{\beta_1},\Del_{\beta'_1})=O(1)$; it may of course happen that $\beta_1=\beta'_1$, $\beta_2=\beta'_2$. Call ${\cal F}^0_{\Del_1,\Del_2}$ the set of forests
such that $z_{\beta_i}$, $i=1,2$ belong to fixed boxes $\Del_{\beta_1}:=\Del_1$,
$\Del_{\beta_2}:=\Del_2$; ${\cal F}^0_{\Del_1,\Del_2}$ is empty unless
$d(\Del_1,\Del_2)=O(1)$. Applying explicitly the operator $\frac{d}{d\gamma_{\beta_1}}
\frac{d}{d\gamma_{\beta_2}}$ to the r.-h.s. of (\ref{eq:FN}), one gets

\BEA &&  F_{2,c}((t_1,x_1),(t_2,x_2))=\sum_{\Del_1,\Del_2\in\D^0} \sum_{\F\in {\cal F}^0_{\Del_1,\Del_2}}
 \prod_{\alpha} \frac{1}{2\II\pi} \oint_{\partial B(0,r_{\alpha})} \frac{d\gamma_{\alpha}}{\gamma_{\alpha}^2}   \nonumber\\
&&\qquad \Big{\langle}  {\cal S}_1(\vec{\gamma})^{-1} {\cal S}_2
(\vec{\gamma})^{-1} \ \cdot\ \frac{d}{d\gamma_{\beta_1}} \left({\cal S}_1(\vec{\gamma})\right)\  \frac{d}{d\gamma_{\beta_2}} \left({\cal S}_2(\vec{\gamma}) \right)  \Big{\rangle}
\EEA

\medskip\noindent

%%%%%%%%%%%%%%%%%%%%%%%%%%ש
\vskip 1cm

\bigskip\noindent Compared to the previous subsection, we must now add a supplementary
estimate in the preparatory phase.

\bigskip\noindent {\bf A'.} {\bf (power-counting factors for $\eta$-pairings between strings})
As mentioned in \S 3.2 and  illustrated in \S \ref{subsection:four-point}, $\eta$-pairings produce {\em outer contractions} linking
different strings. Contrary to {\em inner contractions} inside $0$-th scale clusters which contribute to the
two-point function, outer contractions produce $4$-point functions, which
have not been renormalized. We must now show that the power-counting effect of an outer contraction is comparable to that described in (\ref{eq:small-factor}).
For that, consider parallel chronological sequences on two strings,
\BEA && \int_{\Del_1^0} dz'_1\,  \int_{\Del_{2}^0} dz'_{2}\,   \langle \eta(z'_1)\eta(z'_{2})\rangle_{\vec{s}(\vec{w})} \nonumber\\
&&\qquad  \left( \cdots A^{j_1}(\cdot,\cdot)\langle j_1|\ B^{j_1}(\cdot,z'_1)
|j_1\rangle\, A^{j'_1}(z'_1,\cdot) \langle j'_1|\  \cdots \right)
 \nonumber\\
&&\qquad  \left( \cdots A^{j_{2}}(\cdot,\cdot) \langle j_{2}|\  B^{j_{2}}(\cdot,z'_{2}) |j_{2}\rangle\ A^{j'_{2}}(z'_{2},\cdot) \langle j'_2|\,  \cdots \right)
\EEA
in which  vertex integration points $z'_1,z'_{2}$ are located in neighboring
boxes so that the average $\langle \eta(z'_1)
\eta(z'_{2})\rangle_{\vec{s}(\vec{w})}$ does not vanish.
Ladder diagrams considered in \S \ref{subsection:four-point}, see (\ref{eq:ladder}), ,
are of this type. Following the chronological integration procedure of {\bf D.}, we replace (supposedly already integrated) outgoing legs
$A^{j'_1},A^{j'_{2}}$ by $1$ and integrate over $z'_1,z'_{2}$. 
Since  $d(\Del^0_1,\Del^0_2)=O(1)$,  we may
just as well assume (up to a volume prefactor $O(1)$) that $\Del^0_1=\Del^0_{2}$.
We are free to choose the ordering of the strings and may therefore suppose that 
$j_1\le j_{2}$. 
Thanks to the exponential decay of $B^{j_1},B^{j_{2}}$, the 
space-time integration $\int dz'_1\, \int dz'_{2}$ costs a volume factor
$O(2^{j_1(1+\frac{d}{2})})$. On the other hand, 
were $z'_1,z'_{2}$ {\em not} constrained to be located in the same scale $0$
box, we would get instead a volume factor $O\left( \prod_{i=1}^{2} (2^{j_i(1+\frac{d}{2})}) \right)$. The overall
gain is therefore bounded up to a constant by
\BEQ 2^{-j_2(1+\frac{d}{2})} \le 
 \prod_{i=1}^{2} 2^{-5j_i/4} 
\label{eq:PW2} \EEQ
if $d\ge 3$, which is our {\em second key power-counting estimate}.  {\em This shows
that we have produced a small factor $O(2^{-\frac{5}{4}j})$ per low-momentum field
$G^j$, or equivalently $O((1+t-t')^{-5/4})$ per low-momentum field $G((t,x),(t',x'))$,
$t-t'\gtrsim 1$. }

\medskip\noindent {\em Remark.} In the case  $d=3$, this upper bound is optimal (for $j_2=j_1+O(1)$), and
not quite as good as the $O(2^{-\frac{3}{2}j})$ factor due to renormalization, compare
with (\ref{eq:small-factor}). However,  one 
easily shows that
the resulting small factor is actually comparable or smaller than (\ref{eq:small-factor}) if $d\ge 4$. In any case, in order to be able to integrate (see {\bf D.}) we
simply need a small factor $O(2^{-(1+2\kappa)j})$ per low-momentum field $G^j$,
with $\kappa>0$. In the KPZ$_2$ case ($d=2$), one finds $\kappa=0$; this
border case is no more super-renormalizable in the infra-red:  four-point 
functions are superficially divergent in the QFT terminology, which leads to a
floating (i.e. scale-dependent) coupling constant $g$. 

%%%%%%%%%%%%%%%%%%%%%%%%%%%%%%%%%%%%
\vskip 1cm
\bigskip\noindent
So much for {\bf A'.}  Resuming now our previous discussion,
and proceeding as in \S 5.2,  one can prove in exactly the same way that $F_{2,c}=O(1)$. The only difference is that (compare with the discussion below 
 (\ref{eq:PW1bis})), using our second power-counting estimate leaves in store
 in the worst case only
 $(1+t_i-t_{i+1})^{-(\frac{1}{4})^-}$ per low-momentum $G$.
 
\medskip\noindent There remains to see how one gets the prefactor $O(\eps^{\frac{d}{2}-1})$ and the scaling function $K_{eff}$. For that, we remark, proceeding
as in {\bf D.}, that (see (\ref{eq:A-low})) 
$ \frac{d}{d\gamma_{\beta_1}} {\cal S}_i(\vec{\gamma})$, $i=1,2$ is  equal to 
\BEQ 
\Big\{ \tilde{G}_{eff}((t_i,x_i),(t'_i,x'_i)) + O((t_i-t'_i)^{-(\half)^-}) \, 
G_{\nu^{(0)}/c}((t_i,x_i),(t'_i,x'_i))\Big\} \ \cdot\ R_{\eta}^{(0)}(\tau^0=0) \gamma_{\beta_1}V_{\beta_i}(\cdot)(
\cdot,\cdot)\cdots  \label{eq:5.28}
\EEQ 
where $t'_i\in\Del_i\equiv[t^+_{\Del_i}-1,t^+_{\Del_i})\times\bar{\Del}_i$. Then the
main term of 
${\cal S}_i(\vec{\gamma})^{-1}$ is 
\BEQ \int dy_i\, \tilde{G}_{eff}((t_i,x_i),(0,y_i))e^{\frac{\lambda}{\nu^{(0)}}h_0}(y_i)=1+O(\lambda e^{\frac{\lambda}{\nu^{(0)}}||h_0||_{\infty}}) t_1^{-d/2} ||h_0||_{L^1}. \EEQ
Error terms take into account: vertex insertions along any of the two strings,
costing either the already accounted for $O((t_i-t'_i)^{-1/2})$ or
$O((t'_i)^{-1/2})$  for the two numerators, and $O(t_i^{-1/2})$ for the
two denominators;  $\eta$-pairings between strings (see {\bf A'.}),
by construction at times $\le t^+_{\Del_1}+O(1)=t^+_{\Del_2}+O(1)$, costing
$O(((t^+_{\Del_1})^{-(\frac{1}{4})^-})^2)=O((t^+_{\Del_1})^{-(\half)^-})$; corrections in $O(t_i^{-d/2})$ or $O((t'_i)^{-d/2})$
due to the initial condition.

 \medskip\noindent Concluding: replacing the sum over boxes $\Del_i=\Del_{\beta_i}$,
 $\Del_{\beta'_i}$,  and the integral over $z_{\beta_i},z_{\beta'_i}$, $i=1,2$, by
 $O(1)$ times a single integral over a single  space-time variable $(t,x)$ located
 at distance $O(1)$ of all of these, and rescaling
the coordinates, we get asymptotically in the limit $\eps\to 0$ if $t_1\ge t_2$
\BEA &&   F_{2,c}((\eps^{-1}t_1,\eps^{-1/2}x_1),(\eps^{-1}t_2,\eps^{-1/2}x_2))
\equiv \frac{(g^{(0)})^2}{D^{(0)}} \langle h(\eps^{-1}t_1,\eps^{-1/2}x_1) h(\eps^{-1}
t_2,\eps^{-1/2}x_2)\rangle_c\nonumber\\
\label{eq:5.22a}\\
&&\sim_{\eps\to 0} F(\lambda)\   (g^{(0)})^2 \ \cdot\  \eps^{-1-d/2} \int_{0}^{t_2} dt \int dx\  
    \cdot\  \eps^{d/2}p_{\nu_{eff}(t_1-t)}(x_1-x)
 \cdot\,  \eps^{d/2}p_{\nu_{eff}(t_2-t)}(x_2-x) \nonumber\\
 &&\sim_{\eps\to 0} F(\lambda) \frac{(g^{(0)})^2}{D^{(0)}} \eps^{\frac{d}{2}-1}  \ \langle h(t_1,x_1)
h(t_2,x_2)\rangle_{0;\nu_{eff},D^{(0)}}  \label{eq:5.22b}
\EEA
up to error terms smaller by a factor  $O(\eps^{(\half)^-})$, for some function $F(\lambda)=1+O(\lambda^2)$ 
independent of the coordinates $(t_1,x_1),(t_2,x_2)$. Letting 
\BEQ D_{eff}:=F(\lambda)D^{(0)},\EEQ 
and
comparing (\ref{eq:5.22a}) with (\ref{eq:5.22b}), one sees that
\BEQ  \langle h(\eps^{-1}t_1,\eps^{-1/2}x_1) h(\eps^{-1}
t_2,\eps^{-1/2}x_2)\rangle_c \sim_{\eps\to 0} \eps^{\frac{d}{2}-1} \langle
 h(t_1,x_1)
h(t_2,x_2)\rangle_{\lambda=0;\nu_{eff},D_{eff}}, \EEQ
where the coefficient $D_{eff}=D^{(0)}(1+O(\lambda^2))$ is interpreted as the {\em effective noise
strength}.

\medskip\noindent The leading term  for $D_{eff}-D^{(0)}$ may be computed as follows. Following
the expansion in the number of vertices made in \S 4.1, see in particular (\ref{eq:leading-nu-eff}),  the
main term in\\ $\langle h(\eps^{-1}t_1,\eps^{-1/2}x_1) h(\eps^{-1}
t_2,\eps^{-1/2}x_2)\rangle_c$ is obtained from (\ref{eq:5.28}) by simply contracting 
a vertex $V_{\beta_1}\equiv B((t_1,x_1),(t'_1,x'_1)) g^{(0)}(\eta(t'_1,x'_1)-v^{(0)})A((t'_1,x'_1),\cdot)$ on the first string with a vertex
$V_{\beta_2}$ on the second string. Next comes the leading-order correction, obtained
by double-contracting $n=2$ vertex contributions on each string, yielding as in 
(\ref{eq:leading-nu-eff0})
\BEA &&  \Big\langle \Big[ (g^{(0)})^2 \eta(t'_1,x'_1) A((t'_1,x'_1),\cdot)B(\cdot,(t''_1,x''_1))
\eta(t''_1,x''_1) \Big]  \ \cdot \nonumber\\
&&\Big[ (g^{(0)})^2 \eta(t'_2,x'_2) A((t'_2,x'_2),\cdot)B(\cdot,(t''_2,x''_2))
\eta(t''_2,x''_2) \Big]  \Big\rangle
\EEA
Displacing the four outer $B-$ and $A-$ propagators $B(\cdot,(t'_1,x'_1))$, 
$A((t''_1,x''_1),\cdot)$, $B(\cdot,(t'_2,x'_2))$, $A((t''_2,x''_2),\cdot)$ to the
same point $(t'_1,x'_1)$,  integrating over $(t''_1,x''_1),(t'_2,x'_2),(t''_2,x''_2)$
and taking the limit $t'_1\to +\infty$
yields an effective contribution
\BEA &&  C_4:= (g^{(0)})^4 \lim_{t'_1\to +\infty} \int_0^{t'_1} dt''_1 \int dx''_1 \ \int_0^{+\infty} dt'_2
\int dx'_2\ \int_0^{t'_2} dt''_2\int dx''_2 \nonumber\\
&&\qquad (\omega\ast\omega)(t'_1-t'_2,x'_1-x'_2) \, G(t'_1-t''_1,x'_1-x''_1)\, G(t'_2-t''_2,x'_2-x''_2)\, (\omega\ast\omega)(t''_1-t''_2,x''_1-x''_2) \nonumber\\
\label{eq:leading-D-eff} \EEA 
Neglected terms involving e.g. $A((t'_1,x'_1),\cdot)-A((t''_1,x''_1),\cdot)$ involve
a low-momentum gradient, whence an extra $O(\eps^{1/2})$ which vanishes in the 
scaling limit.
Then $C_4$ is added to the main term which (after displacing outer $B$- and $A$-propagators)
becomes $C_2:=(g^{(0)})^2 \int_0^{+\infty} dt'_2 \int dx'_2\,  (\omega\ast\omega)(t'_1-t'_2,x'_1-x'_2)$. Thus $\frac{D_{eff}}{D^{(0)}}-1$ is given to leading
order by the quotient $C_4/C_2=O((g^{(0)})^2)=O(\lambda^2)$.

%%%%%%%%%%%%%%%%%%%%%%%

\subsection{Higher-order KPZ truncated functions}

%%%%%%%%%%%%%%%%%%%%

We must still prove that higher-order truncated functions
\BEQ \langle h(t_1,x_1)\cdots h(t_N,x_N)\rangle_c=:(\frac{\nu^{(0)}}{\lambda})^N
F_{N,c}((t_1,x_1),\ldots,(t_N,x_N)),
\EEQ
 $(N>2)$ are negligible in the
large scale limit because the KPZ field is asymptotically Gaussian, with
correlations given by $K_{eff}=F_{2,c}$. To be specific we prove this
for $N=4$, but the reader may easily adapt the following arguments to arbitrary $N$.
Let $F_{4,c}((t_i,x_i)_{i\le 4}):=\langle \log(w(t_1,x_1))\cdots \log(w(t_4,x_4))\rangle_c$.  The "replica trick" of \S \ref{subsection:hh} extends, with now
$4$ replicas of $\eta$,
\BEQ F_{4,c}:=\frac{1}{4} \Big\langle \prod_{\ell=1}^4 \sum_{k=0}^3
e^{\II k\pi/2} \log w(t_{\ell},x_{\ell}|\eta_{k+1}) \Big\rangle, \EEQ
a classical formula immediately generalized to arbitrary $N$ as
\BEQ F_{N,c}:=\frac{1}{N} \Big\langle \prod_{\ell=1}^N \sum_{k=0}^{N-1}
e^{2\II k\pi/N} \log w(t_{\ell},x_{\ell}|\eta_{k+1}) \Big\rangle, \EEQ
 originally proved
by P. Cartier\footnote{J. Lascoux, private communication.}.
 Then the connected function $F_{4,c}$ is obtained by selecting in (\ref{eq:FN}) those
contributions for which there is a permutation $\sigma$ of the index set
 $\{1,\ldots,4\}$, and for each $i=1,2,3$,  paired vertex insertions $V_{\beta_i}$, $V_{\beta'_i}$ on  strings number $\sigma(i),\sigma(i+1)$. Proceeding as in \S 5.4, we obtain a $O(1)$ denominator of 
 order $3\times 2=6$, multiplied by an expression bounded by (after
 coordinate rescaling) $(\eps^{\frac{d}{2}-1})^3$ instead of the expected overall
 scaling 
\BEA && \sum_{\mbox{pairings}\  \sigma} K_{eff}((\eps^{-1} t_{\sigma(1)},\eps^{-1/2} x_{\sigma(1)}),
 (\eps^{-1} t_{\sigma(2)},\eps^{-1/2} x_{\sigma(2)})) \ \cdot\ \nonumber\\
&&\qquad \cdot\ 
 K_{eff}((\eps^{-1} t_{\sigma(3)},\eps^{-1/2} x_{\sigma(3)}),(\eps^{-1} t_{\sigma(4)},\eps^{-1/2} x_{\sigma(4)}))=O(
 (\eps^{\frac{d}{2}-1})^2)
 \EEA
  for a four-point function.
 
%%%%%%%%%%%%%%%%%%%%%%%

\subsection{A remark on lower large-deviations for $h$}

%%%%%%%%%%%%%%%%%%%%%%%%%%%שש

Similar computations can be made for $\langle w^{-N}(t,x)\rangle$, where
$N=1,2,\ldots$, $N=O(1)$. Compared with the previous subsections, we now get
a product $\prod_{i=1}^N \Big({\cal S}_i(\vec{\gamma})\Big)^{-1}$
instead of $\prod_{i=1}^N \log({\cal S}_i(\vec{\gamma})$. It is 
easy to see that we get in the end
\BEQ \langle w^{-N}(t,x)\rangle=O(1).\EEQ
Using Markov's inequality e.g. for $N=1$ implies then for $A>0$
\BEQ \proba[h(t,x)<-A]=\proba[w^{-1}(t,x)>e^{\frac{\lambda}{\nu^{(0)}}A}]=O(1)\  e^{-\frac{\lambda}{\nu^{(0)}}A},\EEQ
an exponential lower large-deviation estimate for $h(t,x)$.

\medskip\noindent This is however disappointing with respect to the expected lower Gaussian large-deviation
\BEQ \proba[h(t,x)<-A]\lesssim e^{-cA^2}, \EEQ
proved using Gaussian concentration inequalities in Carmona-Hu\cite{CarHu}, Theorem 1.5 in
a deterministic setting. It is plausible that their results extend to our setting
by generalizing to regularized white noise  classical large deviation results for Lipschitz
functions of vector-valued Gaussian random variables, see e.g. \cite{AneBla},
\S 7.3.

%%%%%%%%%%%%%%%%%%%%%%%%%%%%%
%%%%%%%%%%%%%%%%%%%%%%%%%%%%%%%%%%%%

\section{Appendix. Cluster expansions}

%%%%%%%%%%%%%%%%%%%%%%%%%%%%
%%%%%%%%%%%%%%%%%%%%%%%%%%

%%%%%%%%%%%%%%%%%%%%%%%%%%

\subsection{Horizontal cluster expansion}

%%%%%%%%%%%%%%%%%%%%%%%%%
%%%%%%%%%%%%%%%%%%%%%%

The {\em cluster expansion} between boxes of scale 0 is performed according to
the classical Bridges-Kennedy-Abdesselam-Rivasseau (BKAR) procedure (see
\cite{AbdRiv1}, \cite{Unt-mode}, or \cite{MagUnt2}, \S 2.1 and 2.2), which
we now briefly describe, following \cite{MagUnt2}.   We apply it to the $A^0$ and $B^0$ kernels,  and also  to the covariance kernel $C_{\eta}(\cdot,\cdot):=
\langle \eta(\cdot)\eta(\cdot)\rangle$ of the noise. The
 effect of the cluster expansion on the $A$'s and $B$'s is to "cut" all propagators
 between scale $0$ boxes belonging to different polymers. The
effect of the cluster expansion on the $\eta$'s is to make independent the 
$\eta$-fields produced in scale $0$ boxes belonging to different polymers. As a result of those two operations, different polymers have been made totally independent, 
which makes it possible to extract averaged quantities such as counterterms. Since
the covariance kernel of $\eta$ has finite range (with our cut-off conventions 
$\langle \eta(t,x)\eta(t',x')\rangle=0$ except if $(t,x),(t',x')$ 
belong to the same unit box in $\Del^0$ or to neighboring boxes), the cluster expansion on the $\eta$'s is hardly noticeable -- in particular when it comes to bounds --, yet necessary.

\medskip\noindent Let ${\cal O}\subset\D^0$, and $|{\cal O}|:=
\cup_{\Del\in {\cal O}} \Del \subset\R_+\times\R^d$ its support. We say that two boxes
$\Del,\Del'\in{\cal O}$, $\Del\not=\Del'$, are {\em linked} if (i) either 
$\Del=[k,k+1)\times\bar{\Del}$, $\Del'=[k,k+1)\times\bar{\Del}'$, $\bar{\Del},\bar{\Del}'\in \bar{\D}^0$,  or (ii) 
$\Del=[k,k+1)\times \bar{\Del}$, $\Del'=[k-1,k)\times \bar{\Del}'$ or conversely $\Del=[k-1,k)\times \bar{\Del}$, $\Del'=[k,k+1)\times \bar{\Del}'$. By construction,
 there
exists $(t,x)\in\Del$, $(t',x')\in\Del'$, such that $A^0((t,x),(t',x'))\rangle\not=0$ or $A^0((t',x'),(t,x))\not=0$  if and only if $\Del=\Del'$ or
$\Del,\Del'$ are linked. Similarly, if $\langle \eta(t,x)\eta(t',x')\rangle\not=0$, then
$\Del=\Del'$ or $\Del,\Del'$ are linked (and, furthermore, $d(\Del,\Del')=O(1)$).   Denote by
$L({\cal O})$ the set of linked pairs $\{\Del,\Del'\}$.  Then, for every
{\em link weakening} of ${\cal O}$, i.e. for every function
$\vec{s} \ : \ L({\cal O})\to [0,1]$, extended trivially on the diagonal by
letting $\vec{s}_{\Del,\Del}\equiv 1$ ($\Del\in{\cal O}$), we define
\BEQ B^0(\vec{s})((t,x),(t',x'))=s_{\Del^0_{t,x},
\Del^0_{t',x'}} B^0((t,x),(t',x')), \EEQ
\BEQ   A^{0}(\vec{s}) ((t,x),(t',x'))=s_{\Del^0_{t,x},
\Del^0_{t',x'}} A^0((t,x),(t',x')) \EEQ
\BEQ \langle \eta(t,x)\eta(t',x')\rangle_{\vec{s}}:=s_{\Del^0_{t,x},
\Del^0_{t',x'}}  \langle \eta(t,x)\eta(t',x')\rangle \EEQ
if $(\Del_{t,x},\Del_{t',x'})\in L({\cal O})$, $0$ else. Thus the effect of the
function $\vec{s}$ is to {\em weaken off-diagonal elements of the propagator/covariance
kernel}.

\medskip\noindent
We need some terminology before we get to the point. In the 
following discussion, ${\cal O}$ is fixed. A scale $0$ {\em forest} $\F^0$
is a finite number of boxes $\Del\in{\cal O}$, seen as {\em vertices},
connected by  {\em links}, without loops. A (non-oriented) link $\ell$ connects 
$\Del_{\ell}$ to $\Del'_{\ell}$. Space-time variables ranging in $\Del_{\ell}$, resp. 
$\Del'_{\ell}$ are generally denoted as $(t_{\ell},x_{\ell})$, resp. $(t'_{\ell},x'_{\ell})$, or for short $z_{\ell}$, resp. $z'_{\ell}$.  {\em Non-isolated components} of $\F^0$, i.e. connected components of $\F^0$ containing $\ge 2$ boxes are called {\em trees}, or (specifically in this statistical physics
context) {\em polymers}. The (finite) set  of vertices of polymers
is denoted by $V(\F^0)$.  The set of all $0$-th scale cluster forests  is denoted by ${\cal F}^0({\cal O})$, or simply ${\cal F}^0$ if ${\cal O}=\D^0$.
If there exists a link between $\Del$ and $\Del'$, then we write
$\Del\sim_{\F^0}\Del'$, or simply (if no ambiguity may arise)
$\Del\sim\Del'$.

\medskip

Now  the following formula --
called BKAR formula -- holds: let $F=F(A^0,B^0|\eta)$ be some random function of 
the $A^0$ 's and $B^0$ 's, then

\begin{Proposition}[BKAR formula](see \cite{MagUnt2}, Proposition 2.6)  \label{prop:BKAR}

\BEQ  \langle F(A^0,B^0|\eta) \rangle =\sum_{\F^0\in{\cal F}^0}  \left(
\prod_{\ell\in L(\F^0)}\int_0^1   dw_{\ell}\right) \left(  \left(\prod_{\ell\in L(\F^0)} \frac{d}{ds_{\ell}}\right) 
\langle F(A^0(\vec{s}(\vec{w})),B^0(\vec{s}(w)))|\eta \rangle_{\vec{s}(\vec{w})} \right) \label{eq:BKAR}  \EEQ
$s_{\Del,\Del'}(\vec{w})$,
$\Del\not=\Del'$  being the infimum of
the $w_{\ell}$ for $\ell$ running over the unique path from $\Del$ to $\Del'$ in $\F^0$
if $\Del\sim_{\F^0}\Del'$, and  $s_{\Del,\Del'}(\vec{w})=0$ else.

\end{Proposition}
 The above formula is obtained by iterating the following step-by-step procedure.
Choose some box $\Del_1\in\D^0$, and Taylor-expand simultaneously with respect
to the parameters $(s_{\ell})_{\ell}$ where $\ell$ ranges in the set $L_1(\D^0)$ of all pairs
$\{\Del_{\ell},\Del'_{\ell}\}$ such that $\Del_1=\Del_{\ell}$ or $\Del'_{\ell}$. One
obtains: 
\BEA && F(\vec{s}\big|_{L(\D^0)\setminus L_1(\D^0)};\vec{s}\big|_{L_1(\D^0)}=1)=F(\vec{s}\big|_{L(\D^0)\setminus L_1(\D^0)};\vec{s}\big|_{L_1(\D^0)}=0)
\nonumber\\
&&\qquad  + \sum_{\ell_1\in L_1({\cal O})} \int_0^1 dw_1\,  \partial_{s_{\ell_1}}F(\vec{s}\big|_{L(\D^0)\setminus L_1(\D^0)};\vec{s}\big|_{L_1(\D^0)}=w_1). 
\EEA

\medskip\noindent  The following elementary relation is shown in \cite{Unt-mode},  
\BEQ \frac{d}{ds_{\ell}}\langle F(\eta)\rangle_{\vec{s}(\vec{w})} = 
\int_{\Del_{\ell}} dz_{\ell} \int_{\Del'_{\ell}} dz'_{\ell} 
\, \langle \eta(z_{\ell})\eta(z'_{\ell})\rangle_{\vec{s}=1} \ \cdot\ 
\Big{\langle}  \frac{\del}{\del\eta(z_{\ell})} \frac{\del}{\del \eta(z'_{\ell})}
F(\eta) \Big{\rangle}_{\vec{s}(\vec{w})}.  \label{eq:d/ds-eta}
\EEQ
In other words,  an $s$-derivative  acting on
an averaged quantity $\langle F(\eta)\rangle_{\vec{s}(\vec{w})}$ has the effect of  producing an explicit
pairing $\langle \eta(z_{\ell})\eta(z'_{\ell})\rangle_{\vec{s}=1}$, with the
original covariance kernel,  between
two arbitrary points belonging resp. to one box and to the other box.

\medskip\noindent As explained before,
each choice of forest $\F^0$ yields an {\em explicit  connection} through
$A^0$,  $B^0$- or $\eta$-pairings of all boxes within a given connected component (tree), 
and {\em disconnects} boxes $\Del,\Del'$ lying in different connected
components since $B^{0}(\vec{s}(w)) ((t,x),(t',x'))=A^{0}(\vec{s}(w)) ((t,x),(t',x'))=
\langle \eta(t,x)\eta(t',x')\rangle_{\vec{s}(\vec{w})}=0$ for
$(t,x)\in\Del, (t',x')\in\Del'$.

%%%%%%%%%%%%%%%%%%%%%%%%

\subsection{Mayer expansion}

%%%%%%%%%%%%%%%

\bigskip\noindent For the Mayer expansion (see \S \ref{subsection:two-point}), we
choose another set of objects $\cal O$ and a different way of implementing the $s$-dependence, and apply a slightly different formula. Namely, we let ${\cal O}\equiv
{\cal O}(\F^0, \{\mu_{\Del}\}_{\Del})$ be the
set of scale $0$ {\em polymers}, i.e. of {\em non-isolated connected components of $\F^0$ with their
external structure}, depending on the differentiation orders $\{\mu_{\Del}\}_{\Del}$,  produced by the vertical cluster expansion. Among these polymers, there are polymers with exactly 
{\em two} external legs, making up a subset ${\cal O}_1\equiv {\cal O}_1(\F^0,\{\mu_{\Del}\}_{\Del})\subset{\cal O}$. The complementary set ${\cal O}_2\equiv {\cal O}_2(\F^0,\{\mu_{\Del}\}_{\Del}):={\cal O}\setminus{\cal O}_1$ is made up of polymers
with $>2$ external legs, which require no renormalization. 
The following variant of BKAR's formula, found originally in \cite{AbdRiv1}, is
 stated in the present form in \cite{MagUnt2}. We now denote by $\{\P_{\ell},\P'_{\ell}\}$ a pair of polymers connected by a link $\ell\in L({\cal O})$.

 \begin{Proposition}[restricted 2-type cluster or BKAR2 formula] \label{prop:BK2}

Assume ${\cal O}={\cal O}_1\amalg {\cal O}_2$. Choose as initial object an object
$o_1\in {\cal O}_1$ of type 1, and stop the Brydges-Kennedy-Abdesselam-Rivasseau expansion as soon as a
link to an object of type 2 has appeared. Then choose a new object of type $1$, and so
 on. This  leads to a  restricted expansion, for which {\em only} the  link variables
  $z_{\ell}$, with $\ell\not\in {\cal O}_2\times {\cal O}_2$, have been weakened. The
   following closed formula holds. Let $\vec{S}:L({\cal O})\to [0,1]$ be a link
   weakening of ${\cal O}$, and $F=F((\vec{S}_{\ell})_{\ell\in L({\cal O})})$ 
 a smooth function. Let ${\cal F}_{res}({\cal O})$ be the set of  forests
    $\G^0$ on $\cal O$, each component of which is (i) either a tree of
objects of type $1$, called {\em unrooted tree}; (ii)or a {\em rooted tree} such that only the root is of type $2$. Then
\BEQ F(1,\ldots,1)=\sum_{\G^0\in {\cal F}_{res}({\cal O})}
\left( \prod_{\ell \in L(\G^0)}
\int_0^1 dW_{\ell} \right) \left( \left(\prod_{\ell\in L(\G^0)}  
\frac{\partial}{\partial
S_{\ell} } \right) F(S_{\ell}(\vec{W})) \right),  \label{eq:Mayer-BKAR}
 \EEQ

where $S_{\ell}(\vec{W})$ is either $0$ or the minimum of the $w$-variables running along the
unique path in $\bar{\G}^0$ from $\P_{\ell}$ to $\P'_{\ell}$, and $\bar{\G}^0$ is the
forest obtained from $\G^0$ by merging all roots of $\G^0$ into a single vertex.

\end{Proposition}

\medskip\noindent The way functions of the type $\langle F(A^0(\vec{s}(\vec{w})),B^0(\vec{s}(\vec{w}))|\eta)\rangle_{\vec{s}(\vec{w})}$ are made
$S$-dependent is explained in \ref{subsection:two-point}. Differentiating w.r. to
an $S$-parameter $S_{\P_1,\P_2}$  produces a factor\\ $\Big[ \prod_{\Del_1\in\P_1,
\Del_2\in\P_2,(\Del_1,\Del_2)\not\in
{\bf\Del}_{ext}(\P_1)\times
{\bf\Del}_{ext}(\P_2)} {\bf 1}_{\Del_1\not=\Del_2} \Big]-1$, which upon expansion yields a sum over all overlap possibilities
between boxes of $\P_1$ and boxes of $\P_2$ except those containing the external
legs. Each contribution comes with a sign $(-1)^n$, where $n$ is the number of
overlapping boxes (see Fig. 3 for a representation of this rule).
See also Fig. 4 below
illustrating a more elaborate case with $n=2$.

\begin{figure}[H]  \label{fig:Mayer1}
  \centering
  \caption{Mayer subtraction rule for one overlapping box.}
  \vskip -3cm
   \includegraphics[scale=0.4]{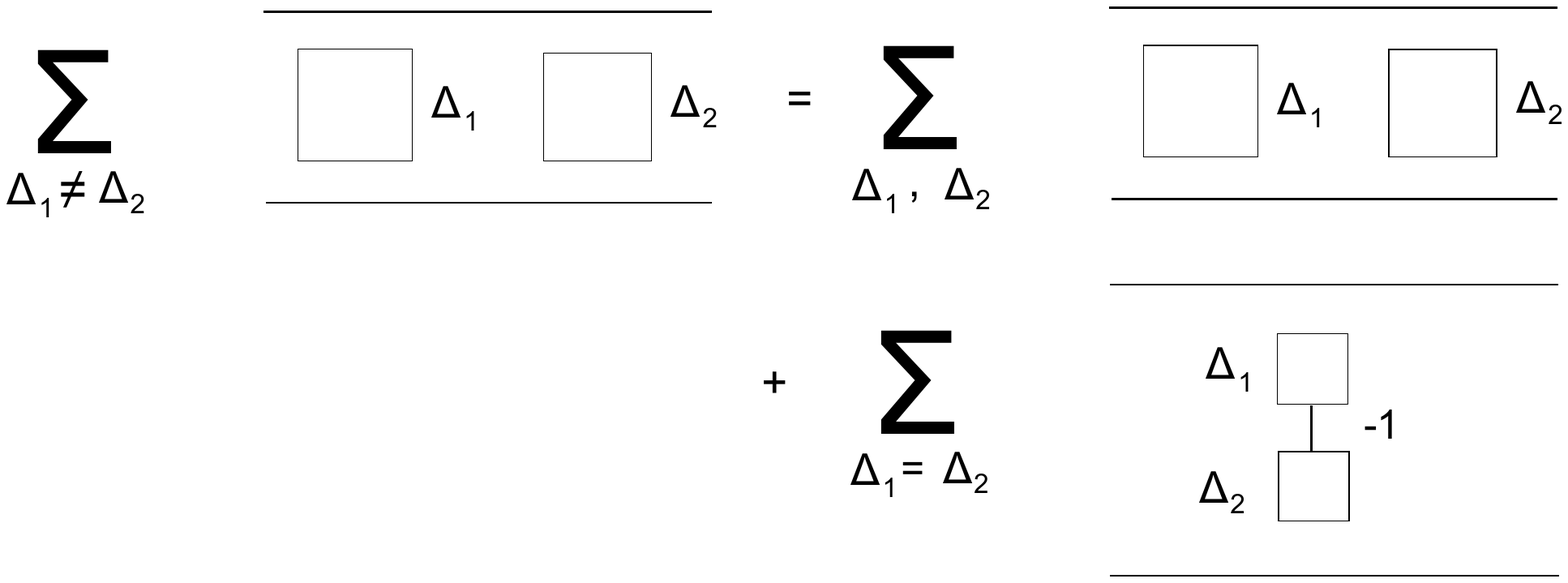}
\end{figure}

\begin{figure}[H]  \label{fig:Mayer}
  \centering
  \caption{Mayer expansion.}
  \vskip -3cm
   \includegraphics[scale=0.9]{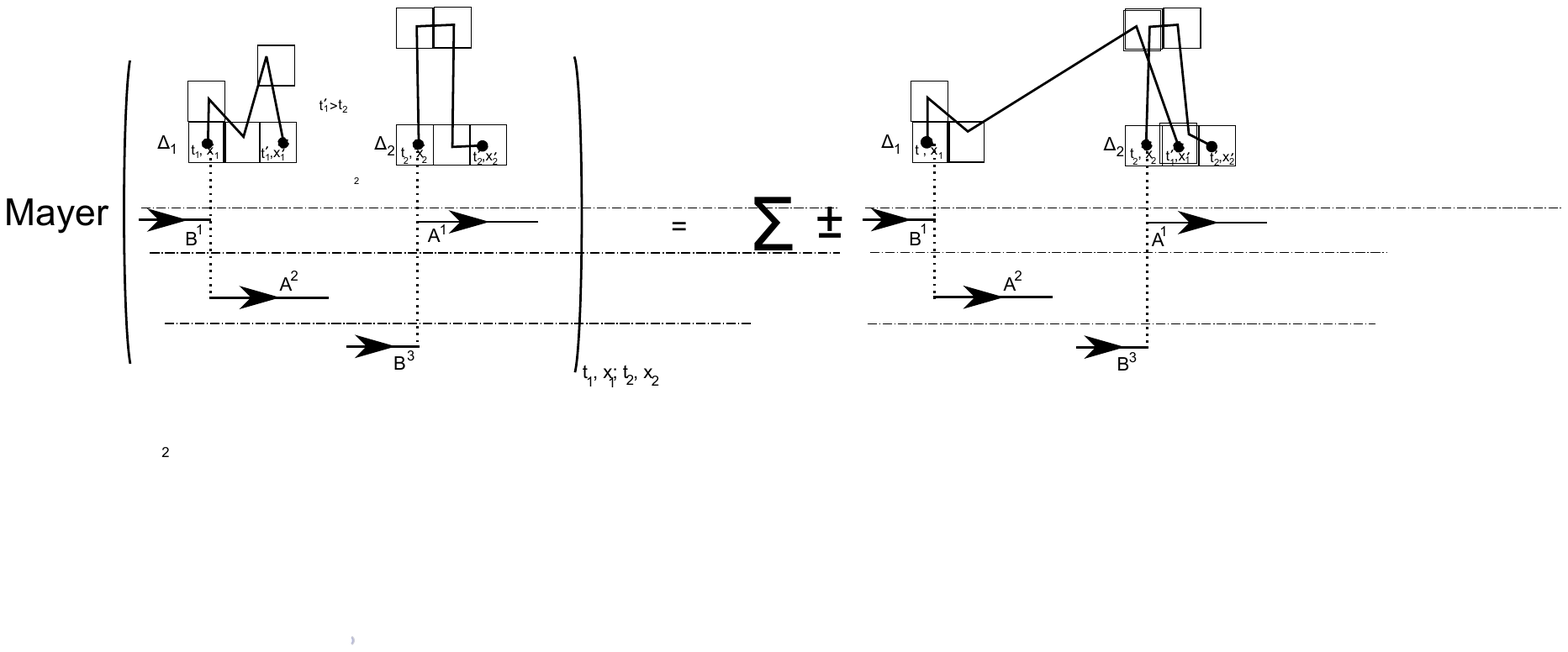}
\end{figure}

 The above procedure leads, as discussed in a much more involved,
multi-scale
context e.g. in \cite{MagUnt2}, Proposition 2.12,   to some mild combinatorial factors, which
we discuss briefly. Recall that (by Cayley's theorem) the number of trees over $\P_1,\ldots,\P_n$ with fixed
coordination numbers $(n(\P_i))_{i=1,\ldots,n}$ equals $\frac{n!}{\prod_i (n(\P_i)-1)!}$. Choose a tree $\T$ component of $\G^0$. Start from the leaves of $\T$ and go down the branches
inductively. Let $\P_1,\ldots,\P_{n(\P')-1}$ be the leaves attached onto one and the
same vertex $\P'$. Choose $n(\P')-1$ (possibly non distinct) boxes $\Del_1,\ldots,
\Del_{n(\P')-1}\in\D^0$ of $\P'$ (there
are $|\P'|^{n(\P')-1}$ possibilities), and assume that $\Del_i\in\P_i$. For each choice of polymer $\P'$, this gives a supplementary factor $O((C  |\P'|)^{n(\P')-1})$, to be
multiplied by $\frac{1}{(n(\P')-1)!}$ coming from Cayley's theorem. Summing over
$n(\P')=2,3,\ldots$ yields $e^{C |\P'|}-1\le e^{C|\P'|}$.  Summing over all boxes takes care automatically of the sum over all permutations of the polymers, which takes down the $n!$ factor.  Since bounds produced in section \ref{sec:bounds} are in
$O((g^{(0)})^m)$, where $m=\sum_{\P\in{\cal O}} |\P|$ is the number of boxes obtained by the cluster expansion, the latter large factor is compensated by a simple
redefinition of coupling constant $g^{(0)}\rightsquigarrow e^C g^{(0)}=O(g^{(0)})$ in
the bounds.

%%%%%%%%%%%%%%%ש

%%%%%%%%%%%%%%%%%%%%%%%
%%%%%%%%%%%%%%%%%%ש

\section{Appendix. The effective propagator}

%%%%%%%%%%%%%%%%%%%%%%
%%%%%%%%%%%%%%%%%%%%%%%ש

The effective propagator $\tilde{G}_{eff}$ obtained in \S \ref{subsection:hh}  by
 resumming  $\nu$-counterterms along a string, see (\ref{eq:Gtilde-eff}) below,
 is shown in this section to be very well approximated at large scale by the Green kernel
 $(\partial_t-\nu_{eff}\Del)^{-1}$. 
 
\noindent We first need a technical lemma.

\begin{Lemma} \label{lem:technical}
There exists some constant $C>0$ such that, for every $\vec{\kappa}=(\kappa_1,\ldots,
\kappa_d)$ and $t>t'$, $x,x'$:
\BEQ \Big|\nabla^{\vec{\kappa}} G_{\nu^{(0)}}((t,x),(t',x'))\Big| \le 
C^{|\vec{\kappa}|+1}\, \Big(\lambda \sqrt{\nu^{(0)}(t-t')} \Big)^{-|\vec{\kappa}|} \, \Gamma(|\vec{\kappa}|/2))\,  G_{\nu^{(0)}-O(\lambda^2)}((t,x),(t',x')).
\EEQ
\end{Lemma}

\noindent{\bf Proof.} The spatial Fourier transform of $\nabla^{\vec{\kappa}} G\equiv \nabla^{\vec{\kappa}}G_{\nu^{(0)}}$ is
\BEQ  \widehat{\nabla^{\vec{\kappa}} G}(t-t',\vec{\xi})=(\II\vec{\xi})^{\vec{\kappa}} \, \hat{K}_{t-t'}(\vec{\xi}), \qquad 
\hat{K}_{t-t'}(\vec{\xi}):=e^{-(t-t')\nu^{(0)}(\vec{\xi},\vec{\xi})}.\EEQ
Let $\vec{\xi}_0:=\frac{x-x'}{2\nu^{(0)}(t-t')}$.  Then

\BEA && \nabla^{\vec{\kappa}}G(t-t',x-x')=(2\pi)^{-d} \int_{\R^d} d\vec{\xi}\, (\II\vec{\xi})^{\vec{\kappa}}
\, \hat{K}_{t-t'}(\vec{\xi}) e^{\II (x-x',\xi)} \nonumber\\
&&=(2\pi)^{-d} \int_{\R^d+\II \vec{\xi}_0}   d\vec{\xi}\, (\II\vec{\xi})^{\vec{\kappa}}
\, \hat{K}_{t-t'}(\vec{\xi}) e^{\II (x-x',\xi)} \nonumber\\
&&= (2\pi)^{-d} e^{-|x-x'|^2/4\nu^{(0)}(t-t')} \int_{\R^d} d\vec{\xi}\,  (\II\vec{\xi}- \vec{\xi}_0)^{\vec{\kappa}} \, 
e^{-\nu^{(0)}(t-t')|\vec{\xi}|^2} \nonumber\\
&&=(\nu^{(0)}(t-t'))^{-|\vec{\kappa}|/2} G_{\nu^{(0)}}((t,x),(t',x')) \int_{\R^d} d\vec{\zeta}\,  (\II\vec{\zeta}-\vec{\zeta}_0)^{\kappa}
\, e^{-|\vec{\zeta}|^2}
\EEA
where now $\vec{\zeta}:=\sqrt{\nu^{(0)}(t-t')}\, \vec{\xi}$, $\vec{\zeta}_0:=\frac{x-x'}{2\sqrt{\nu^{(0)}(t-t')}}$ are non-dimensional parameters.
Rewrite $G_{\nu^{(0)}}((t,x),(t',x'))$ as $G_{\nu^{(0)}+O(\lambda^2)}((t,x),(t',x')) \, \cdot\, O(1)\, e^{-O(\lambda^2)|\vec{\zeta}_0|^2}$. Then, for all
$\vec{\kappa}'\le\vec{\kappa}$,
$$|\vec{\zeta}_0|^{|\vec{\kappa'}|} e^{-\lambda^2|\vec{\zeta}_0|^2}\lesssim \frac{|\vec{\zeta}_0|^{\vec{\kappa'}}}{\lambda^{|\vec{\kappa}|'}|\vec{\zeta}_0|^{|\vec{\kappa}'|}/\Gamma(
\frac{|\vec{\kappa'}|}{2}+1)} =\lambda^{-|\vec{\kappa}'|} 
\Gamma(\frac{|\vec{\kappa'}|}{2}+1) $$
and $\int_{\R^d} d\vec{\zeta}\,  \vec{\zeta}^{\vec{\kappa}-\vec{\kappa}'} \, e^{-|\vec{\zeta}|^2}
= O(C^{|\vec{\kappa}|}) \Gamma(|\vec{\kappa}-\vec{\kappa}'|/2).$ One concludes by
using the binomial formula.
\hfill\eop

\vskip 1cm
\bigskip\noindent Let us now come to the point. Recall $\Del^{\to 0}=
\bar{\chi}^{(0)}\ast\Del$  (see (\ref{eq:Del->0})) is an ultra-violet regularization of $\Del$.

\begin{Lemma}  \label{lem:7}
Let $\del\nu:=\nu_{eff}-\nu^{(0)}$,
\BEQ G_{eff}:=(\partial_t-\nu_{eff}\Del)^{-1}, \EEQ
\BEQ \tilde{G}_{eff}:=A^{\to 1} \left(1-\del\nu B^{\to 1} \Del^{\to 0} A^{\to 1} \right)^{-1} B^{\to 1}=(1-\del\nu\,  G^{\to 1}\Del^{\to 0})^{-1} G^{\to 1}. \EEQ
 Then:
\begin{enumerate}
\item There exists $\tilde{\nu}^{(0)}=\nu^{(0)}+O(\lambda^2)$ and a constant $C>0$
such that
\BEQ \tilde{G}_{eff}((t,x),(t',x'))\le C G_{\tilde{\nu}^{(0)}}((t,x),(t',x')).\EEQ

Furthermore, if $1\le j\le j'\le j''$, 
\BEA && \left(\nabla^{\kappa'} A^j \langle j| \,  \left(1-\del\nu B^{\to 1} \Del^{\to 0} A^{\to 1} \right)^{-1} \nabla^{\kappa''} B^{j'} \, |j'\rangle
\right)
((t,x),(t',x')) \nonumber\\
&&\qquad\lesssim  2^{-\frac{j'}{2}(|\kappa'|+|\kappa''|)} 2^{-(j'-j)} G_{\tilde{\nu}^{(0)}/c}((t,x),(t',x')), \qquad
t-t'\approx 2^{j'}  \label{eq:jj'j''1} \EEA

\BEA &&  \left(\nabla^{\kappa'} A^j \langle j| \,  \left(1-\del\nu B^{\to 1} \Del^{\to 0} A^{\to 1} \right)^{-1} \nabla^{\kappa''} B^{j'} \, |j'\rangle
\right)
((t,x),(t',x')) \nonumber\\
&&\qquad\lesssim   2^{-\frac{j'}{2}(|\kappa'|+|\kappa''|)}2^{-(j''-j)} 2^{-(j''-j')} G_{\tilde{\nu}^{(0)}/c}((t,x),(t',x')), \qquad
t-t'\approx 2^{j''} \label{eq:jj'j''2}  \nonumber\\
\EEA
with $c=1$ if $\kappa'=\kappa''=0$, and $c=\half$ else.

 \item  For every $\kappa'<\half$, the following holds: if $t-t'\approx 1$,
\BEQ (\tilde{G}_{eff}-G_{eff})((\eps^{-1}t,\eps^{-1/2}x),(\eps^{-1}t',\eps^{-1/2}x')) 
\sim_{\eps\to 0} O(\eps^{2\kappa'}) G_{\tilde{\nu}^{(0)}}((\eps^{-1}t,\eps^{-1/2}x),(\eps^{-1}t',\eps^{-1/2}x'))  \label{eq:Gtilde-eff} \EEQ
uniformly for 
\BEQ |x-x'|^2=o(\eps^{-1/2}(t-t')) \label{eq:beyond-regime} \EEQ
 (see (\ref{eq:beyond-regime2}) below)  with $\tilde{\nu}^{(0)}=\nu^{(0)}+O(\lambda^2)$.  
 
\end{enumerate} 
\end{Lemma}

\medskip\noindent  Thus $\tilde{G}_{eff}$ is equal to $G_{eff}$ with an excellent approximation at large scale which holds well beyond the normal regime $\frac{|x|^2}{t}
\lesssim 1$ (one can compare with \cite{FS} where extended heat-kernel asymptotics are
 shown for a lattice regularization instead). Eq. (\ref{eq:jj'j''1},\ref{eq:jj'j''2})
show that the  bounds on  $\tilde{G}_{eff}=\sum_{j,j'\ge 1} A^{j} \, \langle j| \,
 \cdot\,  \left(1-\del\nu B^{\to 1} \Del^{\to 0} A^{\to 1} \right)^{-1} \, \cdot\, 
  B^{j'}\, |j'\rangle$, expressed as a product of a resolvent by two propagators
$A,B$  as in (\ref{A1-VB}), decrease exponentially with the difference of the "scales"
of these three operators, thus yielding bounds that can be resummed adequately. 

\medskip\noindent {\bf Proof.} We concentrate on 2., with 1. proved on the way. Introduce $G_{1,eff}:=(\partial_t-\nu^{(0)}\Del-
\del\nu \Del^{\to 0})^{-1}=A (1-\del\nu B \Del^{\to 0}A)^{-1} B$,  and write for short
$\chi^0$ instead of $\bar{\chi}^{(0)}$. 

\begin{itemize}
\item[(i)] First, {\em  $G_{1,eff}((\eps^{-1}t,\eps^{-1/2}x),(\eps^{-1}t',\eps^{-1/2}x'))=(1+O(\eps)) G_{eff}((\eps^{-1}t,\eps^{-1/2}x),(\eps^{-1}t',\eps^{-1/2}x'))$ if
(\ref{eq:beyond-regime}) holds.} Namely, the spatial Fourier transform of 
$G_{eff}-G_{1,eff}$ is
\BEQ \hat{G}_{eff}(\eps^{-1}(t-t'),\vec{\xi})-\hat{G}_{1,eff}(\eps^{-1}(t-t'),\vec{\xi})=
\hat{K}_{\eps^{-1}(t-t')}(\vec{\xi}),\EEQ
\BEQ \mbox{with}\,  \qquad \hat{K}_{t-t'}(\vec{\xi})= e^{-(t-t')\nu_{eff}(\vec{\xi},\vec{\xi})} 
\ \cdot\ \left(1-e^{-(t-t')\del\nu\,  (\widehat{\chi^0}(\xi)-1)\,  (\vec{\xi},\vec{\xi})} \right).\EEQ
Since $\chi^0$ is compactly supported, its Fourier transform $\widehat{\chi^0}$ extends to
an entire function satisfying: $|\widehat{\chi^0}(\xi)|\lesssim e^{C|\Im(\xi)|}$. If $\chi^0(\cdot)$ is chosen to be isotropic
(which we assume), then $\nabla(\widehat{\chi^0})(0)=0$. Since $\int \chi^0=1$ and $\chi^0$ is smooth, 
$|\widehat{\chi^0}(\vec{\xi})-1|=O_{\vec{\xi}\to 0}(|\vec{\xi}|^2)=
O_{\vec{\xi}\to 0}((t-t')|\vec{\xi}|^2)$ and $|\widehat{\chi^0}(\vec{\xi})|\, |\vec{\xi}|^2=
O_{|\Re(\vec{\xi})|\to \infty}(e^{C|\Im(\vec{\xi})|})=O((t-t')e^{C\sqrt{t-t'}\, |\Im(\vec{\xi})|}).$
{\em  Let 
\BEQ \vec{\xi}_0:= \frac{x}{2\nu_{eff}(t-t')}, \qquad \rho_0:=
\frac{1}{\sqrt{2\nu_{eff}(t-t')}}
(1+\frac{|x|}{\sqrt{2\nu_{eff}(t-t')}}).
\EEQ
}
 Note that, provided 
 \BEQ |x|^2=O(\eps^{-\half+2(\half-\kappa)}(t-t'))\qquad  \mbox{with}\  \frac{1}{4}<\kappa<\half
 \label{eq:beyond-regime2} \EEQ
  -- which is compatible with our hypothesis
(\ref{eq:beyond-regime}) if one lets $\kappa\to(\half)^-$ -- , and $|\vec{\xi}|\lesssim \eps^{\kappa}
\rho_0$ -- whence $|\xi|\ll 1$ -- the error term in the exponential, $\frac{t-t'}{\eps}\del\nu\, |1-\widehat{\chi^0}(\vec{\xi})|\, (\vec{\xi},\vec{\xi})=O(\lambda^2) \eps^{4\kappa-1} \rho_0^4$ is a $O(1)$.  
 Hence 
\BEA && \int_{B(0,\eps^{\kappa}\rho_0)} d\vec{\xi}\, \hat{K}_{\eps^{-1}(t-t')}(\vec{\xi}) e^{\II (x-x',\vec{\xi})}\nonumber\\
&&\  = 
\int_{B(0,\eps^{\kappa}\rho_0)+\II \eps^{1/2}\vec{\xi}_0} d\vec{\xi}\,  \hat{K}_{\eps^{-1}(t-t')}(\vec{\xi}) e^{\II(x-x',\vec{\xi})} + \partial I(x-x') \nonumber\\
&&\ = O(\frac{t-t'}{\eps}\del\nu) \ e^{-|x-x'|^2/4\nu_{eff}(t-t')}
 \int_{B(0,\eps^{\kappa}\rho_0)} d\vec{\xi}\, |\vec{\xi}|^4\, e^{-\eps^{-1}\nu_{eff}(t-t')(\vec{\xi},\vec{\xi})}
  + \partial I(x-x') \nonumber\\
&&\ \sim_{\eps\to 0} O(\lambda^2 \eps^{2\kappa'})\, 
 G_{eff}((\eps^{-1}t,\eps^{-1/2}x),(\eps^{-1}t',\eps^{-1/2}x')) + \partial I(x-x')
 \label{eq:7.15}
\EEA
 
where   $2\kappa'=-1-\frac{d}{2}+ \kappa(d+4) \to_{\kappa\to\half} 1$ and
\BEQ \partial I(x-x'):=\int_{\partial B(0,\eps^{\kappa}\rho_0)\times [0,\II\eps^{1/2}\vec{\xi}_0]} d\vec{\xi}\,  \hat{K}_{\eps^{-1}(t-t')}(\vec{\xi}) e^{\II(x-x',\vec{\xi})} \EEQ 

and
\BEA &&
|\partial I(x-x')|,\ \Big|\int_{|\vec{\xi}|\gtrsim\eps^{\kappa} \rho_0} d\vec{\xi}\, K_{\eps^{-1}(t-t')}(\vec{\xi}) e^{\II (x-x',\vec{\xi})} \Big| \nonumber\\
&& \qquad=O\Big(\eps^{d/2}
\int_{|\vec{\zeta}|\gtrsim \eps^{-(\half-\kappa)} (1+\frac{|x|}{2\nu_{eff}(t-t')})} d\vec{\zeta}\, e^{-\nu_{eff}|\vec{\zeta}|^2}  \Big) \nonumber\\
\EEA
are negligible with respect to $G_{eff}((\eps^{-1}t,\eps^{-1/2}x),(\eps^{-1}t',\eps^{-1/2}x')).$

\medskip\noindent On the other hand,
\BEA
&& A(1-\del\nu B\Del^{\to 0}A)^{-1} B^0|0\rangle = A^0\, \langle 0| \,  B^0 |0\rangle +\del\nu AB\Del^{\to 0}A^0 \, \langle 0|\,  B^0 |0\rangle + \cdots \nonumber\\
&&= G^0 + \del\nu \, G_{1,eff}\,  \Del^{\to 0}G^0.
\EEA

\item[(ii)] Next, $\left(1-\del\nu B^{\to 1} \Del^{\to 0} A^{\to 1} \right)^{-1}
\simeq 1$ at large scale. Namely, expanding $\left(1-\del\nu B^{\to 1} \Del^{\to 0} A^{\to 1} \right)^{-1}$ into a series, we get a geometric series in $\del\nu\, \Del^{\to 0} A^{\to 1}B^{\to 1}$. Write $A^{\to 1}B^{\to 1}\equiv G^{\to 1}$ as $\int_0^{+\infty} \tilde{\chi}^{\to 1}(t) e^{\nu^{(0)} t\Del} dt$, where (in the notations of Definition \ref{def:A}) $\tilde{\chi}^{\to 1}:=\sum_{j=1}^{+\infty} (\chi\ast\chi)^j$ is "one minus a bump function", i.e.
$\tilde{\chi}^{\to 1}\big|_{[0,c]}=0, \tilde{\chi}^{\to 1}\big|_{[c^{-1},+\infty)}=1$
for some $c>0$.   Then,
since $\Del^{\to 0}$ commutes with the $G^{\to 1}$'s, 
\BEA &&( \Del^{\to 0}G^{\to 1})^2((t,x),(t',x')) \nonumber\\
&&\qquad=\int_{t'}^{(t+t')/2} dt''\, \int dx''\, 
(\Del^{\to 0})^2 G^{\to 1}((t,x),(t'',x'')) G^{\to 1}((t'',x''),(t',x')) \nonumber\\
&&\qquad +
\int_{(t+t')/2}^{t} dt''\, \int dx''\, 
 G^{\to 1}((t,x),(t'',x'')) (\Del^{\to 0})^2 G^{\to 1}((t'',x''),(t',x'))
 \nonumber\\
&&\qquad =O(1)\  (t-t') (\Del^{\to 0})^2 G^{\to 1}((t,x),(t',x')).\EEA
We call this the {\em commutation trick}. Recall $|\del\nu|=O(\lambda^2)$.
Iterating yields by using Lemma \ref{lem:technical}
\BEA && \sum_{n\ge 1} (\del \nu)^n\  \Big| (\Del^{\to 0}G^{\to 1})^n ((t,x),(t',x'))  \Big|
\lesssim \sum_{n\ge 1} (\del\nu)^n \frac{(t-t')^{n-1}}{(n-1)!} \  \Big|(\Del^{\to 0})^n G^{\to 1}((t,x),(t',x')) \Big| \nonumber\\
&&\qquad \lesssim \sum_{n\ge 1}  \frac{(t-t')^{n-1}}{(n-1)!} 
2^{-n} \Gamma(n) (t-t')^{-n} \tilde{G}^{\to 1}((t,x),(t',x')) \nonumber\\
&&\qquad = O(1)  \ (t-t')^{-1}
\tilde{G}^{\to 1}((t,x),(t',x'))
\EEA
where $\tilde{G}^{\to 1}=G^{\to 1}_{\nu^{(0)}+O(\lambda^2)}$, and

\BEA  && \Big[\left(1-\del\nu B^{\to 1} \Del^{\to 0} A^{\to 1} \right)^{-1} -1\Big] ((t,x),(t',x'))  \nonumber\\
&&\qquad
=\del\nu   B^{\to 1} \Big( \sum_{n\ge 0} (\del \nu)^n  (\Del^{\to 0}G^{\to 1})^n
\Big) \Del^{\to 0} A^{\to 1} ((t,x),(t',x')) \nonumber\\
&&\qquad =O( \del\nu)\,  (t-t')^{-1}\,  \tilde{G}^{\to 1}((t,x),(t',x')
\EEA
where $\tilde{G}^{\to 1}$ has again been possibly rescaled. Using for a third time the
commutation trick, one finally gets
\BEQ \Big|\Big(A \Big[\left(1-\del\nu B^{\to 1} \Del^{\to 0} A^{\to 1} \right)^{-1} -1\Big] B\Big)
((t,x),(t',x')) \Big| =O(1)\  G^{\to 1}((t,x),(t',x')). \label{eq:7.22} \EEQ

\medskip\noindent On the other hand, the orthonormality of the basis $(|j\rangle)_{j\ge 0}$ implies immediately
\BEQ \Big(A^0 \langle 0|\,  \Big[\left(1-\del\nu B^{\to 1} \Del^{\to 0} A^{\to 1} \right)^{-1} -1\Big] B\Big)
((t,x),(t',x'))= G^0((t,x),(t',x')).\EEQ

\medskip\noindent Point 1. is a particularization of (\ref{eq:7.22}).  If 
$t-t'\approx 2^{j'}$ and $j\not=j'$, then 
\BEQ A^j \langle j| \left(1-\del\nu B^{\to 1} \Del^{\to 0} A^{\to 1} \right)^{-1} B^{j'}
|j'\rangle ((t,\cdot),(t',\cdot))=
\del\nu G^j \Del^{\to 0} G^{j'}((t,\cdot),(t',\cdot)+\cdots
\EEQ has an extra $2^{-(j'-j)}$-prefactor due to a reduced volume of integration in time. If $t-t'\approx 2^{j''}\gg 2^{j'}$, then
the leading term in the series vanishes, so that
\BEA &&  A^j \langle j| \left(1-\del\nu B^{\to 1} \Del^{\to 0} A^{\to 1} \right)^{-1} B^{j'}
|j'\rangle ((t,\cdot),(t',\cdot)) \nonumber\\
&&\qquad=\left((\del\nu)^2 G^j \Del^{\to 0}G^{\to 1}\Del^{\to 0}G^{j'}+\ldots\right)((t,\cdot),(t',\cdot)), 
\EEA
 where the middle propagator $G^{\to 1}$ has scale $ j''+O(1)$, leading for the same
 reason to an extra $2^{-(j''-j)}2^{-(j''-j')}$-prefactor. Gradients
 $\nabla^{\kappa'},\nabla^{\kappa''}$ are easily turned into prefactors by using
 elementary heat kernel estimates  as in Lemma \ref{lem:multi-scale-estimates-A-B} (i).

\item[(iii)] Let us now bound
\BEA  && D:=  A \Big[ (1-\del\nu B\Del^{\to 0}A)^{-1}-(1-\del\nu B^{\to 1}\Del^{\to 0}
A^{\to 1})^{-1} \Big] B \nonumber\\
&& = A \Big[ (1-\del\nu B\Del^{\to 0}A)^{-1} \ \cdot\ \del\nu \left(B\Del^{\to 0}A-B^{\to 1}
\Del^{\to 0} A^{\to 1}\right) (1-\del\nu B^{\to 1}\Del^{\to 0}
A^{\to 1})^{-1}  \Big] B\nonumber\\
&& = A \Big[ (1-\del\nu B\Del^{\to 0}A)^{-1}\ \cdot\ \del\nu \left( B^0|0\rangle\,  \Del^{\to 0} A +B^{\to 1} \Del^{\to 0} A^0 \langle 0| \, \right) (1-\del\nu B^{\to 1}\Del^{\to 0}
A^{\to 1})^{-1} \Big] B. \nonumber\\
\EEA

Thus
\BEA && D=\del\nu \Big(G^0 + \del\nu \, G_{1,eff}\,  \Del^{\to 0}G^0\Big) \Del^{\to 0}
\Big( A(1-\del\nu\, B^{\to 1}\Del^{\to 0}A^{\to 1})^{-1} B \Big) \nonumber\\
&& \qquad + \del\nu\, \tilde{G}_{1,eff} \Del^{\to 0}G^0
\EEA
where the kernel 
\BEA && \tilde{G}_{1,eff}(\cdot,\cdot)=A(1-\del\nu\, B\Del^{\to 0}A)^{-1}B^{\to 1}(\cdot,\cdot)=AB^{\to 1}(\cdot,\cdot)+
\del\nu AB\Del^{\to 0}AB^{\to 1}(\cdot,\cdot)+\ldots \nonumber\\
&&\qquad = G^{\to 1}(\cdot,\cdot)+\del\nu\,  G_{1,eff}\Del^{\to 0}AB^{\to 1}(\cdot,\cdot)\EEA
 is bounded  (using again and again the commutation trick) by
$O(1)\, G_{\nu^{(0)}+O(\lambda^2)}(\cdot,\cdot)$. Hence  

\BEA  &&|D((t,x),(t',x'))| \lesssim  (t-t')^{-1} G_{\nu^{(0)}+O(\lambda^2)}((t,x),(t',x')).  \label{eq:D} \EEA

\item[(iv)] Finally, 
\BEA && G_{1,eff}-\tilde{G}_{eff}-D\nonumber\\
&&\qquad=A \left(1-\del\nu B^{\to 1} \Del^{\to 0} A^{\to 1} \right)^{-1} B - A^{\to 1}  \left(1-\del\nu B^{\to 1} \Del^{\to 0} A^{\to 1} \right)^{-1} B^{\to 1} \nonumber\\
&&\qquad= A^0 \langle 0|\,   \left(1-\del\nu B^{\to 1} \Del^{\to 0} A^{\to 1} \right)^{-1} B +
A^{\to 1}  \left(1-\del\nu B^{\to 1} \Del^{\to 0} A^{\to 1} \right)^{-1} B^0 |0\rangle \nonumber\\
&&\qquad =G^0
\label{eq:GGD} \EEA
and $G^0((\eps^{-1}t,\eps^{-1/2}x),(\eps^{-1}t',\eps^{-1/2}x'))=0$ for $\eps$ small
enough.

\end{itemize} \hfill \eop

\noindent{\bf Remark.} Using a suitably   chosen cut-off $\chi^0$ with vanishing first momenta
(obtained e.g. by subtracting the beginning of the Taylor expansion of its Fourier transform
near zero), 
i.e. such that $\int dx\, x_{i_1}\cdots x_{i_p} \chi^0(x)=0$ for $1\le i_1,\ldots,i_p\le d$
and $p=2,3,\ldots,n-1$ one gets $\nabla^p(\widehat{\chi^0})(0)=0$, $2\le p\le n-1$, which
makes it possible to reduce the prefactor $O(\lambda^2 \eps^{1^-})$ in (\ref{eq:7.15})
to $O(\lambda^2)$ times an arbitrary large power of $\eps$.

%%%%%%%%%%%%%%%%%%%%%%%%

%%%%%%%%%%%%%%%%%%%%%%%%%%%%%%%%%%
%%%%%%%%%%%%%%%%%%%%%%%%%%%%%%%%\`u
%%%%%%%%%%%%%%%%%%%%%%\`u\`u
%%%%%%%%%%%%%%%%%%%%%%%
%%%%%%%%%%%%%%%%%%%%%%%%%
%%%%%%%%%%%%%%%%%%%%%%%
\end{document}